\newif\ifpdf
\ifx\pdfoutput\undefined
\pdffalse 
\else
\pdfoutput=1 
\pdftrue \fi

\newif\iffinal
\finalfalse	
\finaltrue 

\documentclass[reqno,twoside,11pt]{amsart}
\iffinal\else\usepackage[notref,notcite]{showkeys}\fi
\iffinal\else\IfFileExists{pdfsync.sty}{\usepackage{pdfsync}}{}\fi
\usepackage{cite}
\usepackage{amsmath}
\usepackage{amsfonts}
\usepackage{amssymb}
\usepackage{epsfig}
\usepackage{verbatim}

\usepackage{mathrsfs}
\usepackage{mathpazo}
\usepackage{mathabx}

\usepackage{amsthm,mathtools}
\usepackage{xcolor}
\usepackage{ esint }

\DeclareFontFamily{OT1}{eusb}{} \DeclareFontShape{OT1}{eusb}{m}{n} {<5> <6> <7> <8> <9> <10> <11> <12> <14.4> eusb10}{}
\DeclareMathAlphabet{\eusb}{OT1}{eusb}{m}{n}

\DeclareFontFamily{OT1}{eusm}{} \DeclareFontShape{OT1}{eusm}{m}{n} {<5> <6> <7> <8> <9> <10> <11> <12> <14.4> eusm10}{}
\DeclareMathAlphabet{\eusm}{OT1}{eusm}{m}{n}

\DeclareFontFamily{OT1}{eufm}{} \DeclareFontShape{OT1}{eufm}{m}{n} {<5> <6> <7> <8> <9> <10> <11> <12> <14.4> eufm10}{}
\DeclareMathAlphabet{\mathfrak}{OT1}{eufm}{m}{n}

\DeclareFontFamily{OT1}{fraktura}{}
\DeclareFontShape{OT1}{fraktura}{m}{n} {<5> <6> <7> <8> <9> <10> <11> <12> <13> <14.4> [1.1] eufm10}{}
\DeclareMathAlphabet{\fraktura}{OT1}{fraktura}{m}{n}

\DeclareFontFamily{OT1}{cmfi}{} \DeclareFontShape{OT1}{cmfi}{m}{n} {<5> <6> <7> <8> <9> <10> <11> <12> <13> <14.4> [0.9] cmfi10}{}
\DeclareMathAlphabet{\cmfi}{OT1}{cmfi}{b}{n}

\DeclareFontFamily{OT1}{cmss}{} \DeclareFontShape{OT1}{cmss}{m}{n} {<5> <6> <7> <8> <9> <10> <11> <12> <13> <14.4> cmss10}{}
\DeclareMathAlphabet{\cmss}{OT1}{cmss}{m}{n}

\newtheoremstyle{thm}{1.5ex}{1.5ex}{\itshape\rmfamily}{} {\bfseries\rmfamily}{}{2ex}{}

\newtheoremstyle{def}{1.5ex}{1.5ex}{\rmfamily\sl}{} {\bfseries\rmfamily}{}{2ex}{}

\newtheoremstyle{rem}{1.3ex}{1.3ex}{\rmfamily}{} {\itshape}
{} {1.5ex}{}

\newenvironment{proofsect}[1] {\vskip0.1cm\noindent{\rmfamily\itshape#1.}}{\qed\vspace{0.15cm}}

\theoremstyle{thm}
\newtheorem{theorem}{Theorem}[section]
\newtheorem{lemma}[theorem]{Lemma}
\newtheorem{proposition}[theorem]{Proposition}

\newtheorem*{Main Theorem}{Main Theorem.}
\newtheorem{corollary}[theorem]{Corollary}

\theoremstyle{def}
\newtheorem{definition}[theorem]{Definition}

\theoremstyle{rem}
\newtheorem{remark}[theorem]{{Remark}}

\numberwithin{equation}{section}

\renewcommand{\section}{\secdef\sct\sect}
\newcommand{\sct}[2][default]{\refstepcounter{section}
\addcontentsline{toc}{section}
{{\tocsection {}{\thesection}{\!\!\!\!#1\dotfill}}{}}
\vspace{0.7cm}
\centerline{ 
\scshape\arabic{section}.\ #1} \nopagebreak \vspace{0.2cm}}
\newcommand{\sect}[1]{
\vspace{0.4cm} \centerline{\large\scshape\rmfamily #1}
\vspace{0.2cm}}

\renewcommand{\subsection}{\secdef\subsct\sbsect}
\newcommand{\subsct}[2][default]{\refstepcounter{subsection}
\addcontentsline{toc}{subsection}
{{\tocsection{\!\!}{\hspace{1.2em}\thesubsection}{\!\!\!\!#1\dotfill}}{}}
\vspace{0.45\baselineskip} {\flushleft\bf
\arabic{section}.\arabic{subsection}~\bf #1.~}
\\*[3mm]\noindent
\nopagebreak}
\newcommand{\sbsect}[1]{\vspace{0.1cm}\noindent
\textbf{#1.~}\vspace{0.1cm}}

\renewcommand{\subsubsection}{%
\secdef \subsubsect\sbsbsect}
\newcommand{\subsubsect}[2][default]{%
\refstepcounter{subsubsection} 
\addcontentsline{toc}{subsubsection}{{\tocsection{\!\!}
{\hspace{3.05em}\thesubsubsection}{\!\!\!\!#1\dotfill}}{}}
\nopagebreak
\vspace{0.15\baselineskip} \nopagebreak {\flushleft\rmfamily
\itshape\arabic{section}.\arabic{subsection}.\arabic{subsubsection}
\ \rmfamily #1\/.}\ }
\newcommand{\sbsbsect}[1]{\vspace{0.1cm}\noindent
\rmfamily \itshape
\arabic{section}.\arabic{subsection}.\arabic{subsubsection} \
\sffamily #1\/.\ }

\iffinal

\else

\fi

\renewcommand{\caption}[1]{%
\vglue0.5cm
\refstepcounter{figure}
\begin{minipage}{0.9\textwidth}\small {\sc Figure~\thefigure. }#1\end{minipage}}

\newcommand{\supp}{\operatorname{supp}}

\newcommand{\textd}{\text{\rm d}\mkern0.5mu}
\newcommand{\texti}{\text{\rm  i}\mkern0.7mu}
\newcommand{\texte}{\text{\rm  e}\mkern0.7mu}

\newcommand{\Var}{\text{\rm Var}}
\newcommand{\Cov}{\text{\rm Cov}}

\renewcommand{\AA}{\mathcal A}
\newcommand{\BB}{\mathcal B}

\newcommand{\FF}{\mathcal F}

\newcommand{\LL}{\mathcal L}
\newcommand{\MM}{\mathcal M}
\newcommand{\NN}{\mathcal N}

\newcommand{\VV}{\mathcal V}

\newcommand{\ZZ}{\mathcal Z}

\newcommand{\E}{\mathbb E}

\newcommand{\N}{\mathbb N}

\newcommand{\BbbP}{\mathbb P}

\newcommand{\R}{\mathbb R}

\newcommand{\T}{\mathbb T}

\newcommand{\Z}{\mathbb Z}

\newcommand{\scrX}{\mathscr{X}}

\newcommand{\twoeqref}[2]{(\ref{#1}--\ref{#2})}

\def\myffrac#1#2 in #3{\raise 2.6pt\hbox{$#3 #1$}\mkern-1.5mu\raise 0.8pt\hbox{$#3/$}\mkern-1.1mu\lower 1.5pt\hbox{$#3 #2$}}

\newcommand{\cc}{\text{\rm c}}

\newcommand{\wt}{\widetilde}

\DeclareMathOperator{\GFF}{GFF}
\DeclareMathOperator{\DG}{DG}
\newcommand{\laweq}{\,\overset{\text{\rm law}}=\,}

\newcommand{\lawarrow}{{\overset{\text{\rm law}}\longrightarrow}}
\newcommand{\PPP}{\text{\rm PPP}}
\newcommand{\frakz}{\fraktura z}
\newcommand{\frakg}{\fraktura g}
\newcommand{\frakq}{\fraktura q}
\newcommand{\frakc}{\fraktura c}

\newcommand{\MB}{}
\newcommand{\eMB}{\normalcolor}

\usepackage{ulem}
\newcommand{\stkout}[1]{\ifmmode\text{\sout{\ensuremath{#1}}}\else\sout{#1}\fi}

\newcommand{\ZZZ}{\eusb Z}

\begin{document}

\vglue-3mm

\title[Maximum of DG model \hfill]{
A limit law for the maximum of subcritical\\DG model on a hierarchical lattice}
\author[\hfill M.~Biskup, H.~Huang]
{Marek~Biskup \,{\tiny and}\, Haiyu Huang}
\thanks{\hglue-4.5mm\fontsize{9.6}{9.6}\selectfont\copyright\,\textrm{2026}\ \ \textrm{M.~Biskup, H.~Huang. Reproduction for non-commercial use is permitted.}}
\maketitle

\vspace{-5mm}
\centerline{\it Department of Mathematics, UCLA, Los  Angeles, California, USA}
\smallskip


\begin{abstract}
We study the extremal properties of the ``integer-valued Gaussian'' a.k.a.\ DG model on the hierarchical lattice $\Lambda_n:=\{0,\dots,b-1\}^n$ (with $b\ge2$) of depth~$n$. This is a random field~$\varphi\in\mathbb Z^{\Lambda_n}$ with law proportional to $\texte^{\frac12\beta(\varphi,\Delta_n\varphi)}\prod_{x\in\Lambda_n}\#(\textd\varphi_x)$, where~$\Delta_n$ is the hierarchical Laplacian,~$\beta$ is the inverse temperature  and~$\#$ is the counting measure on~$\Z$.
Denoting $\beta_\cc:=2\pi^2/\log b$ and $m_n:=\beta^{-1/2}[(2\log b)^{1/2}n-\frac32(2\log b)^{-1/2}\log n]$, for~$0<\beta<\beta_\cc$ we prove that, along increasing sequences of~$n$  such that the fractional part of~$m_{n}$ converges to an~$s\in[0,1)$, the centered maximum $\max_{x\in\Lambda_n}\varphi_x-\lfloor m_n\rfloor$ tends in law to a discrete variant of a randomly shifted Gumbel law with the shift depending non-trivially on~$s$. The convergence extends to the extremal process whose law tends to a decorated Poisson point process with a random intensity measure. The proofs rely on renormalization group analysis which enables a tight coupling of the DG model to a  hierarchical  Gaussian Free Field.  The interval~$(0,\beta_\cc]$ marks the full range of values of~$\beta$ for which the renormalization group iterations tend to a ``trivial'' fixed point. 
\end{abstract}

%


\section{Introduction}
\label{sec1}\noindent
\vglue-3mm
\subsection{Background}
The extremal properties of logarithmically correlated random fields have been a subject of considerable interest in recent years. A picture that has emerged from the analysis of specific examples is that, in these systems, the suitably centered maximum tends in law to a randomly shifted Gumbel random variable while the associated extremal process tends to a decorated (Gumbel) Poisson point process with a random intensity measure. The randomness of the shift/intensity arises from a dependence structure that, in these systems, persists all the way up to the macroscopic scale.

The salient examples where such limit theorems have been proved include the Branching Brownian motion (Arguin, Bovier and Kistler~\cite{ABK1,ABK2,ABK3}, A\"idekon, Berestycki, Brunet and Shi~\cite{ABBS}), critical Branching Random Walks (A\"idekon~\cite{Aidekon}, Madaule~\cite{Madaule2}) and the Gaussian Free Field on~$\Z^2$ (Bramson, Ding and Zeitouni~\cite{BDingZ}, Biskup and Louidor~\cite{BL1,BL2,BL3}). Further evidence of universality has arrived in the studies of more general logarithmically correlated Gaussian processes (Madaule~\cite{Madaule1}, Ding, Roy and Zeitouni~\cite{DRZ}, Schweiger and Zeitouni~\cite{SZ}) including the four-dimensional membrane model (Schweiger~\cite{Schweiger}). Non-Gaussian processes have been treated as well, e.g., the characteristic polynomial of a random matrix ensemble (Paquette and Zeitouni~\cite{PZ}), the local time of simple random walk on a regular tree (Biskup and Louidor~\cite{BL4}, Abe and Biskup~\cite{AB}) and a class of  interacting fields  on a torus  (Bauerschmidt and Hofstetter~\cite{BH}, Hofstetter~\cite{H},  Barashkov, Gunaratnam and Hofstetter~\cite{BGH}).

The goal of the present paper is to address the limit law of the maximum for a process that distinguishes itself from the above by the field taking only integer values. Specifically, we are interested in the extremal properties of the model of a random interface that is  referred to as the ``integer-valued Gaussian field'' or the ``DG model'' (with DG standing for ``Discrete Gaussian'') in the literature. 

Generally, the \textit{DG model} is a $\Z$-valued process $\{\varphi_x\colon x\in\Lambda\}$ indexed by vertices of a finite graph~$\Lambda$ whose law takes the form
\begin{equation}
\label{E:1.1}
\frac1{\Sigma(\beta)}\texte^{\frac12\beta(\varphi,\Delta\varphi)}\prod_{x\in\Lambda}\#(\textd\varphi_x),
\end{equation}
where~$\Delta$ is the Laplacian associated with~$\Lambda$, the parameter~$\beta>0$ is the inverse temperature, $\#$ is the counting measure on~$\Z$ and~$\Sigma(\beta)$ is a normalization constant. The inner product $(\cdot,\cdot)$ is that in~$\ell^2(\Lambda)$. A suitable boundary condition must be imposed for the Laplacian to ensure that the normalization constant is finite. 

When the counting measure is replaced by the Lebesgue measure in \eqref{E:1.1}, the resulting law is that of the $1/\sqrt\beta$-multiple of the \textit{Gaussian Free Field} (GFF); we can thus think of the DG model as the GFF conditioned on taking integer values. (This is the reason why the word Gaussian appears in the name of this model, which is otherwise not Gaussian at all.) The GFF is a well studied process so a natural question is then: Under what conditions is the GFF a good approximation to the DG model? 

As it turns out, the answer to this varies depending on the underlying graph and the parameter~$\beta$. This is witnessed particularly by the most interesting case of~$\Lambda$ being a large box in~$\Z^2$ where the DG model exhibits a so-called \textit{roughening transition} (first established by Fr\"ohlich and Spencer~\cite{FS}): For~$\beta$ small the DG model is asymptotically GFF-like at large spatial scales (Bauerschmidt, Park and Rodriguez~\cite{BPR1,BPR2}) but,  as shown by a Peierls-type argument,  not so at all for~$\beta$ large. When characterized by infinite vs finite limit variance of the field, the two regimes meet at a single value of~$\beta$.  This was shown using Ginibre-type inequalities by Fr\"ohlich and Park~\cite{FP} although the conclusion is hard to glean from there. Aizenman, Harel, Peled and Shapiro~\cite[Theorem~1.2]{AHPS} state this with a different proof. 

\subsection{Overview}
While our prime interest rests with the behavior of the maximum and extremal values of the DG model in finite subsets of~$\Z^2$, extracting sharp results in the lattice model appears presently too hard.  (The only available sharp asymptotic is in Lubetzky, Martinelli and Sly~\cite{LMS} for large~$\beta$; for~$\beta$ small we only have bounds on the maximum from Wirth~\cite{W} and Garban and Sep\'ulveda~\cite{GS}.)  We will therefore resort to the hierarchical version of the DG model that shares many important features with the model on~$\Z^2$ but is more amenable to analysis. In a way, the hierarchical DG model is  a similar  approximation to the model on~$\Z^2$ as the Branching Random Walk is to the lattice~GFF. The understanding of the extremal behavior of Branching Random Walks has been instrumental for  later  results on the lattice GFF and we expect the same here as well. 

Our control of the DG model extends throughout the subcritical regime of inverse temperatures~$\beta$. In the hierarchical model, this is defined as the regime in which the iterative approach we rely on converges to a ``trivial'' fixed point and 
the model thus behaves as the GFF at large spatial scales. In order to achieve this within the framework of the renormalization group theory, which is our main tool, we had to develop a novel way to control the iterations of the effective potential that avoids linearization. 

Using the iterative analysis, we are then able to verify that the subcritical hierarchical DG model belongs to the universality class of the continuum-valued GFF: the suitably centered maximum is asymptotically a randomly shifted Gumbel and the extremal process is a (Gumbel) Poisson process with random intensity measure. The subcritical DG model is in fact so close to the GFF that the random shift/intensity is the \textit{same} as for the GFF except for a deterministic correction that depends on the scale of the system and can be attributed to a ``rounding error'' caused by the integer-valued nature of the underlying field. Due to a scale-dependence of this correction, the stated limit laws only exist along particular subsequences.

Our analysis of the extremal behavior relies on a coupling of the DG model to the GFF enabled by the iterative/renormalization group method. The idea to interpret renormalization iterations via coupling and use it to study extremal properties of the field has  in the meantime appeared elsewhere (e.g.,~\cite{BH,H,BGH})  but  our approach here is different.  Indeed,  we derive the limit result directly from the corresponding conclusion for the GFF, rather than by mimicking the proofs that did it for the GFF (which is what was done in~\cite{BH} on which the relevant conclusion in~\cite{H,BGH} is based). In addition, we aim directly at the full extremal process.

\section{Main results}
\noindent 
Let us move to the specific details of the problem at hand, starting with definitions and statement of convergence of the maximum.
\subsection{The model and convergence of the maximum}
\noindent
Given naturals~$n\ge0$ and~$b\ge2$ , the hierarchical DG model is defined over the hierarchical lattice
\begin{equation}
 \Lambda_n:=\{0,\dots,b-1\}^n 
\end{equation}
of depth~$n$. (For the connection with~$\Z^d$ we take~$b:=L^d$ for some natural~$L\ge2$, in which case~$\Lambda_n$ can be identified with a box of side-length~$L^{n}$.) The vertices of~$\Lambda_n$ are thus sequences  of length~$n$  taking values in~$\{0,\dots,b-1\}$. 

The hierarchical distance $d(x,y)$ between vertices $x=(x_1,\dots,x_n)$ and~$y=(y_1,\dots,y_n)$ is defined as 
\begin{equation}
 d(x,y):=\min\bigl\{j=0,\dots,n\colon(\forall i=1,\dots,n-j\colon x_i=y_i)\bigr\}.
\end{equation}
  The set  $\Lambda_n$ splits into~$b$ copies of~$\Lambda_{n-1}$ that are distance one from one another, each of which then splits into~$b$ copies of~$\Lambda_{n-2}$ at unit distance from one another, etc. This geometric structure is the basis of iterative approaches to hierarchical models.

The hierarchical Laplacian~$\Delta_n$ on~$\Lambda_n$ is naturally associated with the simple random walk on a $b$-ary tree of depth~$n$ killed when the walk exits the tree through the root. More precisely,~$\Delta_n$ is the generator of the Markov chain induced by observing the walk only on the leaves of the tree which are naturally identified with~$\Lambda_n$. 
The action of~$\Delta_n$ on~$f\colon\Lambda_n\to\R$ takes the explicit form
\begin{equation}
\label{E:1.3a}
\Delta_n f(x) =-\frakc_{n+1} f(x)+\sum_{k=1}^n \sum_{y\in \BB_k(x)}\frakc_k\bigl[\,f(y) - f(x)\bigr],
\end{equation}
where $\BB_k(x):=\{y\in\Lambda_n\colon d(x,y)\le k\}$ is the set of vertices within distance~$k$ from~$x$  and $\{\frakc_k\}_{k=1}^{n+1}$ is the sequence defined by
\begin{equation}
\frakc_k:=\frac1{b^k}\Biggl[\biggl(\,\sum_{j=0}^{k-1}b^j\biggr)^{-1}\!\!\!\!-\!\biggl(\,\sum_{j=0}^{k}b^j\biggr)^{-1}\Biggr],\quad k=1,\dots,n,
\end{equation}
and
\begin{equation}
\frakc_{n+1}:=\biggl(\,\sum_{j=0}^{n}b^j\biggr)^{-1}.
\end{equation}
 As is readily checked,~$\Delta_n$ is symmetric and, since the killing mechanism generates a non-trivial mass term  $\frakc_{n+1}f(x)$, it is  negative definite on~$\ell^2(\Lambda_n)$. 

The law of the hierarchical DG model on~$\Lambda_n$ is a probability measure~$P_{n,\beta}$ on~$\Z^{\Lambda_n}$ that takes the form
\begin{equation}
\label{E:1.1a}
P_{n,\beta}(\textd\varphi):=\frac1{\Sigma_n(\beta)}\texte^{\frac12\beta(\varphi,\Delta_n\varphi)}\prod_{x\in\Lambda_n}\#(\textd\varphi_x).
\end{equation}
Here $(\cdot,\cdot)$ denotes the canonical inner product in~$\ell^2(\Lambda_n)$ and~$\#$ is the counting measure on~$\Z$. 
The measure $P_{n,\beta}$ is well defined for all~$\beta>0$ thanks to the fact that, in light of~$\Delta_n$ being symmetric and negative definite, the normalization constant~$\Sigma_n(\beta)$ is finite. Two constants will play an important role throughout; namely,
\begin{equation}
\beta_\cc:=\frac{2\pi^2}{\log b}\quad\text{and}\quad\alpha:=\sqrt{2\log b}.
\end{equation}
Our first result is then:

\begin{theorem}[Limit law of the maximum]
\label{thm-1.1}
Given~$\beta>0$, let
\begin{equation}
\label{E:1.6a}
m_n:=\frac1{\sqrt\beta}\biggl[\sqrt{2\log b}\,n-\frac32\frac1{\sqrt{2\log b}}\log n\biggr].
\end{equation}
There exists a random variable~$\ZZ$ with~$\ZZ\in(0,\infty)$ a.s.\ and, for all~$\beta\in(0,\beta_\cc)$, there exists a constant $\hat c_\beta(0)\in(0,\infty)$ with the following property: For all sequences $\{n_k\}_{k\in\N}$ of natural numbers such that $n_k\to\infty$ and $s:=\lim_{k\to\infty}(m_{n_k}-\lfloor m_{n_k}\rfloor)$ exists,
\begin{equation}
\label{E:1.6}
P_{n_k,\beta}\Bigl(\,\max_{x\in\Lambda_{n_k}}\varphi_x\le \lfloor m_{n_k}\rfloor + u\Bigr)\,\,\underset{k\to\infty} \longrightarrow\,E\Bigl(\texte^{-\hat c_\beta(s)\ZZ\texte^{-\alpha\sqrt\beta\,u}}\Bigr),\quad u\in\Z,
\end{equation}
where $\hat c_\beta(s):=\hat c_\beta(0)\texte^{\alpha\sqrt\beta\, s}$.
\end{theorem}

The non-constancy of~$s\mapsto\hat c_\beta(s)$ is due to the discrete nature of the DG model. Apart from that, the result is  nearly  identical to that for the associated GFF which  is a Gaussian process on~$\Lambda_n$ with law given by
\begin{equation}
P'_{n,\beta}(\textd\varphi):=\frac1{\Sigma_n'(\beta)}\texte^{\frac12\beta(\varphi,\Delta_n\varphi)}\prod_{x\in\Lambda_n}\textd\varphi_x. 
\end{equation}
 (Another way to think of GFF is a (tree-indexed) Branching Random Walk with step distribution $\NN(0,1/\beta)$.) 
Indeed, the limit law for GFF is obtained from \eqref{E:1.6} by  replacing~$\lfloor m_n\rfloor$ by~$m_n$ and~$\hat c_\beta(s)$ by $(\alpha\sqrt\beta)^{-1}$, see \eqref{E:7.53a}; \rm the  conclusion  will then hold for all~$u\in\R$. The random variable~$\ZZ$, which translates to the random shift of the Gumbel law mentioned above, is the same for both processes and is independent of~$\beta$ and~$s$ (and the subsequence achieving~$s$). These considerations imply: 

\begin{corollary}
\label{cor-1.2}
Write~$\varphi'$ for the samples of the GFF under~$ P_{n,\beta}'$. For all $\beta\in(0,\beta_\cc)$ there exists~$a>0$ such that 
\begin{equation}
\sup_{u\in\R}\,\inf_{r\in[-a,a]}\,\biggl|\,P_{n,\beta}\bigl(\,\max_{x\in\Lambda_{n}}\varphi_x\le 
u\bigr)
- P_{n,\beta}'\Bigl(\,\max_{x\in\Lambda_{n}}\varphi_x'\le 
u+r\Bigr)\biggr|\,\,\underset{n\to\infty}\longrightarrow\,\,0.
\end{equation}
In particular, for each~$\beta\in(0,\beta_\cc)$, the maxima of the DG model and  the GFF   can be coupled to within a bounded distance, with probability tending to one as $n\to\infty$.
\end{corollary}

\subsection{Extremal process}
As already alluded to, our proof of Theorem~\ref{thm-1.1} is based on a coupling of the DG model to the GFF that  naturally arises from the renormalization group analysis of the DG model. Unfortunately, the coupling does not necessarily preserve the maximizer(s) and so we will have to work with the full extremal process of the DG model for which we thus get a limit law as well. 

In order to make precise statements, we need more definitions. First, with each vertex $x=(x_1,\dots,x_n)\in\Lambda_n$ we can associate the real number
\begin{equation}
\label{E:1.9a}
 [x]_n:=\sum_{i=1}^n x_i b^{-i}. 
\end{equation}
The resulting map~$x\to[x]_n$  embeds~$\Lambda_n$ injectively into~$[0,1]$. 

Second, in order to talk about point processes and their weak convergence, let~$\MM(\scrX)$ denote the space of Radon measures on a topological space~$\scrX$.
We endow $\MM(\scrX)$ with the topology of vague convergence which is generated by integrals of the measure against continuous, compactly-supported functions. This in turn permits us to talk about random elements of~$\MM(\scrX)$ and weak convergence thereof. We will also write~$\MM_\N(\scrX)$ for the set of integer-valued measures in~$\MM(\scrX)$ which, assuming~$\scrX$ is Hausdorff, are necessarily concentrated on  locally-finite sets  of points. \rm

Third, naturally associated with each sample~$\varphi$ of the DG model is a random element of $\MM_\N([0,1]\times\R)$ called the extremal process
\begin{equation}
\eta_n:=\sum_{x\in\Lambda_n}\delta_{[x]_n}\otimes\delta_{\varphi_x-\lfloor m_n\rfloor},
\end{equation}
where~$m_n$ is as in \eqref{E:1.6a}.
This process records the  position (via  $x\mapsto[x]_n$) and values of the DG model reduced by~$\lfloor m_n\rfloor$. (Integrating~$\eta_n$ against compactly supported functions will pick up only the values that are close to~$m_n$.)  

Fourth and finally, let $\PPP(\nu)$, for a given $\sigma$-finite measure~$\nu$, denote the Poisson point process with intensity~$\nu$ which, if~$\nu$ is itself random, is sampled conditionally on~$\nu$. (The process is defined on a space that carries the measure~$\nu$.) We then have:

\begin{theorem}[Limit extremal process]
\label{thm-1.2}
There exists an a.s.-finite random Borel measure~$Z$ on~$[0,1]$ with $Z(A)>0$ a.s.\ for all non-empty (relatively) open~$A\subseteq[0,1]$ and, for all $\beta\in(0,\beta_\cc)$,  there exist a constant $\tilde c_\beta\in(0,\infty)$ and  a probability measure~$\nu_\beta$ on $\MM_\N(\Z)$  satisfying
\begin{equation}
\label{E:2.14ii}
\nu_\beta\biggl(\Bigl\{\zeta\in\MM_\N(\Z)\colon\,
\zeta(\{0\})\ge1\,\,\wedge\,\,\zeta((0,\infty))=0 \Bigr\}\biggr)=1
\end{equation}
 such that the following holds:  For all sequences $\{n_k\}_{k\in\N}$ of natural numbers such that $n_k\to\infty$ and $s:=\lim_{k\to\infty}(m_{n_k}-\lfloor m_{n_k}\rfloor)$ exists, 
\begin{equation}
\label{E:1.9}
\eta_{n_k}
\,\,\,\underset{k\to\infty}\lawarrow\,\,\,
\sum_{i,j\ge1}\delta_{x_i}\otimes\delta_{h_i+t_j^{(i)}},
\end{equation}
where $\{(x_i,h_i)\}_{i\ge1}$ enumerates the points in a sample from
\begin{equation}
\text{\rm PPP}\Bigl(\tilde c_\beta\,\texte^{\alpha\sqrt\beta\,s}\,Z\otimes \sum_{n\in\Z}\texte^{-\alpha\sqrt\beta\, n}\delta_n\Bigr),
\end{equation}
and $\{t_j^{(i)}\}_{j\ge1}$ enumerates the points in the $i$-th member of the sequence $\{t^{(i)}\}_{i\ge1}$ of i.i.d.\ samples from~$\nu_\beta$ that are independent  of~$Z$  and $\{(x_i,h_i)\}_{i\ge1}$.  The set
\begin{equation}
\bigl\{h_i+t_j^{(i)}\colon i,j\ge1\bigr\}\cap[u,\infty)
\end{equation}
is finite almost surely for each~$u\in\R$.  Almost every sample of~$Z$ is diffuse. 
\end{theorem}

The convergence \eqref{E:1.9} means that, for any $f\colon [0,1]\times\Z\to\R_+$ that is continuous in the first coordinate and supported on a finite set of values in the second coordinate, the quantity $E(\texte^{-\int f\textd\eta_{n_k}})$ tends to
\begin{equation}
\label{E:2.18a}
E\,\biggl(\exp\Bigl\{-\tilde c_\beta\, \texte^{\alpha\sqrt\beta\, s}\int
Z(\textd x)\otimes\nu_\beta(\textd\zeta)\sum_{n\in \Z}\texte^{-\alpha\sqrt\beta\, n}\bigl(1-\texte^{-\int f(x,n+\cdot)\textd\zeta})\Bigr\}\biggr)
\end{equation}
in the limit as~$k\to\infty$, where the expectation is over~$Z$. 
Theorem~\ref{thm-1.1} will actually be derived from this  formulation  by noting that the event $\max_{x\in\Lambda_n}\varphi_x\le\lfloor m_n\rfloor+u$ translates into the extremal process~$\eta_n$ not charging the set~$[0,1]\times(u,\infty)$.  Tightness arguments and a  calculation then  show  that the objects in \eqref{E:1.6} are given by
\begin{equation}
\label{E:1.13u}
\ZZ:=Z\bigl([0,1]\bigr)
\end{equation}
and 
\begin{equation}
\label{E:1.14u}
\hat c_\beta(s):=\tilde c_\beta\,\texte^{\alpha\sqrt\beta\, s}\sum_{n\in\Z}\texte^{-\alpha\sqrt\beta\,n}\,\nu_\beta\biggl(\Bigl\{\zeta\in\MM_\N(\Z)\colon\zeta\bigl((-n,\infty)\bigr)\ge1\Bigr\}\biggr),
\end{equation}
where,  as we prove,  the sum on the right converges. (A detailed derivation will be given in the proof of Theorem~\ref{thm-1.2} and Theorem~\ref{thm-1.1}.) 

As  it turns out,  the~$Z$ measure \textit{coincides} with the random intensity governing the scaling limit of the extremal process of the normalized GFF (i.e., that with $\beta:=1$); see Theorem~\ref{thm-6.1}. The cluster law~$\nu_\beta$ of the DG model is derived from the corresponding cluster law~$\nu$ of the normalized GFF, albeit through a non-constructive limit procedure that makes the connection inexplicit.  While we can show that the samples from~$\nu_\beta$ can have arbitrarily large total mass, the lack of explicit handle prevented us from showing that $\zeta(\Z)=\infty$ $\nu_\beta$-a.s., which we nonetheless believe to be the case.

\subsection{Discussion}
We proceed with a discussion of the above statements. 
Perhaps the most surprising aspect of our results is the aperiodic dependence   on the underlying scale~$n$, which causes a limit law to exist only for subsequences where $m_n-\lfloor m_n\rfloor$ converges. To understand why this happens,   note that near-maximal values of the GFF occur at ``heights'' of the asymptotic form
\begin{equation}
m_n+h_i+t^{(i)}_j,\quad i,j\ge1,
\end{equation}
where $\{h_i\}_{i\ge1}$ are (Gumbel) Poisson points and $\{t_j^{(i)}\}_{j\ge1}$ are correlated (and non-positive) but independent of~$\{h_i\}_{i\ge1}$ and independent for different~$i$.   The coupling to the DG model shows that the extremal values of the DG model closely follow those of the GFF but are constrained to be integers. As a result, the DG model extremal points occur    at integer ``heights'' at or near the values
\begin{equation}
\lfloor m_n+ h_i+t^{(i)}_j\rfloor,\quad i,j\ge1.
\end{equation}
%
  Here, the fractional part of~$m_n$ effectively acts as a uniform shift of $\{h_i\}_{i\ge1}$ and remains relevant in the limit, albeit only as a ``rounding error.'' (A change of variables propagates this into the prefactor in \eqref{E:1.6} and \eqref{E:2.18a}.)  It is as if the DG model maximum wanted to concentrate around~$m_n$, even though~$m_n$ itself is not an integer.

 Moving over to more general aspects,  our results demonstrate a very close relationship between the extremal properties of the DG model and those of the GFF throughout the parameter regime~$\beta\in(0,\beta_\cc)$. This should be regarded as a  particular  statement of \textit{universality}. Notwithstanding, the reader may naturally wonder what happens when~$\beta$ does not lie in this interval.

In our proofs, $\beta_\cc$ marks the largest value of~$\beta$ at which the coarse-graining/re\-norma\-lization-group iterations converge to a ``trivial'' fixed point. This is usually expressed via the behavior of the two-point correlation function or the so-called fractional charge. We instead articulate it in much stronger terms by constructing a coupling of the DG model and the GFF that keeps the two fields  (at any given point)  within a tight distance of each other. If a suitably discretized version of the GFF is used, the two fields even agree on a positive fraction of all vertices.

At~$\beta=\beta_\cc$ the convergence to a ``trivial'' fixed point still takes place but only at a polynomial rate which precludes control of the coupling at the desired level.  As we show in a follow up paper~\cite{BiH}, for~$\beta>\beta_\cc$ (with $\beta-\beta_\cc$ small) the iterations converge to a ``non-trivial'' fixed point. This means that the  ``correct'' approximation of the DG model is not the Branching Random Walk with step distribution~$\NN(0,1/\beta)$ but rather  a Branching Markov Chain  whose steps have a Gaussian density of~$\NN(0,1/\beta)$ multiplied by a positive $1$-periodic function that depends on the current state of the chain. 

Unlike Branching Random Walks, the extremal behavior of Branching Markov Chains has been controlled only in a few cases (e.g., the local time of simple random walk on regular trees~\cite{BL4}) and no  robust  theory exists.
 That being said, our results in~\cite{BiH} show that~$\varphi$ is logarithmically correlated throughout. Indeed, we have
\begin{equation}
\Cov(\varphi_x,\varphi_y) = \begin{cases}
\sigma^2(\beta)\bigl[n-d(x,y)\bigr]+O(1),\qquad&\text{if }\beta\ne\beta_\cc,
\\
\beta_\cc^{-1}\bigl[n-d(x,y)\bigr]-\bar c\log\bigl(\frac {n}{1+d(x,y)}\bigr)+O(1),\qquad&\text{if }\beta=\beta_\cc,
\end{cases}
\end{equation}
where~$\sigma^2(\beta)$ is the quantity that obeys
\begin{equation}
\sigma^2(\beta)\begin{cases}
=1/\beta,\qquad&\text{if }\beta\le\beta_\cc,
\\
<1/\beta,\qquad&\text{if }\beta>\beta_\cc,
\end{cases}
\end{equation}
and~$\bar c$ is a positive (explicit) constant while~$O(1)$ is bounded uniformly in~$n$ and~$x,y\in\Lambda_n$.
Using this we can 
 offer a sophisticated guess concerning the maximum of the DG model for~$\beta>\beta_\cc$: The limit law of the maximum remains formally the same discrete-Gumbel like (i.e., \eqref{E:1.6} applies), albeit with the various objects --- namely, the centering~$m_n$, the random variable~$\ZZ$ and the tail exponent~$\alpha\sqrt\beta$ --- now having a different $\beta$-dependence than in the subcritical regime.

Our guess is based on the assumption that the Markov chain underlying the approximating Branching Markov Chain process behaves closely to a random walk albeit now with the process at time~$n$ having variance~$\sigma(\beta)^2 n$. Borrowing on  arguments from the theory of Branching Random Walks, this should mean that the probability that the maximum exceeds value~$r$ is asymptotic~to
\begin{equation}
\label{E:1.14}
b^n\frac1n\frac c{\sqrt n}\texte^{-\frac{r^2}{2\sigma(\beta)^2 n}},
\end{equation}
where~$b^n$ accounts for the possible positions of the maxima, the factor $1/n$ arises from a ``ballot theorem'' that accounts for the entropic-repulsion effect caused by conditioning the remaining vertices to have value below~$r$ and the remainder of the term (with~$c$ a positive factor  independent of~$r$)  is the Gaussian density with variance~$\sigma(\beta)^2 n$. 

For the specific choice
\begin{equation}
r:=\sigma(\beta)\biggl[\sqrt{2\log b}\,n-\frac32\frac1{\sqrt{2\log b}}\log n\biggr]+ u
\end{equation}
the quantity in \eqref{E:1.14} is proportional to~$\texte^{-\tilde\alpha(\beta)u}$ with
\begin{equation}
\tilde\alpha(\beta):=\frac{\sqrt{2\log b}}{\sigma(\beta)}.
\end{equation}
We thus expect that the sole effect on the centering and the tail exponent to be the replacement of~$\beta$ by~$\sigma(\beta)^{-2}$. The situation with the random variable~$\ZZ$ is less clear as it arises directly, via the so-called derivative martingale, from the approximating Branching Markov Chain and so it likely depends non-trivially on~$\beta$. (The lack of this dependence for~$\beta<\beta_\cc$ is due to the $\beta$-dependent factors scaling out explicitly.)

 The situation at~$\beta=\beta_\cc$ is quite unclear due to the second-order correction to the covariance arising in this case and so we refrain from making any conjectures. A rigorous control of this case will undoubtedly be quite interesting. We do expect the correction to the variance to propagate to the centering sequence for the maximum.  

The reader unfamiliar with hierarchical models  has likely wondered why the covariance structure (and presumably also the extremal behavior) remains qualitatively  close in the supercritical and subcritical regimes. After all, the subcritical regime marks the delocalized phase of the model and the supercritical regime  should thus  correspond  to the localized one. While this  is  the case in the lattice setting, it is not in the hierarchical version.  This is partly due to the hierarchical Laplacian being long-range, as the  coefficients in \eqref{E:1.3a} decay as $\frakc_k\asymp b^{-2k}$ which  translates to polynomial decay in distance in~$\ell^\infty$-distance under the embedding of~$\Lambda_n$ into~$\Z^2$.  However, the situation is more subtle because, as shown by Fr\"ohlich and Zegarlinski~\cite{FZ} and later Marchetti~\cite{Marchetti}, the DG model in~$d=1$ with $1/r^2$-decay of interaction does exhibit both the localized and the delocalized phases. See Garban~\cite{Garban} and Coquille, Dario and Ny~\cite{CDN} for recent work on the one-dimensional DG model.

We close the discussion by noting that the coupling argument presented here is capable of handling other extremal questions one may consider for the hierarchical DG model. One of these is that of ``intermediate'' level sets, a.k.a., the thick points, which are those where the field exceeds $\lambda \beta^{-1/2} \sqrt{2\log b}\, n$ for~$\lambda \in(0,1)$.  Precise control of these has been achieved in a number of contexts; e.g., for the lattice GFF in Biskup and Louidor~\cite{BL5} and random-walk local time in Jego~\cite{Jego} and Abe, Biskup and Lee~\cite{AB2,ABL}. 

\subsection{Outline}
The remainder of this paper is organized as follows: In Section~\ref{sec2} we explain the two main technical inputs underlying our results; namely, the iterative (a.k.a.\ renormalization group) approach to the DG model and the ensuing coupling of the DG model to the~GFF. These steps are where the restriction to $\beta\in(0,\beta_\cc)$ originates from. Details of  these  are worked out in Sections~\ref{sec3} (iterations) and~\ref{sec4} (coupling). In the remaining two sections we then address the proof of the above theorems starting from a proof of tightness (Section~\ref{sec5}) and then convergence of the extremal process (Section~\ref{sec6}).


\section{Two main new ideas}
\label{sec2}
\noindent
 Here we commence the exposition of our proofs by developing our version of the renormalization group approach to the hierarchical DG model and stating its main conclusions (which come as Theorems~\ref{thm-3.4} and~\ref{thm-4.1}). 

\subsection{Iterative description of hierarchical DG model}
\label{sec-2.1}\noindent
 Let us start  by explaining the structure of the hierarchical DG model that drives the iterative description thereof. This will allow us to give several good reasons for the particular choice of the hierarchical Laplacian in \eqref{E:1.3a}. 

The study of hierarchical models was initiated in Dyson~\cite{Dyson}. The hierarchical setting is generally more friendly to coarse-graining arguments than their lattice counterpart which is why hierarchical models have served as a useful preliminary tool in real-space renormalization group theory (Brydges~\cite{Brydges-notes}). While the hierarchical version is for many models just a convenient approximation, other models --- such as the two-dimensional GFF or long-range percolation (see, e.g., Biskup and Krieger~\cite{BK} or Hutchcroft~\cite{Hutchcroft1}) --- seem to exhibit  an approximate  hierarchical structure naturally.

In the past, the hierarchical DG model has been studied mainly in its ``soft'' variant called the Sine-Gordon model (Marchetti and Perez~\cite{MP}, Guidi and Mar\-chetti~\cite{GM}, Benfatto and Renn~\cite{BR}) with calculations and results often presented for its dual version of the lattice Coulomb gas. We will work solely with the DG model although the forthcoming derivations apply to a whole class of hierarchical models including the Sine-Gordon model. (We will give appropriate definitions in Section~\ref{sec-2.4}.) 

The iterative description of the DG model is based on an integral identity whose statement requires additional notation.
First, let us add suffix~$n$ to the notation of the canonical inner product in~$\ell^2(\Lambda_n)$ and write it as $(\cdot,\cdot)_n$. Recall also our notation
\begin{equation}
\BB_k(x):=\{y\in\Lambda_n\colon d(x,y)\le k\}
\end{equation}
for the ball of (hierarchical) radius~$k$ centered at~$x$. Note that~$\BB_k(x)$ has~$b^k$ elements. Introduce the map $m\colon\Lambda_n\to\Lambda_{n-1}$ defined by
\begin{equation}
\label{E:2.2a}
m(x):=(x_1,\dots,x_{n-1})\quad\text{when}\quad x=(x_1,\dots,x_n).
\end{equation}
We then have:

\begin{lemma}
Given a natural $n\ge1$ and reals $\lambda_0,\dots,\lambda_{n}>0$, consider the linear operator~$ \LL_n$  on~$\ell^2(\Lambda_n)$ that acts on $f\colon\Lambda_n\to\R$ as
\begin{equation}
\label{E:2.3a}
 \LL_n  f(x):=\sum_{k=1}^n \frac{\lambda_{k-1}^{-1}-\lambda_{k}^{-1}}{b^{k}}\sum_{y\in  \BB_k(x)}\bigl[\,f(y) - f(x)\bigr] - \lambda_{n}^{-1}f(x).
\end{equation}
For each~$n\ge1$ there exists $\frakz_n\in(0,\infty)$ such that for all~$\varphi\colon\Lambda_n\to\R$,
\begin{equation}
\label{E:2.4a}
\texte^{\frac12(\varphi,\, \LL_n\varphi)_n} = \frakz_n\int_{\R^{\Lambda_{n-1}}} \texte^{-\frac1{2\lambda_0}\sum_{x\in\Lambda_n}(\varphi_x-\psi_{m(x)})^2}\texte^{\frac12(\psi,\, \LL'_{n-1}\psi)_{n-1}}\prod_{z\in\Lambda_{n-1}}\textd\psi_z ,
\end{equation}
where $\LL_{n-1}'$  is defined using reals $\lambda_0',\dots,\lambda'_{n-1}$ that are uniquely determined from
\begin{equation}
\label{E:2.5b}
\lambda'_{0}=\frac{\lambda_{1}-\lambda_0}b
\end{equation}
and
\begin{equation}
\label{E:2.6b}
\lambda'_{k}-\lambda'_{k-1}=\frac{\lambda_{k+1}-\lambda_{k}}b,\quad k=1,\dots,n-1.
\end{equation}
\end{lemma}

\begin{proofsect}{Proof}
The formula \eqref{E:2.4a} can be understood informally as a convolution identity for multivariate Gaussians except that interpreting these Gaussians on the same space makes some of them singular and so talking about densities is difficult.
We will instead use the fact that it suffices  to  prove equality under integration of \eqref{E:2.4a} against~$\texte^{(\varphi,f)_n}$ with some~$\frakz_n$ independent of~$f\in\ell^2(\Lambda_n)$, for all such functions~$f$. 

Using the substitution $\wt\varphi:=\varphi+\LL_n^{-1} f$,  which is well defined  as~$\LL_n$ will be shown to be invertible,  on the left-hand side we get
\begin{equation}
\int_{\R^{\Lambda_n}}\texte^{\frac12(\varphi,\,\LL_n\varphi)_n+(\varphi,f)_n}\prod_{x\in\Lambda_n}\textd\varphi_x = \frakz'\,\texte^{-\frac12(f,\,\LL_n^{-1} f)_n},
\end{equation}
where~$\frakz'$ is a positive and finite quantity independent of~$f$. Letting~$f'\colon\Lambda_{n-1}\to\R$ be the function defined by $f'(y):=\sum_{x\in\Lambda_n} f(x)1_{\{m(x)=y\}}$, on the right-hand side the substitutions $\wt\varphi:=\varphi- \psi_{m(\cdot)}-\lambda_0 f$ followed by $\wt\psi:= \psi+ (\LL_{n-1}')^{-1} f'$  give
\begin{multline}
\quad\int_{\R^{\Lambda_n}\times\R^{\Lambda_{n-1}}}\texte^{-\frac1{2\lambda_0}\sum_{x\in\Lambda_n}(\varphi_x-\psi_{m(x)})^2}\texte^{\frac12(\psi,\,\LL_{n-1}'\psi)_{n-1}+(\varphi, f)_n}\prod_{z\in\Lambda_n}\textd\varphi_z\prod_{z\in\Lambda_{n-1}}\textd \psi_z
\\
=\frakz''\,\texte^{\frac12\lambda_0 (f,f)_n-\frac12(f',\,(\LL_{n-1}')^{-1} f')_{n-1}},
\quad
\end{multline}
where~$\frakz''$ is again independent of~$f$ and where we used that $(\psi_{m(\cdot)},f)_n=(\psi,f')_{n-1}$. To prove \eqref{E:2.4a}, we thus have to show that
\begin{equation}
\label{E:2.9a}
(f,\,\LL_n^{-1} f)_n = -\lambda_0 (f,f)_n+\bigl(f',\,(\LL_{n-1}')^{-1} f'\bigr)_{n-1}
\end{equation}
holds for all~$f\in\ell^2(\Lambda_n)$.

Define the operator~$Q_k$ by $Q_kf(x):=\frac1{b^k}\sum_{y\in \BB_k(x)}f(y)$ for~$k=1,\dots,n$ and let~$Q_0$ be the identity and~$Q_{n+1}:=0$. Note that $Q_k$ is, for each~$k=0,\dots,n+1$, an orthogonal projection (in~$\ell^2(\Lambda_n)$) such that $Q_kQ_\ell = Q_\ell Q_k = Q_k$ whenever~$\ell\le k$. In particular,  for all~$k,\ell=0,\dots,n$ we have
\begin{equation}
\label{E:2.8a}
(Q_\ell-Q_{\ell+1})(Q_k-Q_{k+1})=\delta_{k,\ell}(Q_k-Q_{k+1})
\end{equation}
and so $\{Q_{k}-Q_{k+1}\colon k=0,\dots,n\}$ are orthogonal projectors on orthogonal subspaces that, in light of $\sum_{k=0}^n(Q_k-Q_{k+1})$ being the identity, span all of~$\ell^2(\Lambda_n)$.

With the help of the convention $(\lambda_{-1})^{-1}:=0$, the operator~$ \LL_n$  now takes the form
\begin{equation}
\begin{aligned}
 \LL_n  &= \sum_{k=0}^{n}(\lambda_{k-1}^{-1}-\lambda_{k}^{-1})(Q_k-Q_0) - \lambda_n^{-1}Q_0
\\&=\sum_{k=0}^{n}(\lambda_{k-1}^{-1}-\lambda_{k}^{-1})Q_k
=-\sum_{k=0}^n\lambda_k^{-1}(Q_k-Q_{k+1}).
\end{aligned}
\end{equation}
This  and  \eqref{E:2.8a} imply
\begin{equation}
\LL_n^{-1} = -\sum_{k=0}^n\lambda_k(Q_k-Q_{k+1}).
\end{equation}
In particular, $\LL_n$  is invertible and negative definite.

The previous formula gives
\begin{equation}
\label{E:2.13a}
\LL_n^{-1} = -\lambda_0 Q_0-\sum_{k=1}^n(\lambda_{k}-\lambda_{k-1})Q_k.
\end{equation}
Observing that, for each~$k\ge1$ we have $(f,Q_k f)_{n} = b^{-1}(f',Q_{k-1}f')_{n-1}$, using \twoeqref{E:2.5b}{E:2.6b} we now readily verify \eqref{E:2.9a} and thus conclude the proof.
\end{proofsect}

\subsection{What makes the hierarchical Laplacian special?}
The operator defined in \eqref{E:2.3a} is an instance of a hierarchical operator; i.e., one that commutes with all the~$Q_j$ defined in the previous proof. As noted above, $\LL_n$  is invertible and negative definite for all $\lambda_0,\dots\lambda_n>0$. If in fact $0<\lambda_0\le\dots\le\lambda_n$, then we can write
\begin{equation}
\label{E:2.14a}
\lambda_k^{-1}:=\lambda_0^{-1} P(\tau>k),\quad k=0,\dots,n,
\end{equation}
for a random variable~$\tau$ taking values in~$\{1,2\dots,n+1\}$. As
\begin{equation}
\lambda_{k-1}^{-1}-\lambda_k^{-1}=\lambda_0^{-1}P(\tau=k)
\end{equation}
 we can then view~$\LL_n$  as the generator of a Markov chain on~$\Lambda_n$ induced (by observing only visits to the leaves) by a Markov chain on the full $b$-ary tree~$\T_n:=\bigcup_{k=0}^n\Lambda_k$ of depth~$n$ whose transition probabilities have all the symmetries of~$\T_n$ and~$\tau$ represents the maximal height reached by that chain between two successive visits to the leaves; namely, the set~$\Lambda_n$.  Note that  $\tau>n$ marks the event that the chain exits~$\T_n$ through the root and ``dies.'' For~$\LL_n$,  this imitates the effect of Dirichlet boundary conditions.

Any Markov chain  respecting the symmetries of~$\T_n$ is determined by a  collection of numbers~$\{(p_k,q_k)\}_{k=0}^n$, with~$q_0=0$ and~$p_k>0$ and $p_k+bq_k=1$ for each~$k$, that have the following meaning:~$p_k$ is the probability that the chain at a site at height~$k$ above the leaves takes a step towards the root while~$b q_k$ is the (total) probability that it takes a step away from the root. In this parametrization we have
\begin{equation}
\label{E:2.14b}
P(\tau>k) = \biggl[\sum_{\ell=0}^k\prod_{j=1}^\ell\frac{b q_j}{p_j}\biggr]^{-1}, \quad k=0,\dots,n,
\end{equation}
and the generator associated with the process induced by observing the walk only on the leaves is then~$\LL_n$  with $\lambda_k$ as in~\eqref{E:2.14a}. (The prefactor~$\lambda_0$ just changes the speed with which the process moves in continuous time.) 

The above representations permit us to explain the special role that the hierarchical Laplacian~$\Delta_n$ plays among the operators of the form \eqref{E:2.3a}. Indeed, using the simple random walk on~$\T_n$ in place of the above Markov chain boils down to $p_k=q_k$ for each~$k=1,\dots,n$ and so, using \twoeqref{E:2.14a}{E:2.14b},
\begin{equation}
\label{E:2.15a}
\lambda_k=\lambda_0\sum_{j=0}^{k}b^j,\quad k=0,\dots,n.
\end{equation}
It follows that, in this case,~$\LL_n$  coincides with a $1/\lambda_0$-multiple of the hierarchical Laplacian~$\Delta_n$ defined in \eqref{E:1.3a}.

Another way to see the relevance of the hierarchical Laplacian stems from the following question: What specific choice (if any) of the coefficients $\lambda_0,\lambda_1,\dots$ is preserved by the transformation \twoeqref{E:2.5b}{E:2.6b}? A calculation reveals that this happens when
\begin{equation}
\lambda_k-\lambda_{k-1} = \lambda_0 b^{k},\quad k=1,\dots,n,
\end{equation}
which is again solved uniquely by \eqref{E:2.15a}. It follows that the hierarchical Laplacian reproduces itself under the convolution identity \eqref{E:2.4a}.

Finally, yet another fact singling out the particular form of~$\Delta_n$ is that plugging \eqref{E:2.15a} with~$\lambda_0:=1$ in \eqref{E:2.13a} yields
\begin{equation}
(\delta_x,-\Delta_n^{-1}\delta_y) =\lambda_0\delta_{xy}+\sum_{k=d(x,y)\vee1}^n(\lambda_k-\lambda_{k-1})\frac1{b^k}=n+1-d(x,y),
\end{equation}
where we used that~$(\delta_x,Q_k\delta_y)$ vanishes when~$d(x,y)>k$ and equals~$b^{-k}$ otherwise. In light of the hierarchical distance translating  (roughly)  into the logarithm of the Euclidean distance  under  the identification of~$\Lambda_n$ for~$b:=L^{ 2}$ with a square box of side-length~$L^n$ in~$\Z^{ 2}$, this shows that $\Delta_n^{-1}$ behaves much like the inverse of the Laplacian in the box with Dirichlet boundary conditions.

\subsection{Representation as a tree-indexed Markov chain}
In light of the preceding discussion, we will henceforth focus on the case when the sequence $\lambda_0,\dots,\lambda_n$ takes the form \eqref{E:2.15a} with~$\lambda_0:=\beta^{-1}$, i.e., $\LL_n=\beta\Delta_n$. The integral formula \eqref{E:2.4a} then allows for iterative computation of the partition function
\begin{equation}
\Sigma_n(\beta):=\sum_{\varphi\in\Z^{\Lambda_n}}\texte^{\frac\beta2(\varphi,\Delta_n\varphi)}
\end{equation}
normalizing the measure \eqref{E:1.1a}. As a starting point we observe
\begin{equation}
\begin{aligned}
\Sigma_n(\beta)&=\frakz_n\int_{\R^{\Lambda_{n-1}}} \biggl(\,\sum_{\varphi\in\Z^{\Lambda_n}}\texte^{-\frac\beta2\sum_{x\in\Lambda_n}(\varphi_x-\varphi_{m(x)}')^2}\biggr)\texte^{\frac\beta2(\varphi',\Delta_{n-1}\varphi')_{n-1}}\prod_{z\in\Lambda_{n-1}}\textd\varphi_z'
\\
&=\frakz_n\int_{\R^{\Lambda_{n-1}}} \texte^{-\sum_{z\in\Lambda_{n-1}}bv_0(\varphi_z')}\texte^{\frac\beta2(\varphi',\Delta_{n-1}\varphi')_{n-1}}\prod_{z\in\Lambda_{n-1}}\textd\varphi_z',
\end{aligned}
\end{equation}
where we have set
\begin{equation}
\label{E:2.21}
\texte^{-v_0(z)}:=\sum_{n\in\Z}\texte^{-\frac\beta2(z-n)^2}
\end{equation} 
and used that $\LL_{n-1}' =\beta\Delta_{n-1}$ when $\LL_n=\beta\Delta_n$.
Defining $v_1,\dots,v_n$ recursively via
\begin{equation}
\label{E:2.2}
\texte^{-v_{k+1}(z)}:=E\bigl(\texte^{-b v_k(z+\zeta)}\bigr)\quad\text{for}\quad \zeta\laweq\NN(0,1/\beta),
\end{equation}
 further  applications of \eqref{E:2.4a} then yield
\begin{equation}
\label{E:2.23a}
\begin{aligned}
\Sigma_n(\beta)&=\frakz_n\frakz_{n-1}(2\pi/\beta)^{|\Lambda_{n-1}|/2}\int_{\R^{\Lambda_{n-2}}} \texte^{-\sum_{z\in\Lambda_{n-2}}bv_1(\varphi_z'')}\texte^{\frac\beta2(\varphi'',\Delta_{n-2}\varphi'' )_{n-2}}\prod_{z\in\Lambda_{n-2}}\textd\varphi_z''
\\
&=\dots=\frakz_n\Bigl(\,\prod_{k=0}^{n-1}\frakz_k(2\pi/\beta)^{|\Lambda_{k}|/2}\Bigr)\texte^{-v_{n}(0)},
\end{aligned}
\end{equation}
where the factors $(2\pi/\beta)^{|\Lambda_k|/2}$ arise from the normalization of the Gaussian measure in~\eqref{E:2.2} which is not included in \eqref{E:2.4a}.

The calculation \eqref{E:2.23a} is not of much use in its own right --- after all, the factors $\frakz_k$ are not very explicit and the value of the partition function is not particularly illuminating. What people in mathematical physics usually do is to adapt the procedure to calculate expectations of suitable observables. A problem is that the observables keep chaining through the iterations and so one has to follow their ``flow'' separately. In all of this, not much attention is paid or meaning assigned to the auxiliary fields one integrates over; indeed, they are just tools in a massive inductive argument.

Our aim here is different in that we will keep track of the auxiliary fields and, in fact, use them to give a convenient representation of the DG model. 
Write~$\frakg_{1/\beta}$ for the law of~$\NN(0,1/\beta)$ and, for each~$k \in \{0, \dots, n\}$  and~$\varphi\in\R$, let $\frakq_k(\cdot|\varphi)$ be the Borel probability measure on~$\R$ given in infinitesimal form by
\begin{equation}
\label{E:1.23}
\frakq_{k}(\textd\zeta|\varphi):=\begin{cases}
\texte^{v_{k}(\varphi)-bv_{k-1}(\varphi+\zeta)}\,\frakg_{1/\beta}(\textd\zeta),\qquad&\text{if }k\ge1,
\\
\texte^{v_0(\varphi)-\frac\beta2\zeta^2}\,\#(\varphi+\textd\zeta),\qquad&\text{if }k=0,
\end{cases}
\end{equation}
where, we recall,~$\#$ denotes the counting measure on~$\Z$. Note that \eqref{E:2.21} and \eqref{E:2.2} ensure that each~$\frakq_k$ is a Markov kernel on~$\R$. Let~$m^k(\cdot)$ denote the~$k$-fold iteration of function~$m$ from \eqref{E:2.2a} and $\Lambda_0 := \{\rho\}$, where $\rho$ denote the root of the tree.  We then have:

\begin{lemma}
\label{lemma-2.2}
Given~$n\ge1$, let $\{\zeta_k(x)\colon x\in\Lambda_{k},\, k=0, \dots,n\}$ be a family of random variables with joint law given (again in infinitesimal form) by
\begin{equation}
\label{E:2.25}
 \frakq_n\bigl(\textd \zeta_0(\rho)\,\big|\,0\bigr)\otimes \bigotimes_{k=1}^{n}\bigotimes_{x\in\Lambda_k} \frakq_{n-k}\bigl(\textd\zeta_k(x)\,\big|\varphi_{k-1}(m(x))\bigr).
\end{equation}
 Define, for $k=0,\dots,n$,
\begin{equation}
\label{E:2.26}
\varphi_{k}(x):= \sum_{j=0}^k \zeta_j\bigl( m^{k-j}(x)\bigr),\quad x\in\Lambda_k.
\end{equation}
Then~$\{\varphi_n(x)\colon x\in\Lambda_n\}$ has the law $P_{n,\beta}$ of the DG model as defined in \eqref{E:1.1a}.
\end{lemma}

\begin{proofsect}{Proof} 
After some careful bookkeeping we find that the prefactors in \eqref{E:1.23} mostly cancel out and the measure in \eqref{E:2.25} can be cast as
\begin{equation}
\texte^{ v_{n}(0)  -\frac\beta2\sum_{x\in\Lambda_n}\zeta_n(x)^2}
\biggl[\,\bigotimes_{x\in\Lambda_n}\#\bigl(\varphi_{n-1}(m(x))+\textd\zeta_n(x)\bigr)\biggr]\otimes\biggl[\,\bigotimes_{k=0}^{n-1}  \bigotimes_{z\in\Lambda_k}\frakg_{1/\beta}\bigl(\textd\zeta_k(z)\bigr)\biggr],
\end{equation}
where~$\varphi_{n-1}$ is derived from $\{\zeta_k(x)\colon x\in\Lambda_k,  k=0,  \dots,n-1\}$ via \eqref{E:2.26}.
Under the second part of the product law on the right, the law of~$\varphi_{n-1}$ is that of a GFF on~$\Lambda_{n-1}$. Integrating over the variables $\{\zeta_k(x)\colon x\in\Lambda_k, \,  k=0,  \dots,n-1\}$ with~$\varphi_{n-1}$ fixed thus reduces this to
\begin{equation}
Z_{n-1}^{-1}\texte^{-\frac\beta2\sum_{x\in\Lambda_n}\zeta_n(x)^2}\texte^{\frac\beta2(\varphi_{n-1},\Delta_{n-1}\varphi_{n-1})_{n-1}}
\end{equation}
times the measure
\begin{equation}
\biggl[\,\bigotimes_{x\in\Lambda_n}\#\bigl(\varphi_{n-1}(m(x))+\textd\zeta_n(x)\bigr)\biggr]\otimes\biggl[\,\bigotimes_{z\in\Lambda_{n-1}}\textd\varphi_{n-1}(z)\biggr],
\end{equation}
where~$Z_{n-1}$ is a suitable normalization constant. Substituting~$\varphi_n(x):=\varphi_{n-1}(m(x))+\zeta_n(x)$ and performing the integrals over~$\varphi_{n-1}$ with the help of \eqref{E:2.4a} then reduces this to
\begin{equation}
Z_{n-1}^{-1}\frakz_n^{-1}\texte^{\frac\beta2(\varphi_n,\Delta_n\varphi_n)_n}\bigotimes_{x\in\Lambda_n}\#\bigl(\textd\varphi_n(x)\bigr).
\end{equation}
Modulo the form of the normalization constant, this is the law~$P_{n,\beta}$ from \eqref{E:1.1a}.
\end{proofsect}

A key property of the law \eqref{E:2.25} is that the collection of random variables
\begin{equation}
\label{E:2.31a} 
\bigl\{\varphi_k(x)\colon x\in\Lambda_k,  k=0,  \dots,n\bigr\}
\end{equation}
 is tree-indexed Markovian since conditioning on $\{\varphi_j(x)\colon x\in\Lambda_j,  j=0,  \dots,k-1\}$ determines the law of~$\varphi_k$ in such a way that~$\varphi_k(x)$ depends on~$\varphi_{k-1}(m(x))$ only. The transition probabilities vary with~$k$ and so the chain is ``time''-inhomogeneous.
 
The observation we just made also permits us to think of the collection \eqref{E:2.31a} as a  Branching Markov Chain.  Indeed, the value $\varphi_k(x)$ at ``time''~$k$ splits into~$b$ values $\varphi_{k+1}(x^{(1)}),\dots,\varphi_{k+1}(x^{(b)})$ at ``time'' $k+1$ (where $x^{(1)},\dots,x^{(b)}$ denotes the ``children'' of~$x$ in~$\Lambda_{k+1}$) that are independent samples from a law that depends only on~$\varphi_k(x)$. (In the absence of this dependence, the process would be a Branching Random Walk.)

\subsection{Control of effective potentials}
\label{sec-2.4}\noindent
In order to make a good use of the representation \eqref{E:2.26}, we need to control  the ``evolution'' of the effective potentials~$v_0,v_1,\dots$ from \eqref{E:2.2}. Although we are primarily interested in~$v_0$ arising for the DG model (see~\eqref{E:2.21}), other choices of the ``initial value'' may be considered as well. One of these is the Sine-Gordon model where
\begin{equation}
\label{E:2.31}
 v_0^{\mathrm{SG}}(z)  :=-\kappa\cos(2\pi z),\quad z\in\R.
\end{equation}
This case is interesting because, by varying $\kappa$ over the positive reals, it allows us to smoothly interpolate between the GFF ($\kappa=0$) and the DG model ($\kappa=+\infty$). 

The available analyses of the iterations \eqref{E:2.2}  seem to be invariably based  on linearization (see, e.g., Bauerschmidt and Bodineau~\cite{BB}). With the eyes on the ``trivial'' fixed point, this amounts to replacement of the map~$v_k\mapsto v_{k+1}$ by its ``infinitesimal'' version represented by the linear operator
\begin{equation}
\label{E:2.4}
\LL f(z):=b E\bigl(\, f(z+\zeta)\bigr)
\end{equation}
that can now be studied using methods of linear algebra.
Relying on the 1-periodicity dictated by the setting of the problem, this is achieved by passing to the Fourier representation $f(z) = \sum_{n\in\Z}\hat f(n)\texte^{2\pi\texti nz}$ as
\begin{equation}
\label{E:2.5}
\LL f(z) = \sum_{n\in\Z}\hat f(n)\,b\theta^{n^2}\texte^{2\pi\texti nz},
\end{equation}
where
\begin{equation}
\label{E:2.6}
\theta:=\texte^{-2\pi^2/\beta}.
\end{equation}
Denoting the dual $L^1$-norm by $\Vert f\Vert^\ast_1:=\sum_{n\in\Z}|\hat f(n)|$,  we get 
\begin{equation}
\bigl\Vert\LL (f-\hat f(0))\bigr\Vert^\ast_1\le b\theta\bigl\Vert f-\hat f(0)\bigr\Vert^\ast_1
\end{equation}
showing that the linear map is contractive on mean-zero functions whenever~$b\theta<1$. In terms of~$\beta$, this regime translates into
\begin{equation}
0<\beta<\beta_\cc:=\frac{2\pi^2}{\log b}.
\end{equation}
As is readily checked from \eqref{E:2.5}, the map is not contractive when~$\beta>\beta_\cc$.

In order for the linearized flow \eqref{E:2.4} to serve as a good approximation of the non-linear iterations \eqref{E:2.2}, one needs to control the higher-order terms. A standard approach (used, e.g., in~\cite{BB}) is to absorb these into a small change of the first-order term. Unfortunately, this forces us to start the iterations the closer to the desired limit point the larger is the norm of~$\LL$, thus leading to a loss of uniformity. 

One way to overcome the lack of uniformity is by tracking the ``evolution'' of some of the higher order terms as well. This has its merits and has been made to work for, e.g., the lattice Sine-Gordon model at uniformly small~$\kappa$ (Falco~\cite{Falco1,Falco2}). But, even in this approach, one needs a ``small parameter'' to ensure that the remaining high-order terms can be dealt with perturbatively. Unless we are willing (as is sometimes done) to turn~$b$ into a variable, no such ``small parameter'' seems to be there for the DG model.

As a consequence, none of the existing approaches \cite{MP,GM,BR,BB} extend throughout the subcritical regime of the DG model. As our desire is to work all the way up to~$\beta_\cc$, we will proceed non-perturbatively from the outset. Namely, we will work directly with the non-linear problem \eqref{E:2.2} relying on Fourier representation and combinatorial arguments for the Fourier coefficients. Our analysis applies to a rather general class of ``initial values''~$v_0$; definitely, those conforming to:

\begin{definition}
Let~$\VV$ denote the set of even functions~$v\colon\R\to\R$ for which there exists a positive sequence $\{a(n)\}_{n\in\Z}\in\ell^1(\Z)$ such that
\begin{equation}
\label{E:2.38}
\texte^{-v(z)}=\sum_{n\in\Z}a(n)\texte^{2\pi\texti n z},\quad z\in\R,
\end{equation}
and
\begin{equation}
\sup_{n\ge0}\frac{a(n+1)}{a(n)}<\infty.
\end{equation}
\end{definition}

As we will see, the set~$\VV$ includes the function~$v_0$ from \eqref{E:2.21} for the DG model as well as that in \eqref{E:2.31} corresponding to the Sine-Gordon model. We note that the sequence $\{a(n)\}_{n\in\Z}$ has a meaning in its own right as it represents the \textit{a priori} charge intensity in the Coulomb gas model dual to the GFF with on-site potential~$v$. 

For ``initial values'' in~$\VV$ we then get:

\begin{theorem}[Subcritical regime]
\label{thm-3.4}
Let~$v_0\in\VV$ and suppose that $\theta:=\texte^{-2\pi^2/\beta}$ obeys $b\theta<1$.  Denote  $c_0:=\sup_{n\ge0}\frac{a_0(n+1)}{a_0(n)}$ for~$\{a_0(n)\}_{n\in\Z}$ related to~$v_0$ via \eqref{E:2.38}.  Then the iterations $\{v_k\}_{k\ge0}$ obtained from~$v_0$ via \eqref{E:2.2} satisfy 
\begin{equation}
\label{E:2.40}
\sup_{z,z'\in[0,1]}\bigl|\,v_{k+1}(z)-bv_k(z')\bigr|\le 8 (b\theta)^kbc_0
\end{equation}
 whenever  $k\in\N$ is such that  $(b\theta)^kbc_0\le 1/8$.
\end{theorem}

This theorem, which we will prove in Section~\ref{sec3}, will serve as a key input in all of the remaining derivations in this paper. The statement  focuses on  $(z,z')\mapsto v_{k+1}(z)-bv_k(z')$ as this is what  is relevant for  the law \eqref{E:1.23}. The function~$v_k$ does not actually converge by itself which can be attributed to the zero mode (i.e., $n=0$ term) in \eqref{E:2.5} being expansive.  See Remark~\ref{rem-4.4} for a more precise statement. 

In order to demonstrate the strength of our method, we will show that it is capable of handling  the critical case as well:

\begin{theorem}[Critical regime]
\label{thm-3.7}
Let~$v_0\in\VV$ and suppose $\theta:=\texte^{-2\pi^2/\beta}$ obeys $b\theta=1$.  Denote  $c_0:=\sup_{n\ge0}\frac{a_0(n+1)}{a_0(n)}$ for~$\{a_0(n)\}_{n\in\Z}$ related to~$v_0$ via \eqref{E:2.38}.  Then  there exists $\gamma>0$ such that the iterations $\{v_k\}_{k\ge0}$ obtained from~$v_0$ via \eqref{E:2.2} satisfy 
\begin{equation}
\label{E:2.10}
\sup_{z,z'\in[0,1]}\bigl|v_{k+1}(z)-bv_k(z')\bigr|\le 8[1+\gamma k]^{-1/2}bc_0
\end{equation}
 whenever $k\in\N$ is such that  $[1+\gamma k]^{-1/2}bc_0\le 1/8$.
\end{theorem}

The constant~$\gamma$ depends only on~$v_0$; see \eqref{E:3.35} for an explicit  formula. 
The fact that the bound on the decay in \eqref{E:2.10} is only polynomial reflects an actual phenomenon. Indeed, it is not hard to prove a corresponding lower bound as well.  It is easy to  check that, for~$\beta>\beta_\cc$, the difference $v_{k+1}(z)-b v_k(z')$ does not converge to zero as~$k\to\infty$; that is, unless~$v_0$ is $1/m$-periodic for some natural~$m\ge2$.  As we show in~\cite{BiH} by rather involved computations,  for $\beta>\beta_\cc$ with~$\beta-\beta_\cc$ sufficiently small the difference converges to a non-constant continuous function and the convergence is exponentially fast.  

\subsection{Coupling}
With the iterations under control throughout the subcritical regime, we now address the second key ingredient of our proofs. For easier formulations later, we introduce the following notation. Given~$n\ge1$ and a sample $\{\zeta_k(x)\colon x\in\Lambda_{k},\, k=0,  \dots,n\}$ from the measure \eqref{E:1.23}, for each~$x\in\Lambda_n$ and~$k=0,\dots,n$ denote
\begin{equation}
\label{E:2.42}
\xi_k^{\DG}(x):=  \zeta_{n-k}\bigl(m^k(x)\bigr).
\end{equation}
By Lemma~\ref{lemma-2.2}, the field
\begin{equation}
\label{E:2.43}
\varphi^{\DG}_x:=  \sum_{k=0}^{n}  \xi_k^{\DG}(x)
\end{equation}
on~$\Lambda_n$ is then distributed as the DG model on~$\Lambda_n$.

Similarly, given $n\ge1$ and a sample $\{\zeta_k'(x)\colon x\in\Lambda_{k},\, k=0,  \dots,n\}$ of i.i.d.~$\NN(0,1/\beta)$, for each~$x\in\Lambda_n$ and~$k=0,\dots ,n $ denote
\begin{equation}
\label{E:2.44}
\xi_k^{\GFF}(x):=\zeta_{n-k}'\bigl(m^k(x)\bigr).
\end{equation}
The field $\varphi^{\GFF}$ on~$\Lambda_n$ defined by
\begin{equation}
\label{E:2.45}
\varphi^{\GFF}_x:= \sum_{k=0}^{n}  \xi_k^{\GFF}(x)
\end{equation}
is then distributed as the GFF on~$\Lambda_n$ at inverse temperature~$\beta$.
Our statement of the coupling of these two fields is then as follows:

\begin{theorem}[Coupling]
    \label{thm-4.1}
For all $\beta \in(0,\beta_\cc)$ there exist a constant $C >0$, a positive sequence $\{R_k\}_{k\ge1}$ with $\limsup_{k\to\infty}k^{-1}\log R_k<0$ and, for all natural~$n\ge1$, there exists a coupling of
\begin{align}
\label{E:4.1b}
&\bigl\{\xi_k^{\DG}(x)\colon k=0,\dots,n ,\, x\in\Lambda_n\bigr\},
\\
&
\label{E:4.2b}
\bigl\{\xi_k^{\GFF}(x)\colon k=0,\dots,n,\, x\in\Lambda_n\bigr\},
\end{align}
and a family of zero-one valued random variables
\begin{equation}
\label{E:4.3}
\bigl\{B_k(x)\colon x\in\Lambda_{n-k},\,k=1,\dots,n\bigr\}
\end{equation}
such that the following holds:
\begin{enumerate}
\item[(1)] the random variables in \eqref{E:4.3} are independent with
\begin{equation}
\label{E:3.49u}
P\bigl(B_k(x)=1\bigr) = \texte^{-R_k},\quad x\in\Lambda_{n-k},\,k=1,\dots,n,
\end{equation}
\item[(2)] the families \eqref{E:4.2b} and \eqref{E:4.3} are independent of each other,
\item[(3)] for all $k=1,\dots,n$ and all~$x\in\Lambda_n$,
\begin{equation}
\label{E:4.2a}
\xi_k^{\DG}(x)=\xi_k^{\GFF}(x) \quad\text{\rm on }\bigl\{ B_k(m^{k}(x)) =1\bigr\},
\end{equation}
\item[(4)] for all $k=0,\dots,n$ and all~$x\in\Lambda_n$,
\begin{equation}
        \label{E:4.1}
         P\Bigl(\bigl| \xi_k^{\DG}(x) - \xi_k^{\GFF} (x)\bigr|> C\Bigr)=0.
\end{equation}
\end{enumerate}
\end{theorem}

That a joint coupling of the increments exists is perhaps not surprising due to their Markovian structure and the relative closeness of \eqref{E:2.25} to the product normal law. Considerably less obvious is the fact that the increments can be kept within a uniformly bounded distance of each other; cf \eqref{E:4.1}. This is quite useful in estimates of the difference of the sums \eqref{E:2.43} and \eqref{E:2.45} which can thus be controlled by the sheer number of places where equality on the left of \eqref{E:4.2a} fails. In light of \eqref{E:4.2a}, this number is dominated by the number of Bernoulli variables in \eqref{E:4.3} that equal zero.

The sequence $\{R_k\}_{k\ge1}$ governing the errors in \eqref{E:3.49u} is closely related to the error in~\eqref{E:2.40}; see \eqref{E:4.48}. An exponential decay is key but the rate of that decay is not, which is another reason why our approach works all the way up to~$\beta_\cc$.  Theorem~\ref{thm-4.1} will serve as the main input in our control of the extremal behavior; besides this the proofs require only  some standard facts from analysis of the GFF extrema and technical estimates. The proof of Theorem~\ref{thm-4.1} comes at the end of Section~\ref{sec4}.


\section{Renormalization group flow}
\label{sec3}\noindent
The formal presentation of our proofs opens up by careful control of the iterations \eqref{E:2.2} resulting in proofs of Theorem~\ref{thm-3.4} and~\ref{thm-3.7}. We start with the former theorem as it deals with the regime that the rest of this paper is focused on.

\subsection{Subcritical regime}
The argument will be based on Fourier techniques that are enabled by $1$-periodicity of the relevant functions. A starting point is the following representation:

\begin{lemma}
Suppose $a_0=\{a_0(n)\}_{n\in\Z}\in\ell^1(\Z)$ is such that~$v_0$ obeys
\begin{equation}
\label{E:4.1ii}
\texte^{-v_0(z)}=\sum_{n\in\Z}a_0(n)\texte^{2\pi\texti n z},\quad z\in\R.
\end{equation}
Then the iterates $\{v_k\}_{k\in\N}$ defined from~$v_0$ via \eqref{E:2.2} admit the representation
\begin{equation}
\label{E:4.2ii}
\texte^{-v_k(z)}=\sum_{n\in\Z}a_k(n)\texte^{2\pi\texti n z} ,\quad z\in\R,
\end{equation}
with  $a_k=\{a_k(n)\}_{n\in\Z}\in\ell^1(\Z)$, for each~$k\in\N$, where the $a_k$'s are defined recursively by 
\begin{equation}
\label{E:3.8}
a_{k+1}(n) := \sum_{\begin{subarray}{c}
\ell_1,\dots,\ell_b\in\Z\\\ell_1+\dots+\ell_b=n
\end{subarray}}
\biggl[\,\prod_{i=1}^b a_k(\ell_i)\biggr]\,\theta^{n^2},\quad n\in\Z,
\end{equation}
for~$\theta$ as in \eqref{E:2.6}.
\end{lemma}

\begin{proofsect}{Proof}
The fact that the $a_k$'s obey \eqref{E:3.8} is checked directly via \eqref{E:2.2} provided the sums in \eqref{E:3.8} converge absolutely. To see the latter, note that \eqref{E:3.8} implies $a_{k+1}(n)\le\Vert a_k\Vert_1^b \theta^{n^2}$ which in light of~$\theta<1$ shows that $a_k\in\ell^1(\Z)$ implies $a_{k+1}\in\ell^1(\Z)$ for all~$k\ge1$. The summability assumption imposed on~$a_0$ thus propagates through the iterations. 
\end{proofsect}

We note that for the initial  effective potential~$v_0$ of the DG model \eqref{E:2.21} gives
\begin{equation}
\label{E:3.9}
 a_0^{\DG}(n) =\theta^{n^2} \sqrt{\frac{2\pi}{\beta}},\quad n\in\Z,
\end{equation}
while 
for the Sine-Gordon model \eqref{E:2.31} we obtain
\begin{equation}
 a_0^{\mathrm{SG}}(n) =\sum_{\ell=0}^\infty\frac{(\kappa/2)^{2\ell+|n|}}{(\ell+|n|)!\ell!},\quad n\in\Z.
\end{equation}
Both of these lie in~$\VV$ as desired. 

We will assume~$v_0\in\VV$ throughout the rest of this section  and write~$\{a_k(n)\}_{n\in\Z}$ for the  coefficients so that \twoeqref{E:4.1ii}{E:4.2ii} holds.  Observe that, as is checked from \eqref{E:2.2}, all the~$v_k$'s are real-valued and even and so
\begin{equation}
0<a_k(-n)=a_k(n)\le a_k(0),\quad n,k\ge0,
\end{equation}
where the strict positivity of~$a_k$ is implied by the strict positivity of~$a_0$ and \eqref{E:3.8}.
The key estimate driving our control of the sequence $\{v_k\}_{k\ge0}$ now comes in:

\begin{lemma}
For all naturals~$k\ge0$ and~$n\ge0$,
\begin{equation}
\label{E:3.11}
\frac{a_{k+1}(n+1)}{a_{k+1}(n)}\le b\,\theta^{(n+1)^2-n^2}\sup_{\ell\ge0}\frac{a_k(\ell+1)}{a_k(\ell)}.
\end{equation}
In particular, $v_0\in\VV$ implies~$v_k\in\VV$ for all~$k\ge1$.
\end{lemma}

\begin{proofsect}{Proof}
The proof is based on a combinatorial argument that is reminiscent of the Peierls argument from statistical mechanics.
Denote
\begin{equation}
\Xi(n):=\bigl\{(\ell_1,\dots,\ell_b)\in\Z^b\colon \ell_1+\dots+\ell_b=n\bigr\}
\end{equation}
and, for~$n\ge0$, construct a map $\phi\colon\Xi(n+1)\to\Xi(n)$ by setting
\begin{equation}
\label{E:3.9x}
\phi\bigl((\ell_1,\dots,\ell_b)\bigr):=(\ell_1,\dots,\ell_{i-1},\ell_i-1,\ell_{i+1},\dots,\ell_b),
\end{equation}
where
\begin{equation}
\label{E:3.13}
i=i(\ell_1,\dots,\ell_b):=\min\{j=1,\dots , b\colon \ell_j>0\}.
\end{equation}
(This is well defined because $n+1>0$ and so, by the pigeon-hole principle, the set on the right is non-empty.)
Note that, since the map $(i(\cdot),\phi(\cdot))\colon\Xi(n+1)\to \{1,\dots,b\} \times\Xi(n)$ is injective, its projection to~$\phi$ is at most $b$-fold degenerate and so
\begin{equation}
\label{E:3.14}
\bigl|\,\phi^{-1}(\{(\ell_1,\dots,\ell_b)\})\bigr|\le b,\quad  (\ell_1,\dots,\ell_b)\in\Xi(n).
\end{equation}
Denoting the supremum on the right of \eqref{E:3.11} by~$c_k$, observe also that
\begin{equation}
\label{E:3.15}
\prod_{i=1}^b a_k(\ell_i)
\le c_k\prod_{i=1}^b a_k(\ell_i')
\end{equation}
holds for each~$(\ell_1',\dots,\ell_b')\in\phi(\Xi(n+1))$ and each $(\ell_1,\dots,\ell_b)\in \phi^{-1}(\{(\ell_1',\dots,\ell_b')\})$.

Abbreviating $\bar\ell:= (\ell_1,\dots,\ell_b)$ and $\bar\ell':=(\ell_1',\dots,\ell_b')$, the inductive definition \eqref{E:3.8} of~$a_{k+1}$ from~$a_k$ now shows
\begin{equation}
\label{E:3.16}
\begin{aligned}
a_{k+1}(n+1)
&
=\sum_{\bar\ell'\in\phi(\Xi(n+1))}\,
\sum_{\bar\ell\in\phi^{-1}(\{\bar\ell'\})}
\biggl[\,\prod_{i=1}^b a_k(\ell_i)\biggr]\,\theta^{(n+1)^2}
\\
&\le c_k \sum_{\bar\ell'\in\phi(\Xi(n+1))}\,
\sum_{\bar\ell\in\phi^{-1}(\{\bar\ell'\})}
\biggl[\,\prod_{i=1}^b a_k(\ell_i')\biggr]\,\theta^{(n+1)^2}
\\
&\le c_k b \theta^{(n+1)^2-n^2}\sum_{\bar\ell'\in\phi(\Xi(n+1))}
\biggl[\,\prod_{i=1}^b a_k(\ell_i')\biggr]\,\theta^{n^2},
\end{aligned}
\end{equation}
where the first inequality follows from \eqref{E:3.15} and the second inequality from \eqref{E:3.14}.
By~\eqref{E:3.8} again, the sum on the right is at most~$a_{k+1}(n)$.
\end{proofsect}

Iterations of \eqref{E:3.11} yield: 
\begin{lemma}
Suppose~$b\theta\le1$ and let~$c_0:=\sup_{n\ge0}\frac{a_0(n+1)}{a_0(n)}$. Then 
\begin{equation}
\label{E:3.17}
\sup_{\ell\ge0}\frac{a_k(\ell+1)}{a_k(\ell)}\le (b\theta)^kc_0,\quad k\ge0,
\end{equation}
and, consequently,
\begin{equation}
\label{E:3.18}
a_k(n)\le\bigl[(b\theta)^{k-1}bc_0\bigr]^n \theta^{n^2} a_k(0),\quad n,k\ge1.
\end{equation}
\end{lemma}

\begin{proofsect}{Proof}
Continuing to write~$c_k$ for the supremum in \eqref{E:3.11}, the fact that $(n+1)^2-n^2\ge1$ for all~$n\ge0$ gives $c_{k+1}\le(b\theta)c_k$. This shows $c_k\le(b\theta)^kc_0$, proving \eqref{E:3.17}. Plugging this in \eqref{E:3.11} yields
\begin{equation}
a_k(n+1)\le\theta^{(n+1)^2-n^2}(b\theta)^{k-1} bc_0 \,a_k(n),\quad k\ge1,\,n\ge0,
\end{equation}
and so, by induction, we get \eqref{E:3.18}. 
\end{proofsect}

Returning to the original setting of the problem, we conclude:

\begin{proofsect}{Proof of Theorem~\ref{thm-3.4}}
The definition of the coefficients shows 
\begin{equation}
\label{E:3.21}
\texte^{-v_{k+1}(z)} = \sum_{n\in\Z} a_{k+1}(n)\texte^{2\pi\texti n z},
\end{equation}
while \eqref{E:3.8} gives 
\begin{equation}
\label{E:3.22}
\texte^{-bv_k(z')} = \sum_{n\in\Z} a_{k+1}(n)\theta^{-n^2}\texte^{2\pi\texti n z'}.
\end{equation}
The main point is to show that these sums are dominated by the~$n=0$ term for~$k$ large. To this end we write \eqref{E:3.21} as $a_{k+1}(0)[1+e_1]$ and \eqref{E:3.22} as $a_{k+1}(0)[1+e_2]$. Assuming~$k$ to be so large that $\eta:= (b\theta)^{k}bc_0<1$, \eqref{E:3.18} shows
\begin{equation}
|e_1|\le 2\sum_{n\ge1}\eta^n\theta^{n^2} \le 2\theta\frac{\eta}{1-\eta}
\end{equation}
while 
\begin{equation}
|e_2|\le 2\sum_{n\ge1}\eta^n \le 2\frac{\eta}{1-\eta}.
\end{equation}
Now
\begin{equation}
\label{E:3.20}
v_{k+1}(z)-bv_k(z')=\log\frac{1+e_2}{1+e_1} = \log\Bigl[1+\frac{e_2-e_1}{1+e_1}\Bigr]
\end{equation}
and
\begin{equation}
\label{E:3.26}
\Bigl|\frac{e_2-e_1}{1+e_1}\Bigr|\le\frac{|e_1|+|e_2|}{1-|e_1|}.
\end{equation}
Under the (stronger) assumption that $\eta\le1/8$ we then have $|e_2|\le 2/7$ and $|e_1|\le(2/7)\theta$ which by the fact that~$\theta\le b^{-1}\le1/2$ gives $|e_1|\le 1/7$. Then the right-hand side of \eqref{E:3.26} is less than~$1/2$ and so we may use that $|\log(1+x)|\le 2|x|$ holds for all~$|x|\le1/2$. A calculation (still under $\eta\le1/8$) shows
\begin{equation}
\frac{|e_1|+|e_2|}{1-|e_1|}\le4\eta
\end{equation}
and so the right-hand side of \eqref{E:3.20} is at most~$8\eta$.
\end{proofsect}

\begin{remark}
\label{rem-4.4}
 Note that the proof shows that
\begin{equation}
\sup_{z\in\R}\,\Bigl|\,v_k(z) +b^{-1}\log a_{k+1}(0)\Bigr|\,\,\underset{k\to\infty}\longrightarrow\,\,0
\end{equation}
with the convergence exponentially fast.  The superexponential decay of $n\mapsto a_k(n)$ indicates that all derivatives of $z\mapsto v_k(z)$, and thus also of $(z,z')\mapsto v_{k+1}(z)-b v_k(z')$, tend to zero exponentially fast. 
\end{remark}


\subsection{Critical regime}
Our next item of concern is the critical regime $b\theta=1$. Here the sequence of coefficients still concentrates to a ``point mass'' at~$n=0$ or, more precisely, $(z,z')\mapsto v_{k+1}(z)-bv_k(z')$ tends to zero as $k\to\infty$, but the convergence is no longer exponentially fast. Notwithstanding, with suitable refinements, this can still be proved using the same ideas that drove the proofs in the subcritical regime. The key is to control the ratio $a_k(1)/a_k(0)$ separately from the other ratios $a_k(\ell+1)/a_k(\ell)$ with $\ell\ge1$:

\begin{lemma}
For each~$k\ge0$,
\begin{equation}
\label{E:3.28}
\frac{a_{k+1}(1)}{a_{k+1}(0)}\le\theta\frac{\displaystyle\frac{a_k(1)}{a_k(0)}}{\displaystyle 1+\binom b2\Bigl(\frac{a_k(1)}{a_k(0)}\Bigr)^2}+(b-1)\theta\sup_{\ell\ge0}\frac{a_k(\ell+1)}{a_k(\ell)}.
\end{equation}
\end{lemma}

\begin{proofsect}{Proof}
We proceed similarly as in \eqref{E:3.16} albeit with additional care paid to certain particular terms. Recall the definition of the map~$\phi$ from \eqref{E:3.9x}. Denoting $\bar 0:=(0,\dots,0)$, each~$\bar\ell\in\phi^{-1}(\{\bar 0\})$ has exactly one index equal to one and the others equal to zero. This shows that the contribution of~$\bar\ell':=\bar 0$ to~$a_{k+1}(1)$ on the first line of \eqref{E:3.16} equals
\begin{equation}
\label{E:3.25u}
\sum_{\bar\ell\in \phi^{-1}(\{\bar 0\})}\,
\biggl[\,\prod_{i=1}^b a_k(\ell_i)\biggr]\,\theta
=b\theta a_k(1)a_k(0)^{b-1}.
\end{equation}
Next consider~$\bar\ell\in\Xi(0)\smallsetminus\{\bar 0\}$ and observe that each~$\bar\ell'\in\phi^{-1}(\{\bar\ell\})$ has at least two non-zero entries with total sum of the entries equal to one. This forces at least one of the entries to be negative; i.e., $A:=\{i=1,\dots,b\colon\ell'_i<0\}\ne\emptyset$. But the definition of~$\phi$ implies~$\ell_i=\ell_i'$ for each~$i\in A$ and, since~$\bar\ell$ and~$\bar\ell'$ differ only at one index, say~$i$, where $\ell_i'=\ell_i+1$, we have
\begin{equation}
\bigl|\phi^{-1}(\{\bar\ell\})\bigr|\le b-1,\quad \bar\ell\in\Xi(0)\smallsetminus\{\bar 0\}.
\end{equation}
The calculation \twoeqref{E:3.15}{E:3.16} then shows
\begin{equation}
\begin{aligned}
\label{E:3.27u}
\sum_{\bar\ell'\in\phi(\Xi(1))\smallsetminus\{\bar0\}}\,
&\sum_{\bar\ell\in\phi^{-1}(\{\bar\ell'\})}
\biggl[\,\prod_{i=1}^b a_k(\ell_i)\biggr]\,\theta
\\
&\le c_k\theta\sum_{\bar\ell'\in\phi(\Xi(1))\smallsetminus\{\bar0\}}\,
\sum_{\bar\ell\in\phi^{-1}(\{\bar\ell'\})}
\biggl[\,\prod_{i=1}^b a_k(\ell_i')\biggr]
\\
&\le c_k(b-1)\theta\sum_{\bar\ell'\in \Xi(0)\smallsetminus\{\bar0\}}
\biggl[\,\prod_{i=1}^b a_k(\ell_i')\biggr]
\le
c_k (b-1)\theta \bigl[a_{k+1}(0)-a_k(0)^b\bigr],
\end{aligned}
\end{equation}
where~$c_k$ again denotes the supremum in \eqref{E:3.11} and where \eqref{E:3.8} was used in the last step.

Combining \eqref{E:3.25u} and \eqref{E:3.27u} along with the fact that~$c_k a_k(0)\ge a_k(1)$ yields
\begin{equation}
\label{E:3.28u}
\begin{aligned}
a_{k+1}(1)
&\le b\theta a_k(1)a_k(0)^{b-1} + c_k (b-1)\theta \bigl[a_{k+1}(0)-a_k(0)^b\bigr]
\\
&= b\theta a_k(1)a_k(0)^{b-1} - (b-1)\theta c_k a_k(0) a_k(0)^{b-1} +c_k (b-1) \theta a_{k+1}(0)
\\
&\le \theta a_k(1)a_k(0)^{b-1}+c_k (b-1) \theta a_{k+1}(0).
\end{aligned}
\end{equation}
Restricting the sum in \eqref{E:3.8} for~$n:=0$ to~$b$-tuples of indices that are either all zeros or contain one one and one negative one and are zero otherwise shows
\begin{equation}
a_{k+1}(0)\ge a_k(0)^b +\binom b2 a_k(0)^{b-2}a_k(1)^2.
\end{equation}
(We could write $b(b-1)$ instead of~$\binom b2$ but the latter is more compact so we keep that.) Using this jointly with \eqref{E:3.28u} then gives \eqref{E:3.28}.
\end{proofsect}

The recursive inequalities \eqref{E:3.28} and \eqref{E:3.11} do not have the same structure but we can put them together with the help of a suitable parameter. This yields:

\begin{lemma}
Suppose~$b\theta=1$ and write~$c_k$ for the supremum in \eqref{E:3.28}. For each~$k\ge0$, define $\alpha_k\in(0,1)$ to be the unique number such that
\begin{equation}
\label{E:3.33}
\alpha_k = \frac1{ 1+\binom b2 c_k^2\alpha_k^2}.
\end{equation}
Then
\begin{equation}
\label{E:3.34}
c_{k+1}\le\frac{b^{-1}\,c_k}{1+\binom b2 c_k^2\alpha_k^2}+(1-b^{-1}) c_k.
\end{equation}
The sequence $\{c_k\}_{k\ge0}$ is non-increasing while~$\{\alpha_k\}_{k\ge0}$ is non-decreasing.
\end{lemma}

\begin{proofsect}{Proof}
The proof is based on the fact that the definition of~$\alpha_k$ gives
\begin{equation}
\label{E:3.37}
(1-b^{-1}+b^{-1}\alpha_k)c_k = \frac{b^{-1}c_k}{ 1+\binom b2 c_k^2\alpha_k^2}+(1-b^{-1})c_k.
\end{equation}
It thus suffices to dominate the ratio~$a_{k+1}(n+1)/a_{k+1}(n)$ by either side of this equality for all~$n\ge0$.

We start with a bound on $a_{k+1}(1)/a_{k+1}(0)$. In the situation when $a_k(1)/a_k(0)\le\alpha_k c_k$, dropping the second term in the denominator in  \eqref{E:3.28} while bounding the numerator by~$\alpha_k c_k$ shows
\begin{equation}
\frac{a_{k+1}(1)}{a_{k+1}(0)}\le b^{-1}\alpha_k c_k+(1-b^{-1}) c_k = (1-b^{-1}+b^{-1}\alpha_k)c_k.
\end{equation}
If in turn $a_k(1)/a_k(0)\ge\alpha_k c_k$, then we  bound  the denominator of \eqref{E:3.28} by $1+\binom b2 c_k^2\alpha_k^2$ from below and dominate the ratio in the numerator by~$c_k$ to get
\begin{equation}
\frac{a_{k+1}(1)}{a_{k+1}(0)} \le\frac{b^{-1} c_k}{1+\binom b2 c_k^2\alpha_k^2}+(1-b^{-1}) c_k.
\end{equation}
Hence $a_{k+1}(1)/a_{k+1}(0)$ is bounded by the right-hand side of \eqref{E:3.34} in both cases.

Moving to the terms $a_{k+1}(n+1)/a_{k+1}(n)$ with $n\ge1$, here \eqref{E:3.11} gives
\begin{equation}
\label{E:3.33u}
\frac{a_{k+1}(n+1)}{a_{k+1}(n)}\le b^{n^2-(n+1)^2+1}c_k.
\end{equation}
For~$n\ge1$ the prefactor is at most $b^{-2}\le1/4$ while $(1-b^{-1}+b^{-1}\alpha_k)\ge 1-b^{-1}\ge\frac12$. Hence, this is less than the right-hand side of \eqref{E:3.34} in this case as well.

In order to prove the monotonicity statements, note that \eqref{E:3.11} along with $\theta<1$ and $(n+1)^2-n^2\ge1$ for~$n\ge0$ gives $c_{k+1}\le b\theta c_k$. For~$b\theta\le1$ this shows that~$k\mapsto c_{k}$ is non-increasing. Using this on the right of \eqref{E:3.33} then rules out that $\alpha_{k+1}<\alpha_k$, proving that $k\mapsto\alpha_k$ is non-decreasing.
\end{proofsect}

We are now ready to give:

\begin{proofsect}{Proof of Theorem~\ref{thm-3.7}}
Assume~$b\theta=1$. We start by processing the bound \eqref{E:3.34} a bit further. Putting the two terms on the right together we get
\begin{equation}
c_{k+1}\le c_k\frac{b^{-1}+(1-b^{-1})+(1-b^{-1}) \binom b2 c_k^2\alpha_k^2}{1+\binom b2 c_k^2\alpha_k^2}
\le c_k\frac{1+r\binom b2 c_k^2\alpha_k^2}{1+\binom b2 c_k^2\alpha_k^2},
\end{equation}
where we abbreviated~$r:=1-b^{-1}$.
Taking reciprocals and squaring gives
\begin{equation}
\frac1{c_{k+1}^2}\ge \frac1{c_k^2}\biggl[1+\frac{[1+\binom b2 c_k^2\alpha_k^2]^2-[1+r\binom b2 c_k^2\alpha_k^2]^2}{[1+r\binom b2 c_k^2\alpha_k^2]^2}\biggr]
\ge\frac1{c_k^2} +\frac{2\binom b2(1-r)\alpha_k^2}{[1+r\binom b2 c_k^2\alpha_k^2]^2}.
\end{equation}
Invoking~$1-r=b^{-1}$ along with~$\alpha_k^2\ge\alpha_0^2$ in the numerator while using that~$r\le1$ along with~$c_k^2\le c_0^2$ and~$\alpha_k^2\le1$ in the denominator then yields 
\begin{equation}
\label{E:3.35a}
\frac1{c_{k+1}^2}\ge \frac1{c_k^2}+\frac\gamma{c_0^2},\quad k\ge0,
\end{equation}
where
\begin{equation}
\label{E:3.35}
\gamma:=c_0^2\frac{2\binom b2 b^{-1}\alpha_0^2}{[1+\binom b2 c_0^2]^2}
\end{equation}
is a constant expressed only using the properties of~$a_0$.

Iterating \eqref{E:3.35a} shows $c_k^{-2}\ge c_0^{-2}+c_0^{-2}\gamma k = c_0^{-2}[1+\gamma k]$ and so
\begin{equation}
c_k\le c_0[1+\gamma k]^{-1/2},\quad k\ge0.
\end{equation}
Using this in \eqref{E:3.11}, for the Fourier coefficients we get
\begin{equation}
a_{k+1}(n+1)\le b^{n^2-(n+1)^2} bc_0 \bigl[1+\gamma k\bigr]^{-1/2}\,a_{k+1}(n),\quad n,k\ge0
\end{equation}
which implies
\begin{equation}
a_{k+1}(n)\le b^{-n^2} \bigl(bc_0 [1+\gamma k]^{-1/2}\bigr)^n\,a_{k+1}(0),\quad n,k\ge0.
\end{equation}
As an inspection of \eqref{E:3.18} reveals, the right-hand side of this bound has a very similar structure as for the subcritical case; indeed, $\theta^{n^2}$ is replaced by $b^{-n^2}$, which is its value at criticality, while the error term $(b\theta)^{k}bc_0$ is replaced by $[1+\gamma k]^{-1/2}bc_0$. Proceeding as in the proof of Theorem~\ref{thm-3.4}, this implies the claim.
\end{proofsect}

\begin{remark}
Also in this case the control extends to all derivatives.
\end{remark}


\section{Coupling}
\label{sec4}\noindent
Here we will construct the coupling between the DG model and the GFF stated in Theorem~\ref{thm-4.1}.
Due to the discrete nature of the DG model and, simultaneously, the continuum nature of the intermediate fields during iterations, the arguments split into two parts depending on whether we are dealing with continuum-valued or discrete-valued steps. 

\subsection{Continuum-valued steps}
Let $\Vert X\Vert_\infty$ denote the $L^\infty$-norm of real-valued random variable~$X$. We start with some general considerations pertaining to coupling of two random variables that admit a density with respect to a normal law.

\begin{proposition}
\label{prop-4.1}
Let~$\mu$ be the law of $\NN(0,\sigma^2)$ with~$\sigma\in(0,\infty)$ and let~$X$ and~$Y$ have laws
\begin{equation}
P\bigl(X\in A)=\int_A f_1(t)\mu(\textd t)\quad\text{and}\quad P\bigl(Y\in A)=\int_Af_2(t)\mu(\textd t)
\end{equation}
for some measurable~$f_1,f_2\colon\R\to(0,\infty)$.
Denote the associated CDFs by $F(t):=P(X\le t)$ and $G(t):=P(Y\le t)$ and let $h\colon\R\to\R$ be defined by
\begin{equation}
\label{E:4.4}
h(t):=G^{-1}\bigl(F(t)\bigr)
\end{equation}
Then
\begin{equation}
h(X)\laweq Y
\end{equation}
and, if~$f_1$, $1/f_1$, $f_2$ and~$1/f_2$ are bounded, then
\begin{equation}
\label{E:4.6}
\bigl\Vert h(X)-X\bigr\Vert_\infty<3\sigma\sqrt{2\log\Bigl(\,\max\bigl\{\Vert f_1\Vert_\infty\Vert 1/f_2\Vert_\infty,\,\Vert f_2\Vert_\infty\Vert 1/f_1\Vert_\infty\bigr\}\Bigr)}.
\end{equation}
The pair $(X,h(X))$ thus defines a coupling of~$X$ and~$Y$ with $X-Y$ bounded in~$L^\infty$.
\end{proposition}

The assumption  that~$f_1, f_2>0$ implies  that~ $F$ and~$G$ are  strictly increasing and so~ $F^{-1}$ and $G^{-1}$ are proper inverses.  A calculation shows
\begin{equation}
\begin{aligned}
P\bigl(h(X)\le t\bigr)&=P\bigl(F(X)\le G(t)\bigr)
\\
&= P\bigl(X\le F^{-1}(G(t))\bigr)
=F\bigl(F^{-1}(G(t))\bigr) =  G(t)
\end{aligned}
\end{equation}
and so~$h(X)$ has the law of~$Y$. 
It remains to prove the bound \eqref{E:4.6} which we will do by proving that $h(t)-t$ is bounded. The argument comes in two lemmas dealing with various separate cases.

\begin{lemma}
\label{lemma-4.3}
Suppose that~$1/f_1$ and~$f_2$ are bounded. Then
\begin{equation}
\label{E:4.8}
h(t)\le t+\frac2t\,\sigma^2\log\bigl(\,\Vert f_2\Vert_\infty\Vert 1/f_1\Vert_\infty\bigr),\quad t>0,
\end{equation}
and
\begin{equation}
\label{E:4.9}
h(t)\ge t-\frac2{|t|}\,\sigma^2\log\bigl(\,\Vert f_2\Vert_\infty\Vert 1/f_1\Vert_\infty\bigr),\quad t<0.
\end{equation}
\end{lemma}

\begin{proofsect}{Proof}
First observe that $\tilde X := X/\sigma$ and $\tilde Y:=Y/\sigma$  have densities  $\sigma f_1(\sigma \cdot)$ and $\sigma f_2(\sigma \cdot)$ with respect to the law of~$\NN(0,1)$. This means that, once we prove the above inequalities in the case $\sigma^2 = 1$, the case of general $\sigma^2$ follows the fact that $\tilde h(t):= G_{\tilde Y}^{-1}(F_{\tilde X}(t)) = \frac{1}{\sigma} h(\sigma t)$. So henceforth we assume that~$\sigma^2=1$.

Writing~$g$ for the probability density of~$\NN(0,1)$, the definition \eqref{E:4.4} translates into
\begin{equation}
\int_{-\infty}^{t}f_1(x)g(x)\textd x = \int_{-\infty}^{h(t)} f_2(x)g(x)\textd x.
\end{equation}
This shows that~$h$ is differentiable and obeys the ODE
\begin{equation}
f_1(t)g(t) = h'(t)f_2\bigl(h(t)\bigr) g\bigl(h(t)\bigr),
\end{equation}
which takes the explicit form
\begin{equation}
h'(t)=\frac{f_1(t)}{f_2(h(t))}\texte^{\frac12(h(t)-t)(h(t)+t)}.
\end{equation}
Suppose now that~$t>0$ and~$\epsilon>0$ are such that
\begin{equation}
\label{E:4.13}
h(t)\ge t+\epsilon\,\,\wedge\,\,\bigl(\,\Vert f_2\Vert_\infty\Vert 1/f_1\Vert_\infty\bigr)^{-1}\texte^{\epsilon t/2}>1.
\end{equation}
The first inequality implies
\begin{equation}
\frac{f_1(t)}{f_2(h(t))}\texte^{(h(t)-t)(h(t)+t)/2} \ge \bigl(\,\Vert f_2\Vert_\infty\Vert 1/f_1\Vert_\infty\bigr)^{-1}\texte^{\epsilon t/2}\texte^{\epsilon h(t)/2}
\end{equation}
and the second shows that the right-hand side exceeds one. This in turn shows that $u\mapsto h(u)-u$ is non-decreasing for $u-t$ small positive. As this reinforces the conditions in \eqref{E:4.13}, we readily infer that
\begin{equation}
\text{\eqref{E:4.13}}\,\,
\Rightarrow\,\, 
h(u)\ge u+\epsilon\,\text{ when }\,u\ge t.
\end{equation}
Abbreviating $c:=(\Vert f_2\Vert_\infty\Vert 1/f_1\Vert_\infty)^{-1}$, this shows
\begin{equation}
h'(u)\texte^{-\epsilon h(u)/2}\ge c\texte^{\epsilon u/2},\quad u\ge t.
\end{equation}
 This integrates to
\begin{equation}
\texte^{-\epsilon h(t)/2}-\texte^{-\epsilon h(u)/2}\ge c\bigl(\texte^{\epsilon u/2}-\texte^{\epsilon t/2}\bigr),\quad u\ge t.
\end{equation}
But this is absurd because the left-hand side is bounded while the right-hand side diverges as~$u\to\infty$. No~$t>0$ and $\epsilon>0$ satisfying \eqref{E:4.13} may thus exist and so $h(t) <  t+\epsilon$ once $c\texte^{\epsilon t/2}>1$. Solving for~$\epsilon$ yields
\eqref{E:4.8}.

Suppose now that~$t<0$ and~$\epsilon>0$ are such that
\begin{equation}
\label{E:4.18}
h(t)\le t-\epsilon\,\,\wedge\,\,\bigl(\,\Vert f_2\Vert_\infty\Vert 1/f_1\Vert_\infty\bigr)^{-1}\texte^{-\epsilon t/2}>1.
\end{equation}
Then
\begin{equation}
\frac{f_1(t)}{f_2(h(t))}\texte^{(h(t)-t)(h(t)+t)/2} \ge \bigl(\,\Vert f_2\Vert_\infty\Vert 1/f_1\Vert_\infty\bigr)^{-1}\texte^{-\epsilon t/2}\texte^{-\epsilon h(t)/2}.
\end{equation}
This again forces $h(u) \le u - \epsilon$ as well as $h'(u)>1$ for all~$u\le t$ and thus proves
\begin{equation}
h'(u)\texte^{\epsilon h(u)/2}\ge c\texte^{-\epsilon u/2},\quad u\le t,
\end{equation}
where~$c$ is as above. Integrating we get
\begin{equation}
\texte^{\epsilon h(t)/2}-\texte^{\epsilon h(u)/2}\ge c\bigl(\texte^{-\epsilon u/2}-\texte^{-\epsilon t/2}\bigr),\quad u\le t,
\end{equation}
which is again absurd when $u\to -\infty$. Hence we conclude \eqref{E:4.9} as well.
\end{proofsect}

Next we need to control the difference $h(t)-t$ in the regimes not covered by Lem\-ma~\ref{lemma-4.3}. Explicitly, we need an upper bound when~$t<0$ and a lower bound when $t>0$. This is the content of:

\begin{lemma}
\label{lemma-4.4}
Suppose that~$1/f_2$ and~$f_1$ are bounded. Then
\begin{equation}
\label{E:4.22}
0\le h(t)\,\,\Rightarrow\,\, h(t)\ge t-\frac2t\,\sigma^2\log\bigl(\,\Vert f_1\Vert_\infty\Vert 1/f_2\Vert_\infty\bigr),\quad t>0,
\end{equation}
and
\begin{equation}
\label{E:4.23}
h(t)\le 0\,\,\Rightarrow\,\, h(t)\le t+\frac2{|t|}\,\sigma^2\log\bigl(\,\Vert f_1\Vert_\infty\Vert 1/f_2\Vert_\infty\bigr),\quad t<0.
\end{equation}
\end{lemma}

\begin{proofsect}{Proof}
We again assume without loss of generality that $\sigma^2 = 1$. Suppose that, for some~$\epsilon>0$ and~$t>\epsilon$, 
\begin{equation}
\label{E:4.24}
0\le h(t)\le t-\epsilon\,\,\wedge\,\,\Vert f_1\Vert_\infty\Vert 1/f_2\Vert_\infty\,\texte^{-\epsilon t/2}<1.
\end{equation}
Then 
\begin{equation}
\frac{f_1(t)}{f_2(h(t))}\texte^{(h(t)-t)(h(t)+t)/2} \le \Vert f_1\Vert_\infty\Vert 1/f_2\Vert_\infty\,\texte^{-\epsilon t/2}\texte^{-\epsilon h(t)/2}.
\end{equation}
This shows $h'(t) < 1$ and thus that $u\mapsto h(u) - u$ is non-increasing for $u-t$ small positive. As this reinforces the conditions in \eqref{E:4.24} and ~$h$ is non-decreasing by definition, we deduce that
\begin{equation}
0\le h(u)\le u-\epsilon,\quad u\ge t.
\end{equation}
But then
\begin{equation}
h'(u)\texte^{\epsilon h(u)/2}\le \Vert f_1\Vert_\infty\Vert 1/f_2\Vert_\infty\,\texte^{-\epsilon u/2},\quad u\ge t.
\end{equation}
which integrates into
\begin{equation}
\texte^{\epsilon h(u)/2}-\texte^{\epsilon h(t)/2}
\le \Vert f_1\Vert_\infty\Vert 1/f_2\Vert_\infty\,\bigl(\texte^{-\epsilon t/2}-\texte^{-\epsilon u/2}\bigr),\quad u\ge t.
\end{equation}
This contradicts the fact that~$h(u)\to\infty$ as~$u\to\infty$. No~$t>0$ and~$\epsilon>0$ satisfying \eqref{E:4.24} can thus exist, proving \eqref{E:4.22}.

Next suppose that, for some ~$\epsilon>0$ and~$t<-\epsilon$, 
\begin{equation}
\label{E:4.29}
t+\epsilon\le h(t)\le 0\,\,\wedge\,\,\Vert f_1\Vert_\infty\Vert 1/f_2\Vert_\infty\,\texte^{\epsilon t/2}<1.
\end{equation}
Then 
\begin{equation}
\frac{f_1(t)}{f_2(h(t))}\texte^{(h(t)-t)(h(t)+t)/2} \le \Vert f_1\Vert_\infty\Vert 1/f_2\Vert_\infty\,\texte^{\epsilon t/2}\texte^{\epsilon h(t)/2}.
\end{equation}
This shows $h'(t) < 1$ and thus that $u\mapsto h(u) - u$ is non-increasing for $u-t$ small negative. As this again reinforces the conditions in \eqref{E:4.29} and ~$h$ is non-decreasing by definition, we deduce that
\begin{equation}
u + \epsilon \le h(u)\le 0,\quad u\le t.
\end{equation}
But then
\begin{equation}
h'(u)\texte^{-\epsilon h(u)/2}\le \Vert f_1\Vert_\infty\Vert 1/f_2\Vert_\infty\,\texte^{\epsilon u/2},\quad u\le t,
\end{equation}
which integrates into
\begin{equation}
\texte^{-\epsilon h(u)/2}-\texte^{-\epsilon h(t)/2}
\le \Vert f_1\Vert_\infty\Vert 1/f_2\Vert_\infty\,\bigl(\texte^{\epsilon t/2}-\texte^{\epsilon u/2}\bigr),\quad u\le t.
\end{equation}
This contradicts the fact that~$h(u)\to-\infty$ as~$u\to-\infty$. No~$t<0$ and~$\epsilon>0$ satisfying \eqref{E:4.29} can thus exist, proving \eqref{E:4.23}.
\end{proofsect}

With Lemmas~\ref{lemma-4.3} and~\ref{lemma-4.4} in hand, we now give:

\begin{proofsect}{Proof of Proposition~\ref{prop-4.1}}
First observe that $\Vert f_1\Vert_\infty\Vert 1/f_2\Vert_\infty \ge 1$ for otherwise $f_1(t) < f_2(t)$ for all $t\in \mathbb R$, contradicting $\int_{\mathbb R} f_1(t) \mu(\textd t) = \int_{\mathbb R} f_2(t) \mu(\textd t) = 1 $. Similarly, $\Vert f_2\Vert_\infty \Vert 1/f_1\Vert_\infty \ge 1$ holds as well. 
This allows us to define
\begin{equation}
A:=2\sigma^2\log\Bigl(\,\max\bigl\{\Vert f_1\Vert_\infty\Vert 1/f_2\Vert_\infty,\Vert f_2\Vert_\infty\Vert 1/f_1\Vert_\infty\bigr\}\Bigr).
\end{equation}
Note that, since~$h$ is continuous, \twoeqref{E:4.22}{E:4.23} forces $h(t)>0$ on~$[\sqrt A,\infty)$ and~$h(t)<0$ on~$(-\infty,-\sqrt A]$. From \twoeqref{E:4.8}{E:4.9} and \twoeqref{E:4.22}{E:4.23} we then get
\begin{equation}
\label{E:4.35}
\sup_{t\colon |t|\ge \sqrt A}\bigl|\,h(t)-t\bigr|\le \sqrt A.
\end{equation}
For~$t\in[-\sqrt A,\sqrt A]$ the monotonicity of~$h$ shows
\begin{equation}
-2\sqrt A\le h(-\sqrt A)\le h(t)\le h(\sqrt A)\le 2\sqrt A
\end{equation}
and so
\begin{equation}
\sup_{t\colon |t|\le\sqrt A}\bigl|\,h(t)-t\bigr|\le 3\sqrt A.
\end{equation}
Along with \eqref{E:4.35}, this gives the desired claim.
\end{proofsect}

\subsection{Discrete-valued steps and proof of Theorem~\ref{thm-4.1}}
We will now move to the proof of Theorem~\ref{thm-4.1}. Here we note that Proposition~\ref{prop-4.1} gives a bound on the coupling distance but does not provide any bound on the probability of equality. This is fixed in:

\begin{lemma}
\label{lemma-4.5}
Let~$\mu=\NN(0,\sigma^2)$ and $f\colon\R\to(0,\infty)$ be such that~$\int f\textd\mu=1$. Let~$X$ have the law~$\mu$ and let~$Y$ have the law~$f\mu$. Assuming~$f$ and~$1/f$ to be bounded and~$\mu(f\ne1)>0$, for all strictly positive $\eta\ge \log \Vert 1/f\Vert_\infty$ there exists a coupling of~$X$,~$Y$ and a zero-one valued random variable~$B$ such that
\begin{equation}
\label{E:4.38}
\Vert X-Y\Vert_\infty\le 3\sigma
\sqrt{2\log\Biggl(\max\biggl\{\frac{\texte^{2\eta}\Vert f\Vert_\infty-1}{\texte^{2\eta}-1},\,\texte^\eta+1\biggr\}\Biggr)},
\end{equation}
\begin{equation}
\label{E:4.43}
X=Y\text{\rm\ on }\{B=1\},
\end{equation}
\begin{equation}
\label{E:4.44}
\text{\rm $B$ and~$X$ are independent}
\end{equation}
and
\begin{equation}
\label{E:4.45}
P(B=1) = \texte^{-2\eta}.
\end{equation}
(When~$\mu(f\ne 1)=0$, we have~$X\laweq Y$ and we can use~$Y:=X$ and~$B:=1$.)
\end{lemma}

\begin{proofsect}{Proof}
Let~$\eta \ge \log \Vert 1/f\Vert_\infty$ be strictly positive. Sample~$X$ from~$\mu$ and let~$B$ be an independent Bernoulli with $P(B=1):=\texte^{-2\eta}$. Then set~$Y:=X$ when~$B=1$ and~$Y:=h(X)$ when $B=0$, for~$h$ as in \eqref{E:4.4} with~$F$ being the CDF of~$\mu$ and~$G$ the CDF of~$\tilde f\mu$ where
\begin{equation}
\tilde f:=\frac{f-\texte^{-2\eta}}{1-\texte^{-2\eta}}.
\end{equation}
As is readily checked, this defines a coupling of~$(X,Y,B)$ with the given marginals and such that \twoeqref{E:4.43}{E:4.45} hold.

In order to prove \eqref{E:4.38}, note that the bound in Proposition~\ref{prop-4.1} gives
\begin{equation}
\label{E:4.47}
\Vert X-Y\Vert_\infty\le 3\sigma\sqrt{2\log\bigl(\max\{\,\Vert \tilde f\Vert_\infty,\,\Vert1/\tilde f\Vert_\infty\}\bigr)}.
\end{equation}
A computation shows
\begin{equation}
\Vert\tilde f\Vert_\infty\le\frac{\Vert f\Vert_\infty-\texte^{-2\eta}}{1-\texte^{-2\eta}}
\end{equation}
while
\begin{equation}
\Vert1/\tilde f\Vert_\infty\le\frac{1-\texte^{-2\eta}}{\texte^{-\eta}-\texte^{-2\eta}}=\texte^\eta+1.
\end{equation}
Inserting these in \eqref{E:4.47}, we get \eqref{E:4.38}.
\end{proofsect}

Lemma~\ref{lemma-4.5} is sufficient to handle the coupling of the increments $\xi_k^{\DG}(x)$ and~$\xi_k^{\GFF}(x)$ for all~$k=1,\dots,n$. For the discrete-valued step (the case $k=0$ of  the DG model) we will instead need a different argument for uniform closeness. This comes in:

\begin{lemma}
\label{lemma-4.6}
For each $\beta>0$ there exists~$C_0>0$ such that the following holds for all~$z\in\R$: If~$X\laweq\NN(0,1/\beta)$ and~$Y$ takes values in~$z+\Z$ with probabilities
\begin{equation}
P\bigl(Y=z+n) = \texte^{v_0(z)-\frac\beta2(z+n)^2},\quad n\in\Z,
\end{equation}
for~$v_0$ as in \eqref{E:2.21}, then there exists a coupling of~$X$ and~$Y$ such that
    \begin{equation}
    \label{E:4.51}
         \Vert X-Y\Vert_{\infty} \le C_0.
    \end{equation}
\end{lemma}

\begin{proofsect}{Proof}
Since the law of~$Y$ does not change if~$z$ is changed by an integer (this also uses the periodicity of $v_0$), we may and will henceforth assume $z \in [0,1)$. The proof uses similar arguments as those underlying Proposition~\ref{prop-4.1} albeit adapted to the discrete setting.

We start by noting that, in light of the strict positivity of the Gaussian density, there exists a unique non-decreasing  sequence  $\{ u_k \}_{k\in\Z}$ of reals such that
\begin{equation}
\label{E:4.52}
      \int^\infty_{z+ u_k}\frac1{\sqrt{2\pi/\beta}}\,\texte^{-\frac\beta2 x^2}\textd x = P(Y \ge z + k),\quad k\in\Z.
\end{equation}
The positivity also implies $ u_k  \to \pm \infty$ as $k\to \pm \infty$. Now let $h\colon\R\to \R$ be such that
\begin{equation}
\label{E:5.47a}
h(x):=z+k\,\,\text{ if }\,\,x\in[z+ u_k ,z+ u_{k+1})
\end{equation}
for $k\in\Z$.
 The definition \eqref{E:4.52} then gives $h(X) \laweq Y$.  

In order to show \eqref{E:4.51},  it suffices  to prove that~$\sup_{k\in\Z}| u_k -k|<\infty$.  We will show this by proving that $|u_k-k|<1$ for~$k$ large positive; the argument for negative~$k$ is analogous and we omit it. 
We start by invoking  the explicit form of the law of~$Y$ along with a substitution and some algebra to turn \eqref{E:4.52} into 
    \begin{equation}
        \int_{0}^\infty \texte^{-\frac{\beta}{2}(x+ u_k  +z)^2} \textd x = \texte^{v_0(z)} \sqrt{2\pi/\beta} \,\, \texte^{-\frac{\beta}{2}(z + k)^2}\biggl(1+  \sum_{j> k} \texte^{-\frac{\beta}{2} (z +j )^2+\frac{\beta}{2}(z + k)^2}\biggr).
    \end{equation}
    Abbreviating $ c(z) := \beta\,\texte^{v_0(z)} \sqrt{2\pi/\beta}$ and $A_\beta(k) :=  \sum_{j> k} \texte^{-\frac{\beta}{2} (z+j )^2+\frac{\beta}{2}(z+ k)^2}$, this becomes 
    \begin{equation}
        \int_{0}^\infty \texte^{-\frac{\beta}{2}x^2 - \beta (z+ u_k ) x} \textd x =  \beta^{-1} c(z)  \,\texte^{-\frac{\beta}{2}(z+ k)^2+ \frac{\beta}{2}(z+ u_k )^2}  \left(1+ A_\beta(k)\right), 
    \end{equation}
 which, assuming that $|u_k|\ge1$ and thus $u_k+z\ne0$, gives 
    \begin{equation}
        \label{E:4.57}
        \int_{0}^\infty \texte^{-\frac{x^2}{2\beta (z+ u_k )^2 } - x} \textd x =  c(z) (z+ u_k ) \,\texte^{\frac{\beta}{2} (2z+ k+u_k)(u_k-k)}\left(1+ A_\beta(k)\right)
    \end{equation}
by another substitution under the integral  and a rewrite of the exponent on the right. 

We will now use \eqref{E:4.57} to control $ u_k -k$ for~$k$ large positive. First note that $ u_k  \to \infty$ as $k\to \infty$ tells us there exists $k_1\in\Z$ such that $ u_k  \ge 1$ for all $k\ge k_1$. Suppose we had $ u_k  \le k-1$ for some $k\ge k_1$. Then \eqref{E:4.57} gives 
    \begin{equation}
        \label{E:4.58}
        \int_{0}^\infty \texte^{-\frac{x^2}{2\beta (1+z)^2}  - x} \textd x \le  c(z) (z+k-1)  \, \texte^{-\frac{\beta}{2} (2(k+z)+ 1)}\left(1+ A_\beta(k)\right).
    \end{equation}
 But the right-hand side is   at most $c_1 k\texte^{-\beta k}[1+A_\beta(k)]$  for $c_1:=\sup_{0\le z<1}c(z)$ and  a direct estimate shows
    \begin{equation}
        A_\beta(k) = \sum_{j>k} \texte^{-\frac{\beta}{2} (j-k)(k+j+2z)} = \sum_{l=1}^\infty \texte^{-\frac{\beta}{2} l(2k+l+2z)} \le  \sum_{l=1}^\infty \texte^{-\frac{\beta}{2} (l^2+2kl)}  < \infty. 
    \end{equation}
 This implies the right-hand side of \eqref{E:4.58} vanishes in the limit as $k\to \infty$  uniformly in $z\in[0,1)$,  which is absurd since the left-hand side remains uniformly positive. It follows that there exists $K_1\ge k_1$ so that $k\ge K_1$ implies $ u_k > k-1$  for all~$z\in[0,1)$. 
    
Next  suppose that  $ u_k  \ge k+1$. Denoting $c_2:=\inf_{0\le z<1}c(z)$, the equality in \eqref{E:4.57} gives 
    \begin{equation}
        1 = \int_0^\infty \texte^{-x} \textd x \ge  c(z)  (z+k+ 1) 
        \,\texte^{\frac{\beta}{2} (2(k+z) + 1)}\ge c_2 k \texte^{\beta k}
    \end{equation}
which is again absurd for~$k$ large because  the right-hand side  diverges as~$k\to\infty$.  Hence  there exists $K_1'$ such that $k\ge K_1'$ implies $ u_k  < k+1$ for all~$z\in[0,1)$.  In particular, $ |u_k-k|<1 $ holds for~$k\ge K:=\max\{K_1,K_1'\}$  uniformly in~$z\in[0,1)$. 
\end{proofsect}

We are now ready to give:

\begin{proofsect}{Proof of Theorem~\ref{thm-4.1}}
Fix~$\beta\in(0,\beta_\cc)$ and $n\ge1$. For~$k=1,\dots,n $, define~$R_k$ by
\begin{equation}
    \label{E:4.48}
    R_k := 2\sup_{z,z'\in[0,1]} \bigl|\,v_{k}(z) - b v_{k-1}(z')\bigr|.
    \end{equation}
Then sample i.i.d.\ random variables \eqref{E:4.2b} with common law~$\NN(0,1/\beta)$ and independent Bernoulli random variables \eqref{E:4.3} with
\begin{equation}
P\bigl(B_k(x)=1\bigr)=\texte^{-R_k},\quad x\in\Lambda_{n-k},\,k=1,\dots,n.
\end{equation}
 Note  that the density governing the definition of~$\xi_k^{\DG}(x)$ from $\xi_k^{\GFF}(x)$ takes the form
\begin{equation}
     f_{k,x}(z):=\texte^{v_{k}(\varphi_k(x)) - b v_{k-1}(z+\varphi_k(x))},
    \end{equation}
where $\varphi_k(x):=\sum_{j=k+1}^{ n }\xi_j^{\DG}(x)$  and where, in particular,  $\varphi_{n}  (x)  = 0$. Lemma~\ref{lemma-4.5} with~$\eta:=R_k/2$ (which exceeds $\log\Vert 1/f_{k,x}\Vert_\infty$)  allows us to recursively construct $\xi^{\DG}_k(x)$ for $x\in\Lambda_n$ and $k= n, \dots,1$ as a deterministic function of the GFF increments \eqref{E:4.2b} and the Bernoulli's \eqref{E:4.3}. For~$k=0$ we in turn use  \eqref{E:5.47a} in the proof of  Lemma~\ref{lemma-4.6}.

The bound \eqref{E:4.38} now shows
\begin{equation}
\label{E:4.67}
        \bigl\Vert\, \xi_k^{\DG}(x) - \xi_k^{\GFF}(x)\bigr\Vert_{\infty} \le 3\beta^{-1/2} \sqrt{2\log (\texte^{R_k} + 1)}
    \end{equation}
for all $k=1,\dots,n$ and~$x\in\Lambda_n$, 
while \eqref{E:4.51}  dominates  the $k=0$ case by a constant~$C_0$.
Theorem~\ref{thm-3.4} asserts that~$R_k\to0$ as~$k\to\infty$ and so the right-hand side of \eqref{E:4.67} is bounded uniformly on~$k$.  
Denoting
    \begin{equation}
        C_1 := 3\beta^{-1/2}  \sqrt{2\log (\texte^{\sup_{k\ge1} R_k} + 1)}, 
    \end{equation}
the  bound \eqref{E:4.1} holds with~$C:=\max\{C_0,C_1\}$.  That $\limsup_{k\to \infty} k^{-1}\log R_k < 0$ follows from the exponential decay of~$R_k$ proved in Theorem~\ref{thm-3.4}. 
\end{proofsect}


\section{Tightness of DG extrema}
\label{sec5}\noindent
With the needed tools assembled, we now finally turn our attention to the extremal behavior of the DG model. The goal of this section is to establish tightness of the maximum (centered by~$m_n$) as well as  of  the level sets (above $m_n-\lambda$ with~$\lambda>0$ fixed). These will be useful in the proofs in Section~\ref{sec6}.

\subsection{Statements and preliminaries}
Recall that~$P_{n,\beta}$ denotes the law of the DG model at inverse temperature~$\beta$ and depth~$n$. We start with the statement of the relevant results to be proved in this section: 

\begin{theorem}
    \label{thm-5.1}
    For all $\beta\in(0,\beta_\cc)$ and with~$m_n$ as in \eqref{E:1.6a},
\begin{equation}
\lim_{u\to \infty}  \sup_{n\ge 1} P_{n,\beta}\Bigl(\bigl|\,\max_{x\in \Lambda_n} \varphi_x^{\DG} -m_n\bigr|\ge u\Bigr) = 0.
\end{equation}
\end{theorem}

\begin{theorem}
    \label{thm-5.2}
Given~$\lambda\in\R$ and a sample~$\varphi^{\DG}$ from~$P_{n,\beta}$, let
\begin{equation}
\label{E:6.2w}
G_n(\lambda)  := \bigl|\{x\in \Lambda_n\colon \varphi_x^{\DG} \ge m_n - \lambda\}\bigr|.
\end{equation}
For all~$\beta\in(0,\beta_\cc)$ and all~$\lambda>0$,
\begin{equation}
\lim_{u\to \infty}  \sup_{n\ge 1} P_{n,\beta}\bigl(G_n(\lambda)\ge u\bigr) = 0.
\end{equation}
\end{theorem}

The proofs will rely on some standard results about the GFF which we will state next. Recall that the GFF at inverse temperature~$\beta$ on the hierarchical lattice is nothing but a Branching Random Walk with steps distributed  as  $\NN(0,1/\beta)$. The tightness  and limit laws of the maxima  of general Branching Random Walks  have  been studied extensively (e.g., by A\"idekon~\cite{Aidekon} and Madaule~\cite{Madaule2}). For Gaussian step distributions, the calculations become very explicit. We summarize what we need in:

\begin{lemma}
    \label{lemma-5.1}
    There exists $a  > 0$ such that for all $t > 0$ and all~$\beta>0$,
    \begin{equation}
        \sup_{n\ge 1} P_{n,\beta}'\Bigl( \bigl|\,\max_{x\in \Lambda_n} \varphi_x^{\GFF} - m_n\bigr| > t\Bigr) \le  \frac{1}{a} \,\texte^{-a\sqrt{\beta}\, t}.
    \end{equation}
\end{lemma}

\begin{proofsect}{Proof}
This is proved, modulo scaling by~$\beta$ and quotations of ``ballot theorems,'' in the review article~\cite[Lemmas~7.3 and~7.15]{B-notes}. Alternatively, consult Mallein~\cite[Theorem~4.1]{Mallein} which gives upper tail for general Branching Random Walks as well as the key technical input $\inf_{n\ge1} P_{n,\beta}'(\max_{x\in \Lambda_n}\varphi_x^{\GFF}\ge m_n)>0$ for the lower tail. 
\end{proofsect}

\begin{lemma}
\label{lemma-5.3}
    There exist $c',C'>0$ and, for all $\beta>0$, all $t > 0$ and all $u < t$, there exists $n_0\ge1$ such that for all $n\ge n_0$ and all $z\in \Lambda_n$, 
    \begin{equation}
        \label{E:5.3}
        P_{n,\beta}'\left(\varphi_z^{\GFF} \ge m_n + u, \max_{x\in \Lambda_n} \varphi_x^{\GFF} \le  m_n + t\right) \le C' b^{-n} \bigl(1+ (t \lor (t-u))^2\beta\bigr) \texte^{-\sqrt{\beta} c' u} .
    \end{equation}
\end{lemma}

\begin{proofsect}{Proof}
This follows by conditioning on~$\varphi_z^{\GFF}$ and applying a ``ballot estimate.'' See, e.g., the proof of \cite[Lemma~6.9]{BL3} for the (harder) case of lattice GFF.
\end{proofsect}

We note that $n_0$ in Lem\-ma~\ref{lemma-5.3} can and will be chosen so that $m_{n_0}+u$ is a large positive number. We also record  an  estimate on the tail of sum of Bernoulli random variables based on the Chernoff bound:

\begin{lemma}
    \label{lemma-5.4}
     Let $\{Z_k\}_{k\ge 0}$ be independent zero-one valued random variables. Abbreviate $p_k: = P(Z_k = 0)$. Then for all natural $k \le \ell$ and $r > p_{k,\ell}:=\sum_{i=k}^{\ell} p_i$,  
\begin{equation}
    \label{E:5.6}
    P\biggl(\,\sum_{i=k}^{ \ell} 1_{\{Z_i = 0\}} > r\biggr) \le \exp\Bigl\{r - p_{k,\ell} + r \log p_{k,\ell}- r\log r\Bigr\}.
\end{equation}
In particular, if $\sum_{i\ge0} p_i < \infty$ then there exists $r_0 \ge 1$ such that for all $r \ge r_0$
\begin{equation}
    \label{E:5.14}
    P\biggl(\,\sum_{i=k}^{\ell} 1_{\{Z_i = 0\}} > r\biggr) \le \exp\biggl\{-\frac{1}{2} r\log r\biggr\}
\end{equation}
holds  for all $\ell\ge k\ge 0$.
\end{lemma}

\begin{proofsect}{Proof}
Given any~$\lambda>0$, the Chernoff bound shows 
\begin{equation}
    \begin{aligned}
        P\biggl(\,\sum_{i=k}^{\ell} &\,1_{\{Z_i = 0\}} > r\biggr) \le \texte^{-r\lambda} E\biggl[\exp\Bigl(\lambda \sum_{i=k}^{\ell} 1_{\{Z_i = 0\}} \Bigr)\biggr] \\
        &= \exp\biggl\{ - r\lambda  + \sum_{i=k}^{\ell} \log \bigl(1+(\texte^{\lambda} - 1)p_i\bigr)\biggr\} 
        \le \exp\biggl\{ - r\lambda  + \sum_{i=k}^{\ell} (\texte^{\lambda} - 1)p_i\biggr\}.
    \end{aligned}
 \end{equation}
    Now set $\lambda := \log\left(r/p_{k,\ell}\right)$ and observe that~$r>p_{k,\ell}$ gives $\lambda> 0$. An explicit evaluation then shows \eqref{E:5.6}.  
 For \eqref{E:5.14} we observe that
    \begin{equation}
 r - p_{k,\ell} + r \log p_{k,\ell}\le r \Biggl[\,1+\log\biggl(\,\sum_{i=1}^\infty p_i\biggr)\Biggr]\le\frac12 r\log r
    \end{equation}
once~$r\ge [\texte \sum_{i\ge0}  p_i]^2$. 
\end{proofsect}


\subsection{Key proposition}
We now move to the statement and proof of a technical proposition that provides the key step in the proof of both Theorem~\ref{thm-5.1} and~\ref{thm-5.2}. 

\begin{proposition}
\label{prop-5.1}
    For each~$\beta\in(0,\beta_\cc)$ there exists $C''>0$ and, for all $t > 0$ and all $u < t$, there exists $n_0\ge1$ such that for all $n\ge n_0$ and all $z\in \Lambda_n$,     \begin{equation}
        \label{E:5.21}
        P\left(\varphi_z^{\DG} \ge m_n + u, \max_{x\in \Lambda_n} \varphi_x^{\GFF} \le  m_n + t\right) \le C'' b^{-n} \Bigl(1+\bigl(t \lor (t-u)\bigr)^2\beta\Bigr) \texte^{-\sqrt\beta\, c' u} ,
    \end{equation}
where~$c'$ is as in Lemma~\ref{lemma-5.3} and~$P$ is the coupling law from Theorem~\ref{thm-4.1}.
\end{proposition}

The proof is based on swapping~$\varphi_z^{\DG}$ for~$\varphi_z^{\GFF}$ which effectively reduces the claim to Lemma~\ref{lemma-5.3}. To control the error incurred by the swap, we have to bound the probability that the difference $\varphi_z^{\DG}-\varphi_z^{\GFF}$ is large on the background of the event that $\varphi_z^{\DG}$ is itself large (as this provides the much needed $b^{-n}$ term). 

The argument relies on the coupling from Theorem~\ref{thm-4.1} but its use is complicated by the restriction on the GFF-maximum (which is needed to beat a polynomial term in~\eqref{E:5.15}). Indeed, this maximum is a global event that involves, in principle, all the coupling variables. That being said, this restriction is not needed if the difference $\varphi_x^{\DG} - \varphi_x^{\GFF}$ is sufficiently large. This will be handled using:

 \begin{lemma}
    \label{lemma-5.5}
    Let $\beta\in(0,\beta_\cc)$ and let~$C$ and $\{R_k\}_{k\ge1}$ be as in Theorem~\ref{thm-4.1}. There exists~$C'''>0$ such that for all~$n\ge1$, all $x\in \Lambda_n$, all $u\in\R$, and all $k \ge 1$,
     \begin{equation}
     \label{E:5.15}
         P\bigl(\varphi_x^{\DG} \ge m_n + u, \varphi_x^{\DG} - \varphi_x^{\GFF} > Ck\bigr) \le C''' nb^{-n} \texte^{-\sqrt\beta\, c' u} \sum_{j=k}^\infty R_j,
     \end{equation}
where~$c'$ is as in Lemma~\ref{lemma-5.3}.
 \end{lemma}
 
 \begin{proofsect}{Proof}
Fix $x\in \Lambda_n$ and note that, since the event involves only the coupling variables on the path from~$x$ to the root, we may write $\xi_j^{\DG}$, $\xi_j^{\GFF}$, and $B_i$ instead of $\xi_j^{\DG}(x)$, $\xi_j^{\GFF}(x)$, and $B_i(x)$, respectively. Define
\begin{equation}
\label{E:5.10}
 \tau =\tau(x):= \sup \bigl\{j=1,\dots,n : B_j = 0\bigr\} 
\end{equation}
 to be the height that (as far as the Bernoulli's \eqref{E:4.3} are concerned) the coupling fails for the first time along the path from the root to~$x$. By Theorem~\ref{thm-4.1} and our choice of~$C$, the event $\varphi_x^{\DG} - \varphi_x^{\GFF} > Ck$ implies $\tau \ge k$. The probability in \eqref{E:5.15} is thus bounded by 
     \begin{equation}
        \label{E:5.16}
         P(\varphi_x^{\DG} \ge m_n + u, \tau \ge k) = \sum_{i=k}^{n} P(\varphi_x^{\DG} \ge m_n + u, \tau = i).
     \end{equation}
Next we note that,  since only the variables on the ``path'' to a single~$x\in\Lambda_n$ are involved, the DG-process and the GFF-process are mutually absolutely continuous. More precisely, let $\FF_k:=\sigma(\xi_i^{\DG}, \xi_i^{\GFF} , B_i\colon i\ge k)$ and abbreviate $\xi_{\ge k} := \sum_{i=k}^{n}  \xi_i$, with the superscript~``DG'' or ``GFF'' added depending on the context. Recall the notation $\frakq_k(\cdot|\varphi)$ from \eqref{E:1.23} and let~$\mu$ denote the law of~$\NN(0,1/\beta)$. Then for each Borel~$A\subseteq\R^k$, \eqref{E:2.25} along with \eqref{E:2.42} and~\eqref{E:2.44}~give 
\begin{equation}
\begin{aligned}
P\biggl(\,(\xi_0^{\DG},\dots,&\xi_{k-1}^{\DG}) \in A \,\bigg| \,\FF_k \biggr) 
         = \int_A \bigotimes_{i=0}^{k-1}\frakq_i(\textd\xi_i|\xi_{\ge i+1})
        \\
        &\le \int_A \Bigl(\,\prod_{j=0}^{k-1} \texte^{R_j/2}\Bigr) \bigotimes_{i=0}^{k-1} \mu (\textd\xi_i)
        \le  c P\biggl((\xi_0^{\GFF},\dots,\xi_{k-1}^{\GFF}) \in A  \,\bigg|\,\FF_k\biggr),
\end{aligned}
\end{equation}
where $R_0 := 2\texte^{\Vert v_0\Vert}$ and $c := \exp\{\frac{1}{2}\sum_{j=0}^{\infty} R_j\} < \infty$.  
On the event $\{|\xi_{\ge k}^{\DG} - \xi_{\ge k}^{\GFF}| \le C\}$, which lies in~$\FF_k$, we thus get 
\begin{equation}
\label{E:5.13}
\begin{aligned}
P\bigl(\varphi_x^{\DG} \ge m_n + u\,\big|\,\FF_k\bigr) &= P\bigl(\varphi_x^{\DG} - \xi_{\ge k}^{\DG} \ge m_n + u - \xi_{\ge k}^{\DG} \,\big|\,\FF_k\bigr)  \\
        &= P\biggl(\,\sum_{j=0}^{k-1} \xi_j^{\DG} \ge m_n + u - \xi_{\ge k}^{\DG} \,\bigg|\,\FF_k\biggr)  \\
        &\le cP\biggl(\,\sum_{j=0}^{k-1} \xi_j^{\GFF} \ge m_n + u - \xi_{\ge k}^{\DG} \,\bigg|\,\FF_k\biggr)  \\
        &\le cP\biggl(\,\sum_{j=0}^{k-1} \xi_j^{\GFF} \ge m_n + u - \xi_{\ge k}^{\GFF} - C \,\bigg|\,\FF_k\biggr)  \\
        & = c P\bigl(\varphi_x^{\GFF} \ge m_n + u - C \,\big|\,\FF_k\bigr),
\end{aligned}
\end{equation}
where we use that~$\xi_{\ge k}^{\DG}$ and $\xi_{\ge k}^{\GFF}$ are~$\FF_k$-measurable to treat them as constants under the conditional probability.

Since $\xi_{\ge i+1}^{\DG} = \xi_{\ge i+1}^{\GFF}$ on~$\{\tau=i\}$, Theorem~\ref{thm-4.1} shows
\begin{equation}
 \{\tau=i\}\subseteq \bigl\{|\xi_{\ge i}^{\DG} - \xi_{\ge i}^{\GFF}| \le C\bigr\}. 
\end{equation}
 This enables the inequality in \eqref{E:5.13} which then yields 
\begin{equation}
\begin{aligned}
P\bigl(\varphi_x^{\DG} \ge m_n + u, \tau = i\bigr) &= E\Bigl[ P\bigl(\varphi_x^{\DG} \ge m_n + u\,\bigl|\, \FF_i \bigr)  1_{\{\tau = i\}} 1_{\{|\xi_{\ge i}^{\DG} - \xi_{\ge i}^{\GFF}| \le C\}} \Bigr]  \\
        &\le c E\Bigl[ P\bigl(\varphi_x^{\GFF} \ge m_n + u-C\,\bigl|\, \FF_i \bigr)  1_{\{\tau = i\}}  \Bigr]   \\
        &= cP\bigl(\varphi_x^{\GFF} \ge m_n + u - C\bigr) P(\tau = i),
\end{aligned}
\end{equation}
where in the final step we used that the Bernoulli's \eqref{E:4.3} are independent of the GFF increments, as stated in Theorem~\ref{thm-4.1}(2).

From \eqref{E:5.16} and 
and explicit calculation we now deduce 
    \begin{equation}
        P\bigl(\varphi_x^{\DG} \ge m_n + u, \varphi_x^{\DG} - \varphi_x^{\GFF} > Ck\bigr) \le C''' nb^{-n} \texte^{-\sqrt\beta\, c' u} P(\tau \ge k),
    \end{equation}
    where $C'''$ absorbs the constant~$c$ along with various constants that arise in the calculation.
 Noting that
\begin{equation}
P(\tau \ge k) = 1- \prod_{j=k}^{n} \texte^{-R_k} \le \sum_{j=k}^{n} R_j \le \sum_{j=k}^\infty R_j,
\end{equation}
the claim follows. 
 \end{proofsect}

We are now ready to give:

\begin{proofsect}{Proof of Proposition~\ref{prop-5.1}}
Given integers~$\ell, r\ge1$ to be chosen later, we write
\begin{equation}
\label{E:6.19u}
\begin{aligned}
\qquad
\bigl\{\varphi_z^{\DG} \ge m_n + &u\bigr\}
\subseteq\bigl\{\varphi_z^{\GFF} \ge m_n + u-\ell\bigr\}
\\
&\cup
\bigl\{\varphi_z^{\DG} \ge m_n + u,\,\ell < \varphi_z^{\DG} - \varphi_z^{\GFF} \le C r \log n\bigr\}
\\
&\qquad\qquad\qquad\cup
\bigl\{\varphi_z^{\DG} \ge m_n + u,\,  \varphi_z^{\DG} - \varphi_z^{\GFF} > C r \log n\bigr\}.
\qquad
\end{aligned}
\end{equation}
For first event on the right, Lemma~\ref{lemma-5.3} gives
\begin{equation}
\begin{aligned}
\label{E:5.28}
     P\Bigl(\,\max_{x\in \Lambda_n} \varphi_x^{\GFF} \le m_n + t, &\, \varphi_z^{\GFF} \ge m_n + u - \ell\Bigr) 
     \\
     &\le C' b^{-n}  \bigl(1+(t \lor (t-u+\ell) )^2 \beta \bigr)\texte^{-\sqrt\beta c' (u-\ell)} \\
     &\le C'\ell^2\texte^{\sqrt\beta c'\ell} b^{-n}  \bigl(1+(t \lor (t-u+1) )^2 \beta \bigr)\texte^{-\sqrt\beta \,c' u} .
\end{aligned}
\end{equation}
For the third event, we can drop the restriction on the maximum and bound the result via Lemma~\ref{lemma-5.5} to get
\begin{equation}
    \label{E:5.22}
    P\bigl(\varphi_z^{\DG} \ge m_n + u, \varphi_z^{\DG} - \varphi_z^{\GFF} > C r \log n \bigr)\le C''' nb^{-n} \texte^{-\sqrt\beta c' u} \sum_{j=\lfloor r \log n\rfloor}^\infty R_j.
\end{equation} 
Since $\sum_{j=k}^\infty R_j$ is exponentially decaying in~$k$, we can choose $ r$ so large that  
\begin{equation}
    \label{E:5.23}
    \sum^\infty_{j=\lfloor r \log n\rfloor} R_j \le \frac{1}{n^2}
\end{equation}
once~$n$ is sufficiently large. Under such circumstances, this term contributes at most a constant times $n^{-1}b^{-n}\texte^{-\sqrt\beta c' u}$.

It thus remains to estimate the contribution of the second event on the right of \eqref{E:6.19u} to the probability in question. Here we first partition the event as
\begin{equation}
\begin{aligned}
\label{E:5.24}
    P\Bigl(&\max_{x\in \Lambda_n} \varphi_x^{\GFF} \le m_n + t,  \varphi_z^{\DG} \ge m_n + u, \ell < \varphi_z^{\DG} - \varphi_z^{\GFF} \le C r \log n  \Bigr)
    \\
    &\le \sum_{j=0}^{\lfloor C r \log n \rfloor} P\Bigl(\,\max_{x\in \Lambda_n} \varphi_x^{\GFF} \le m_n + t,  \varphi_z^{\DG} \ge m_n + u, \varphi_z^{\DG} - \varphi_z^{\GFF} - \ell-j\in(0,1] \Bigr)
\end{aligned}
\end{equation}
Invoking the Bernoulli's in the coupling, the probability on the right is then bounded  further by the quantities 
\begin{equation}
\label{E:6.25}
\begin{aligned}
P&\Bigl(\max_{x\in \Lambda_n} \varphi_x^{\GFF} \le m_n + t,  \varphi_z^{\GFF} \ge m_n + u - (\ell+j+1), \varphi_z^{\DG} - \varphi_z^{\GFF} > \ell+j\Bigr)\\
    &\le P\biggl(\,\max_{x\in \Lambda_n} \varphi_x^{\GFF} \le m_n + t,  \varphi_z^{\GFF} \ge m_n + u - (\ell+j+1),\, C\sum_{k=0}^{n} 1_{\{B_k(z) = 0\}} > \ell+j\biggr)
\end{aligned}
\end{equation}
where we set $B_0(z):= 0$ and invoked \eqref{E:4.1}. 

 Using the  independence stated in Theorem~\ref{thm-4.1}(2), the last probability  in \eqref{E:6.25} factors into  the product
\begin{equation}
\label{E:5.25}
P\Bigl(\,\max_{x\in \Lambda_n} \varphi_x^{\GFF} \le m_n + t,  \varphi_z^{\GFF} \ge m_n + u - (\ell+j+1)\Bigr) P\biggl(\, \sum_{k=0}^{n}  1_{\{B_k(z) = 0\}} > \frac{\ell+j}{C}\biggr). 
\end{equation}
For $\ell$ large, \eqref{E:5.14} gives
\begin{equation}
    \label{E:5.26}
    P\biggl(\, \sum_{k=0}^{n}   1_{\{B_k(z) = 0\}} > \frac{\ell+j}{C}\biggr) \le \exp\left(-\frac{1}{2} C^{-1}(\ell + j) \log (C^{-1}(\ell+j))\right).
\end{equation}
Lemma~\ref{lemma-5.3} in turn bounds the probability on the left of \eqref{E:5.25}  by
\begin{equation}
    \label{E:5.27}
    C' b^{-n} \Bigl(1+\bigl(t \lor (t-u+\ell+j+1)\bigr)^2 \beta\Bigr) \texte^{\sqrt{\beta}\, c' (\ell+j+1)}\texte^{-\sqrt{\beta}\, c' u},
\end{equation}
provided that~$n$ is large enough so that $m_n + u - (\ell+C r \log n +1)$ is a large positive number  (enables Lemma~\ref{lemma-5.3}  for all $j$ between 0 and $\lfloor C r \log n\rfloor$).

Combining these observations, the probability on the left of \eqref{E:5.24} is thus bounded by $C' b^{-n}\texte^{-\sqrt\beta\, c' u}$ times 
\begin{equation}
\begin{aligned}
\sum_{j=0}^\infty \Bigl(1+&\bigl(t \lor (t-u+\ell+j+1)\bigr)^2 \beta\Bigr)  \,\texte^{\sqrt\beta\, c' (\ell+j+1)-\frac{1}{2} C^{-1}(\ell + j) \log (C^{-1}(\ell+j))} \\
&\le 
\Bigl(1+\bigl(t \lor (t-u+1)\bigr)^2\beta\Bigr) 
\sum_{j=\ell+1}^\infty \,j^2 \texte^{\sqrt\beta c' j-\frac{1}{2} C^{-1}(j-1) \log (C^{-1}(j-1))}.
\end{aligned}
\end{equation}
The sum on the right converges  and so the expression is at most $C''(1+(t\vee(t-u))^2\beta)$, for a suitable choice of~$C''$.  Together with \twoeqref{E:5.28}{E:5.23}, this implies the assertion. 
\end{proofsect}

\subsection{Proofs of tightness}
With the technical inputs settled, we are ready to prove the main results of this section. We start with the tightness of the centered DG-maximum.

\begin{proofsect}{Proof of Theorem~\ref{thm-5.1}}
We first address the tightness of the upper tail. (This is the difficult part but we have already done most of the work.) 
    Indeed, for any $t > 0$ we have
\begin{equation}
\begin{aligned}
P\Bigl(\,&\max_{x\in \Lambda_n} \varphi_x^{\DG} \ge m_n + u\Bigr) 
\\
&\le  P\Bigl(\,\max_{x\in \Lambda_n} \varphi_x^{\GFF} \le m_n + t, \max_{x\in \Lambda_n} \varphi_x^{\DG} \ge m_n + u\Bigr) + P\Bigl(\,\max_{x\in \Lambda_n} \varphi_x^{\GFF} > m_n + t\Bigr).
\end{aligned}
\end{equation}
By the union bound and Lemma~\ref{lemma-5.1},
\begin{equation}
\begin{aligned}
\label{E:5.30}
    P\Bigl(\,\max_{x\in \Lambda_n} \varphi_x^{\DG} \ge &\,m_n + u\Bigr) 
    \\
    &\le \sum_{z\in\Lambda_n} P\Bigl(\,\max_{x\in \Lambda_n} \varphi_x^{\GFF} \le m_n + t,  \varphi_z^{\DG} \ge m_n + u\Bigr) + \frac{1}{a}\texte^{-a \sqrt \beta\, t}.
\end{aligned}
\end{equation}
Proposition~\ref{prop-5.1} then shows that, for sufficiently large~$n$,
\begin{equation}
    \label{E:5.31}
    P\left(\max_{x\in \Lambda_n} \varphi_x^{\DG} \ge m_n + u\right) \le C'' \Bigl(1+\bigl(t \lor (t-u) \bigr)^2 \beta\Bigr) \texte^{-\sqrt\beta c' u} + \frac{1}{a}\texte^{-a \sqrt \beta\, t}.
\end{equation}
Setting~$t:=u/2$ then shows that the supremum over~$n\ge1$ of the left-hand side tends to zero as~$u\to\infty$.

It remains to prove the tightness of the lower tail. Let $k\in\{1,\dots,n \}$ and recall that $\xi_{\ge k}^{\DG}(x)$ and $\xi_{\ge k}^{\GFF}(x)$ denote the sums of $\xi_{i}^{\DG}(x)$ and $\xi_{i}^{\GFF}(x)$ for~$i=k,\dots,n $. (These are still parametrized by~$x\in\Lambda_n$.) Then for any~$x\in\Lambda_n$, 
\begin{equation}
\label{E:5.29}
    P\Bigl(\xi_{\ge k}^{\DG}(x)=\xi_{\ge k}^{\GFF}(x)\Bigr) \ge P\Bigl(B_j(x) = 1, j=k,\dots,n \Bigr)\ge \texte^{-\sum_{j\ge k}R_j}.
\end{equation}
Since $\xi_{\ge k}^{\GFF}$ has the law of GFF on~$\Lambda_{n-k}$, Lemma~\ref{lemma-5.1} shows that, for any $t > 0$,
\begin{equation}
\label{E:5.30a}
    P\Bigl(\,\max_{x\in\Lambda_n} \xi_{\ge k}^{\GFF}(x) \ge m_{n-k} - t \Bigr) \ge 1 - \frac{1}{a}\texte^{-a\sqrt \beta\, t}.
\end{equation}
Writing~$X_k$ for the maximizer of $x\mapsto \xi_{\ge k}^{\GFF}(x)$ that is minimal in a natural ordering of~$\Lambda_n$, we then have 
\begin{equation}
\begin{aligned}
P\Bigl(\,\max_{x\in\Lambda_n}\xi_{\ge k}^{\DG}(x)< &\,m_{n-k}-t\Bigr) 
         \\
    &\le P\Bigl(\,\max_{x\in\Lambda_n}\xi_{\ge k}^{\GFF}(x)<m_{n-k}-t \Bigr) + P\Bigl(\xi_{\ge k}^{\DG}(X_k)\ne \xi_{\ge k}^{\GFF}(X_k)\Bigr).
\end{aligned}
\end{equation}
Using \eqref{E:5.29} along with the fact that  $X_k$ is independent of the Bernoulli's as it is determined by the GFF increments, we have 
\begin{equation}
    P\Bigl(\xi_{\ge k}^{\DG}(X_k)\ne \xi_{\ge k}^{\GFF}(X_k)\Bigr) \le 1 - \texte^{-\sum_{j=k}^\infty R_j} \le \sum_{j=k}^\infty R_j.
\end{equation}
Then
\begin{equation}
    \label{E:5.37}
     P\Bigl(\,\max_{x\in\Lambda_n}\xi_{\ge k}^{\DG}(x)< m_{n-k}-t\Bigr) \le \frac{1}{a} \texte^{-a \sqrt \beta\, t} + \sum_{j=k}^\infty R_j
\end{equation}
follows with the help of \eqref{E:5.30a}.

Now for any $x\in\Lambda_n$ and any $\lambda > 0$, the standard estimate for the Gaussian distribution gives
\begin{equation}
\label{E:5.34}
     P\biggl(\, \sum_{j=0}^{k-1} \xi^{\GFF}_j(x) < -\lambda \biggr)\le  \texte^{-\frac{\beta \lambda^2 }{2k}}.
\end{equation}
Let~$Y_k$ be the maximizer of $x\mapsto \xi_{\ge k}^{\DG}(x)$ that is again minimal in the natural ordering of~$\Lambda_n$. Since~$Y_k$ is measurable with respect to the increments ``above'' level~$k-1$  and thus independent of $\sum_{j=0}^{k-1} \xi^{\GFF}_j(\cdot)$,  the bound \eqref{E:5.34} applies even with~$x$ replaced by~$Y_k$.
Combining \eqref{E:5.37} and \eqref{E:5.34}  via  the union bound, along with  the almost sure inclusion 
\begin{equation}
\biggl\{\sum_{j=0}^{k-1} \xi^{\DG}_j(x) < -\lambda - Ck\biggr\} \subseteq \biggl\{\sum_{j=0}^{k-1} \xi^{\GFF}_j(x) < -\lambda\biggr\}
\end{equation}
implied by the coupling,   we see 
\begin{equation}
    P\left(\max_{x\in \Lambda_n} \varphi_x^{\DG} < m_{n-k} - t -\lambda -Ck \right) \le  \frac{1}{a}\texte^{-a \sqrt{\beta}\, t} +\sum_{j=k}^\infty R_j +  \texte^{-\frac{\beta \lambda^2 }{2k}}.
\end{equation}
Note that, for $n>k\ge1$,
\begin{equation}
    m_{n-k} = m_n - c_1\,k + c_2\log\left(\frac{n}{n-k}\right)\ge m_n - c_1\,k 
\end{equation}
for some constants $c_1,c_2>0$ that can be gleaned from \eqref{E:1.6a}.
Setting $\lambda := c_1k$ and $t:=c_1 k$ above, we obtain
\begin{equation}
    \label{E:5.43}
     P\left(\max_{x\in \Lambda_n} \varphi_x^{\DG} < m_n - (3c_1+C)k \right) \le \frac{1}{a} \texte^{-a \sqrt{\beta}\, c_1k} + \sum_{j=k}^\infty R_j + \texte^{-\frac{\beta c_1^2}{2} k}.
\end{equation}
By Theorem~\ref{thm-3.4}, $R_k\to0$ exponentially fast when $\beta\in(0, \beta_c)$. Taking~$k\to\infty$ (after $n\to\infty$) thus yields the lower-tail tightness and thus proves the desired claim. 
\end{proofsect}

Next we will prove tightness of the size of level sets:

\begin{proofsect}{Proof of Theorem~\ref{thm-5.2}}
With the help of Proposition~\ref{prop-5.1} we get for~$n$ large enough,
\begin{equation}
\begin{aligned}
        E\bigl[G_n(\lambda) 1_{\{\max_{x\in \Lambda_n} \varphi_x^{\GFF} \le m_n + t\}}\bigr] &= E\left[ \sum_{z\in \Lambda_n} 1_{\{\varphi_z^{\DG} \ge m_n -\lambda\}} 1_{\{\max_{x\in \Lambda_n} \varphi_x^{\GFF}\le m_n + t\}}\right]\\
        &= \sum_{z\in \Lambda_n} P\left(\varphi_z^{\DG} \ge m_n -\lambda,~\max_{x\in \Lambda_n} \varphi_x^{\GFF} \le m_n + t\right)\\
        &\le C''  \bigl(1+(t+\lambda)^2\beta\bigr) \texte^{\sqrt\beta\, c' \lambda}.
     \end{aligned}
\end{equation}
    Then by Markov inequality and Lemma~\ref{lemma-5.1} for $n$ sufficiently large we get
\begin{equation}
\begin{aligned}
P\bigl(G_n(\lambda) \ge u\bigr) &\le P\Bigl(G_n(\lambda) \ge u, \max_{x\in \Lambda_n} \varphi_x^{\GFF} \le m_n + t\Bigr) + P\Bigl(\,\max_{x\in \Lambda_n} \varphi_x^{\GFF} > m_n + t\Bigr)\\
        &\le \frac{1}{u} E\bigl [G_n(\lambda) 1_{\{\max_{x\in \Lambda_n} \varphi_x^{\GFF} \le m_n + t\}}\bigr] + \frac{1}{a} \texte^{-a\sqrt{\beta}\,t}\\
        &\le \frac{1}{u} C''  \bigl(1+(t + \lambda)^2\beta\bigr) \texte^{\sqrt\beta\, c' \lambda} + \frac{1}{a} \texte^{-a \sqrt{\beta}\, t}.
        \end{aligned}
\end{equation}
Taking $u\to \infty$ followed by~$t\to\infty$ now yields the desired conclusion. 
\end{proofsect}


\section{Limit law}
\label{sec6}\noindent
Here we will finally prove our main results concerning the extremal behavior of the hierarchical DG model. As noted earlier, the proofs will be deduced from the corresponding statements about the extremal behavior of the GFF, a.k.a.\ Branching Random Walk.

\subsection{GFF extremal process}
Recall our notation~$P_{n,\beta}'$ for the law of the GFF on~$\Lambda_n$ at inverse temperature~$\beta$ and depth~$n$. We will write~$E_{n,\beta}'$ for the associated expectation. Recall also the notation~$m_n$ for the centering sequence from \eqref{E:1.6a} which, we note, depends on~$\beta$, and the notation~$[x]_n$ from \eqref{E:1.9a} for the image of~$x$ under the natural embedding of~$\Lambda_n$ into the unit interval~$[0,1]$. We then have:

\begin{theorem}
\label{thm-6.1}
There exists an a.s.-finite random Borel measure~$Z$ on~$[0,1]$ with $Z(A)>0$ a.s.\ for all non-empty (relatively) open~$A\subseteq[0,1]$ and a probability measure~$\nu$ on $\MM_\N(\R)$ such that for all $\beta>0$ and all continuous~$f\colon[0,1]\times\R\to[0,\infty)$ with compact support,
\begin{multline}
\label{E:6.1}
E_{n,\beta}'\biggl(\,\exp\Bigl\{-\sum_{x\in\Lambda_n} f\bigl( [x]_n, \varphi_x^{\GFF} - m_n \bigr)\Bigr\}\biggr)
\\\underset{n\to\infty}\longrightarrow\,\,
E\biggl(\exp\Bigl\{-\int Z(\textd x)\otimes\texte^{-\alpha  h}\textd h\otimes\nu(\textd\chi)\bigl[1-\texte^{-\int f(x,\,\beta^{-1/2}(h+\cdot))\textd\chi}\bigr]\Bigr\}\biggr),
\end{multline}
where~$\alpha:=\sqrt{2\log b}$ and where the expectation on the right is over the law of~$Z$.
 Almost every  sample~$\chi$ from~$\nu$ obeys~$\supp(\chi)\subseteq(-\infty,0]$ with~$\sup\supp(\chi)=0$ and~$\chi(\R)=\infty$.
\end{theorem}

For functions not depending on the ``spatial'' coordinate, this is a special case of Madaule~\cite[Theorem~1.1]{Madaule2}. The inclusion of the ``spatial''  locations can in principle be referred to Mallein~\cite{Mallein-clusters} but, for the case at hand, there is a shorter argument that we find more illuminating. Details will be provided in Section~\ref{sec-7.4}.  

Noting that each sample of the process~$\chi$ can be written as~$\sum_{i\ge1}\delta_{t_i}$ for a real-valued sequence~$\{t_i\}_{i\ge1}$, the integral inside the exponential boils down to
\begin{equation}
\int f\bigl(x,\beta^{-1/2}(h+\cdot)\bigr)\textd\chi = \sum_{i\ge1} f\bigl(x,\beta^{-1/2}(h+t_i)\bigr).
\end{equation}
The fact that $\nu$ is a law on the space of Radon measures entails that the sum is finite for~$\nu$-a.e.\ sample~$\chi$. The conditions on the support of a.e.-sample from~$\nu$ in turn ensure that the representation of the limit process is unique.

In the process notation, Theorem~\ref{thm-6.1} can be restated as
\begin{equation}
\label{E:6.3}
\sum_{x\in\Lambda_n}\delta_{[x]_n}\otimes\delta_{\varphi_x^{\GFF}-m_n}
\,\,\,\,\underset{n\to\infty}\lawarrow\,\,\,\,
\sum_{i,j\ge1}\delta_{x_i}\otimes\delta_{\beta^{-1/2}(h_i+t_j^{(i)})},
\end{equation}
where $\{(x_i,h_i)\}_{i\ge1}$ enumerates the points in a sample from the Poisson point process
\begin{equation}
\text{\rm PPP}\bigl(Z(\textd x)\otimes \texte^{-\alpha h}\textd h\bigr)
\end{equation}
and $\{t_j^{(i)}\}_{j\ge1}$ enumerates the points in the $i$-th member of the sequence $\{t^{(i)}\}_{i\ge1}$ of i.i.d.\ samples from~$\nu$ that are independent of $\{(x_i,h_i)\}_{i\ge1}$. Observe that neither~$Z$ nor~$\nu$ depend on~$\beta$ which reflects on the GFF at inverse temperature~$\beta$ being just a~$\beta^{-1/2}$-multiple of the GFF at unit inverse temperature. 

The weak convergence of probability laws \eqref{E:6.3} is relative to the vague topology on Radon measures on~$[0,1]\times\R$ which is why \eqref{E:6.1}  is stated  for continuous functions with compact support. Unfortunately, this restriction is too strict for our purposes. We therefore state and prove:

\begin{lemma}
\label{lemma-6.2}
Let $f\colon[0,1]\times\R\to[0,\infty)$ be continuous and such that $f(x,h)\le A\texte^{-a h^2}$, for some~$A,a>0$ and all~$x\in[0,1]$ and~$h\in\R$. Then \eqref{E:6.1} holds.
\end{lemma}

\begin{proofsect}{Proof}
Let~$\rho\colon[0,\infty)\to[0,1]$ be non-increasing and continuous with $\rho = 1$ on~$[0,1]$ and~$\rho=0$ outside~$[0,2]$. Let~$f_r(x,h):=f(x,h)\rho(|h|/r)$ and observe that~$f_r\le f$ with~$f_r$ increasing monotonically to~$f$ as~$r$ increases to infinity. 

Next let us write
\begin{equation}
G_n'(\lambda):=\bigl|\{x\in\Lambda_n\colon\varphi^{\GFF}_x\ge m_n-\lambda\}\bigr|
\end{equation}
 for the size of an extremal GFF level set and note that the same argument as in the proof of Theorem~\ref{thm-5.2} using Lemma ~\ref{lemma-5.3} gives
\begin{equation}
\label{E:6.5}
E_{n,\beta}'\bigl[G_n'(\lambda) 1_{\{\max_{x\in \Lambda_n} \varphi_x^{\GFF} \le m_n + t\}}\bigr]\le C'\bigl(1+(t+\lambda)^2 \beta\bigr) \,\texte^{\sqrt \beta \,c'\lambda}
\end{equation}
for large enough~$n$ and~$\lambda>0$. Noting $|f(x,h)-f_r(x,h)|\le A\texte^{-a h^2}1_{\{|h|>r\}}$, for all~$\epsilon>0$ and~$r> t>0$ a simple use of Markov inequality yields
\begin{equation}
\label{E:6.6}
\begin{aligned}
P_{n,\beta}'\biggl(\,\Bigl|\sum_{x\in\Lambda_n} &f\bigl( [x]_n, \varphi_x^{\GFF} - m_n \bigr)-\sum_{x\in\Lambda_n} f_r\bigl( [x]_n, \varphi_x^{\GFF} - m_n \bigr)\Bigr|>\epsilon\biggr)
\\
&\le P_{n,\beta}'\Bigl(\,\max_{x\in \Lambda_n} \varphi_x^{\GFF} > m_n + t\Bigr)
\\
&\qquad+ A\epsilon^{-1}\,E_{n,\beta}'\biggl(\,\sum_{x\in\Lambda_n}\texte^{-a(\varphi_x^{\GFF}-m_n)^2}1_{\{\varphi_x^{\GFF} < m_n-r\}}1_{\{\max_{x\in \Lambda_n} \varphi_x^{\GFF} \le m_n + t\}}\biggr),
\end{aligned}
\end{equation}
where the values with~$\varphi_x^{\GFF} > m_n+r$ do not enter due to our restriction on the maximum and the assumption that~$r>t$. Assuming~$r$ to be a natural, the last expectation is bounded using $1_{\{\varphi_x^{\GFF} < m_n-r \}} = \sum_{\ell=r+1}^\infty 1_{\{-\ell \le \varphi_x^{\GFF} -m_n < -\ell+1 \}}$ by
\begin{equation}
\label{E:6.7}
\begin{aligned}
\sum_{\ell=r+1}^\infty \texte^{-a \ell^2}\,E_{n,\beta}'\bigl[G_n'(\ell) &1_{\{\max_{x\in \Lambda_n} \varphi_x^{\GFF} \le m_n + t\}}\bigr]
\\
&\qquad\le \sum_{\ell=r+1}^\infty C'\bigl(1+(t+\ell)^2 \beta\bigr) \,\texte^{\sqrt \beta c'\ell -a\ell^2},
\end{aligned}
\end{equation}
where \eqref{E:6.5} was used to get the inequality.
 
The sum on the right of \eqref{E:6.7} vanishes as~$r\to\infty$ and, taking~$t\to\infty$ afterwards with the help of Lemma~\ref{lemma-5.1}, so does the probability on the left of \eqref{E:6.6}. Hence we get
\begin{equation}
\lim_{r\to\infty}\limsup_{n\to\infty}\biggl|\,E_{n,\beta}'\bigl(\texte^{-\sum_{x\in\Lambda_n} f( [x]_n, \varphi_x^{\GFF} - m_n )}\bigr)-E_{n,\beta}'\bigl(\texte^{-\sum_{x\in\Lambda_n} f_r( [x]_n, \varphi_x^{\GFF} - m_n )}\bigr)\biggr|=0.
\end{equation}
Since~$f_r$ does have compact support, \eqref{E:6.1} applies to it. It then suffices to show that the right-hand side of \eqref{E:6.1} for~$f_r$ tends to that for~$f$ as~$r\to\infty$. Thanks to the monotonicity of~$r\mapsto f_r$, this follows by the Monotone Convergence Theorem for the inner integral and then the Bounded Convergence Theorem for the expectation.
\end{proofsect}

\subsection{Conversion to GFF extremal process}
Next we will develop an argument that reduces observables involving the extremal process of the DG model to those of a GFF. This is what we will need to deduce Theorem~\ref{thm-1.2} directly from Theorem~\ref{thm-6.1}. 

Pick $\beta\in(0,\beta_\cc)$ and a continuous function ~$f\colon[0,1]\times\R\to[0,\infty)$ with compact support; we will keep these fixed throughout the rest of this subsection. We are interested in the $n\to\infty$ weak limit of the random variable
\begin{equation}
Y_n:=\sum_{x\in\Lambda_n } f\bigl( [x]_n, \varphi_x^{\DG} - \lfloor m_n\rfloor \bigr).
\end{equation}
We will work in the coupling from Theorem~\ref{thm-4.1}. Pick~$k$ with~$1\le k<n$ and recall the notation $\xi_{\ge k}^{\GFF}(x)$ for the sum $\sum_{j=k}^{n} \xi_j^{\GFF}(x)$ which we for later convenience simplify as follows: Given~$z\in\Lambda_{n-k}$, pick any~$x\in\Lambda_n$ with $m^k(x)=z$ and abbreviate
\begin{equation}
\varphi_z^{\GFF}:=\xi_{\ge k}^{\GFF}(x).
\end{equation}
Similarly, denoting $\xi_{\ge k}^{\DG}(x):=\sum_{j=k}^{n}\xi_j^{\DG}(x)$, let
\begin{equation}
\varphi_z^{\DG}:=\xi_{\ge k}^{\DG}(x).
\end{equation}
The definition of the coupling ensures that these objects do not depend on the choice of~$x$ and so we may think of them as indexed by~$\Lambda_{n-k}$. Note that~$\{\varphi_z^{\DG}\colon z\in\Lambda_{n-k}\}$ is continuously distributed and $\{\varphi_z^{\GFF}\colon z\in\Lambda_{n-k}\}$ is a GFF on~$\Lambda_{n-k}$.

Using these objects we now define $g_k\colon[0,1]\times\R\to[0,\infty)$ via
\begin{equation}
\label{E:1.1e}
\texte^{-g_k(v,h)}:=E\Bigl(\texte^{-\sum_{x\in\Lambda_k(z)}f(v,\varphi_x^{\DG}-\varphi_z^{\DG}+h)}\,\Big|\,\varphi_z^{\GFF} = \varphi_z^{\DG} = h\Bigr),
\end{equation}
where $\Lambda_k(z):=\{x\in\Lambda_n\colon m^k(x)=z\}$ and where the conditioning on the (apparently singular) event $\varphi_z^{\GFF} = \varphi_z^{\DG} = h$ is to be interpreted as setting an initial value in the tree-indexed Markov chain defining the coupling.
The  technical aspects of the desired conversion are the content of: 

\begin{proposition}
\label{prop-6.2}
We have
\begin{equation}
\label{E:6.11}
\lim_{k\to\infty}\limsup_{n\to\infty}\,\Biggl|\,
E_{n,\beta}(\texte^{-Y_n}) - E_{n-k,\beta}'\biggl(\exp\Bigl\{-\sum_{z\in\Lambda_{n-k}} g_k\bigl([z]_{n-k},\varphi_z^{\GFF} - \lfloor m_{n}\rfloor\bigr)\Bigr\}\biggr)\Biggr|=0.
\end{equation}
\end{proposition}

The proof of Proposition~\ref{prop-6.2} will require a few lemmas. Let~$k$ be a natural with $1\le k<n$ and let $\tau(x)$ be  the random ``time'' defined in \eqref{E:5.10}. We then modify~$Y_n$ into
\begin{equation}
Y'_{n,k}:=\sum_{x\in\Lambda_n } f\bigl( [m^k(x)]_{n-k}, \varphi_x^{\DG} - \lfloor m_n \rfloor\bigr) 1_{\{\tau(x) < k\} }.
\end{equation}
As we show next, these modifications are not very significant provided~$n$ and~$k$ are taken sufficiently large:

\begin{lemma}
\label{lemma-6.4}
For each~$\epsilon>0$,
\begin{equation}
\lim_{k\to\infty}\,\limsup_{n\to\infty}\,P\bigl(|Y_n-Y_{n,k}'|>\epsilon\bigr)=0,
\end{equation}
where $P$ denote the coupling law from  Theorem~\ref{thm-4.1}. 
\end{lemma}

\begin{proofsect}{Proof}
Let~$Y''_{n,k}$ be given by the same formula as~$Y'_{n,k}$ except with the first argument of~$f$ replaced by~$[x]_n$. Let~$\lambda>0$ be such that~$f$ is supported in~$[0,1]\times[-\lambda,\lambda]$ and let~$\omega(r):=\sup\{|f(x,h)-f(y,h)|\colon |x-y| < r,\,h \in [-\lambda, \lambda]\}$ be the modulus of continuity of~$f$ in the first variable. In light of~$|[m^k(x)]_{n-k}-[x]_n| < b^{-n+k}$ we then have
\begin{equation}
|Y_{n,k}'-Y_{n,k}''|\le \omega(b^{-n+k})G_n(\lambda)
\end{equation}
which tends to zero in probability as~$n\to\infty$ thanks to Theorem~\ref{thm-5.2}.

It thus suffices to prove the claim with~$Y_{n,k}'$ replaced by~$Y_{n,k}''$. Using the coupling measure~$P$ and writing~$A$ for a constant that bounds~$|f|$ uniformly,  we get
\begin{equation}
P\bigl(|Y_n-Y_{n,k}''|>\epsilon\bigr)
\le P\biggl(\,\sum_{x\in\Lambda_n}1_{\{\varphi_x^{\DG} \ge \lfloor m_n \rfloor-\lambda\}}1_{\{\tau(x) \ge k\}}>\epsilon/A\biggr).
\end{equation}
The probability on the right is bounded via Markov inequality  by 
\begin{multline}
\label{E:6.18}
\quad P\Bigl(\,\max_{x\in \Lambda_n} \varphi_x^{\GFF} > m_n + t\Bigr)
\\+
 A\epsilon^{-1}\sum_{x\in\Lambda_n}P\Bigl(\varphi_x^{\DG} \ge \lfloor m_n \rfloor-\lambda,\,\tau(x) \ge k,
\,\max_{x\in \Lambda_n} \varphi_x^{\GFF}\le m_n+t\Bigr).
\quad
\end{multline}
Similarly as in~\eqref{E:6.19u},  the probability under the sum can be bounded by splitting the event $\{\varphi_x^{\DG} \ge \lfloor m_n \rfloor - \lambda\}$ into three events. For each of these three events, we need to show the same type of bound as in right side of \eqref{E:5.21} with additional exponential decay in~$k$ induced by $\{\tau(x) \ge k\}.$ 

The first event only involves GFF so the exponential decay in~$k$ from $\{\tau(x) \ge k\}$ can be extracted immediately from the independence of the GFF and the Bernoulli's ensured by Theorem~\ref{thm-4.1}(2). For the third event, the probability vanishes already in the limit as $n\to \infty$ due to the $n^{-1}$ decay facilitated by the choice of $ r$ in \eqref{E:5.23};  in particular, no rate of  decay in~$k$ is  needed.  For the second event, we follow the estimates up to \eqref{E:5.25} where instead of the second probability on the right we get
\begin{equation}
P\biggl(\,\sum^{n}_{ i=0}   1_{\{ B_i(z) = 0\}} > \frac{\ell+j}{C},\,\tau( z) \ge k\biggr).
\end{equation}
 Here we denote $\eta:=\frac{\ell+j}{C}$ and observe  that, for $\eta$ large enough,  Lemma~\ref{lemma-5.4} gives
\begin{equation}
\begin{aligned}
        P\biggl(\,\sum^{n}_{ i=0} &  1_{\{ B_i( z) = 0\}} > \eta,\, \tau( z) \ge k \biggr)
\\
&\le P\biggl(\,\sum^{k-1}_{ i=0}  1_{\{ B_i( z) = 0\}} > \frac{\eta}{2}\biggr)\,P\bigl(   \tau(z)  \ge k\bigr)
+ P\biggl(\,\sum^{n}_{ i=k}  1_{\{ B_i( z) = 0\}} > \frac{\eta}{2}\biggr )
\\
&\le \exp\biggl\{-\frac{1}{2} \frac{\eta}{2} \log\Bigl(\frac\eta2\Bigr)\biggr\} \biggl(1 - \prod_{ i \ge k} \texte^{- R_i}\biggr)
\\
&\qquad\qquad+ \exp\Biggl\{\frac{\eta}{2} + \frac{\eta}{2}\log \biggl(\, \sum^{n}_{ i=k} (1 - \texte^{- R_i})\biggr) - \frac{\eta}{2} \log\Bigl(\frac\eta2\Bigr)\Biggr\}.
\end{aligned}
\end{equation}
 Both terms above exhibit simultaneously super-exponential decay in~$\eta$ as well as exponential decay in~$k$ and so the probability on the right of \eqref{E:6.18} vanishes in the limit as $n\to \infty$ and then as $k\to \infty$. Taking $t\to \infty$, we get Lemma~\ref{lemma-6.4}.
\end{proofsect}

We also need the following ``reverse-Jensen'' type of inequality:

\begin{lemma}
\label{lemma-6.5}
Given a sequence $\{p_k\}_{k=0}^n$ of positive probabilities, let $\{Y_x\colon x\in\bigcup_{k=0}^n\Lambda_k\}$ be independent zero-one valued random variables such that $\BbbP(Y_x = 1) =p_k$ for all~$k=0,\dots,n$ and all $x\in\Lambda_k$. Writing $x\in\Lambda_n$ as~$x=(x_1,\dots,x_n)$, set $Z_x:=\prod_{i=0}^n Y_{x_i}$, where $x_0 := \rho$ denotes the root of the tree. Then for all~$h\colon\Lambda_n\to[0,\infty)$ and all~$\lambda\in(0,1)$,
\begin{equation}
\E\biggl(\,\exp\Bigl\{-\sum_{x\in\Lambda_{n}}h(x)Z_x\Bigr\}\biggr)\le
\frac1{\lambda^2}\frac{1-q}{q}+\exp\Bigl\{-(1-\lambda) q\sum_{x\in\Lambda_{n}}h(x)\Bigr\},
\end{equation}
where $q:=\prod_{k=0}^n p_k$. 
\end{lemma}

\begin{proofsect}{Proof}
Assume $h$ is not identically zero for otherwise there is nothing to prove. We will use that the Chebyshev inequality for  a  non-negative random variable~$X$ gives
\begin{equation}
\BbbP\Bigl(X\le(1-\lambda)\E X\Bigr)\le\BbbP\Bigl(\bigl|X-\E X\bigr|\ge \lambda\E X\Bigr)\le \frac1{\lambda^2}\frac{\Var(X)}{[\E X]^2}
\end{equation}
whenever~$\E X>0$.  We will apply this to $X:=\sum_{x\in\Lambda_{n}}h(x)Z_x$. The first moment calculation results in 
\begin{equation}
\E\Bigl(\sum_{x\in\Lambda_{n}}h(x)Z_x\Bigr)=q\sum_{x\in\Lambda_{n}}h(x).
\end{equation}
 To compute the variance of~$X$ we note that, for any~$x,y\in\Lambda_{n}$ with the nearest common ancestor in~$\Lambda_{n-\ell}$, we have 
\begin{equation}
\Cov(Z_x,Z_y)
=\Bigl(\,\prod_{k=0}^{\ell-1}  p_k\Bigr)^2\prod_{i=\ell}^n p_i -  \Bigl(\prod_{k=0}^n  p_k\Bigr)^2 = q^2 \biggl(\Bigl(\prod_{k=\ell}^n p_k\Bigr)^{-1}-1\biggr)
 \le q^2[q^{-1}-1]. 
\end{equation}
 Hence we get 
\begin{equation}
\Var\biggl(\sum_{x\in\Lambda_{n}}h(x)Z_x\biggr)
=\sum_{x,y\in\Lambda_{n}}h(x)h(y)\Cov(Z_x,Z_y)
\le
q^2\,\frac{1-q}{q}\Bigl(\,\sum_{x\in\Lambda_{n}}h(x)\Bigr)^2.
\end{equation}
The above Chebyshev bound then  yields 
\begin{equation}
\BbbP\biggl(\,\sum_{x\in\Lambda_{n}}h(x)Z_x
\le (1-\lambda) q\sum_{x\in\Lambda_{n}}h(x)\biggr)
 \le \frac1{\lambda^2}\frac{1-q}{q}, 
\end{equation}
 which then  readily implies the claim.
\end{proofsect}

We are now fully equipped to give:

\begin{proofsect}{Proof of Proposition~\ref{prop-6.2}} 
We will work in the coupling measure~$P$ throughout. With the help of  Lemma~\ref{lemma-6.4}, it suffices to prove \eqref{E:6.11} where $E_{n,\beta}(\texte^{-Y_n})$ is replaced by $E(\texte^{-Y_{n,k}'})$. Denote
\begin{equation}
\FF_k:=\sigma\Bigl(\xi_j^{\DG}(x),\xi_j^{\GFF}(x),B_j(x)\colon j\ge k,\, x\in\Lambda_n\Bigr).
\end{equation}
 The Markovian definition of the coupling along with the rewrite 
\begin{equation}
\varphi_x^{\DG} - \lfloor m_n\rfloor  = \varphi_x^{\DG}- \varphi_z^{\GFF} + \varphi_z^{\GFF}-\lfloor m_n\rfloor
\end{equation}
then casts~$E(\texte^{-Y'_{n,k}}|\FF_k)$ as
\begin{equation}
\label{E:6.14}
\prod_{z\in\Lambda_{n-k}}E\biggl(\exp\Bigl\{-\!\!\sum_{x\in\Lambda_k(z)} 1_{\{\tau(x) < k\}}f\bigl([z]_{n-k},\varphi_x^{\DG}- \varphi_z^{\GFF} + \varphi_z^{\GFF}-\lfloor m_n\rfloor\bigr)\Bigr\}\bigg|\,\FF_k\biggr).
\end{equation}
The event~$\{\tau(x) < k\}$ is~$\FF_k$ measurable and so it either occurs for all~$x\in\Lambda_k(z)$ or none of them. In the former case we have~$\varphi_z^{\DG}=\varphi_z^{\GFF}$. Writing~$\tilde g_{n,k}$ for the function defined by
\begin{equation}
\label{E:1.1f}
\texte^{-\tilde g_{n,k}(v,h)}:=E\Bigl(\texte^{-\sum_{x\in\Lambda_k(z)}f(v,\varphi_x^{\DG}-\varphi_z^{\DG}+h)} \,\Big|\,\varphi_z^{\GFF} = \varphi_z^{\DG} = h+\lfloor m_n\rfloor\Bigr) ,\end{equation}
the expectation in \eqref{E:6.14} can thus be contracted to
\begin{equation}
E(\texte^{-Y'_{n,k}}|\FF_k) = 
\exp\Bigl\{-\sum_{z\in\Lambda_{n-k}} \tilde g_{n,k}\bigl([z]_{n-k},\varphi_z^{\GFF} - \lfloor m_{n}\rfloor\bigr)1_{F(z)}\Bigr\},
\end{equation}
where~$F(z):=\bigcap_{j=k}^{n}  \{ B_{j}(x)=1\}$ for some (and any) $x\in\Lambda_n$ such that~$z=m^k(x)$.

Next observe a key fact that the law of~$\varphi_x^{\DG}-\varphi_z^{\DG}$ ``initiated'' from~$\varphi_z^{\DG}=m+h$  is the same for all~$m\in\Z$ due to the 1-periodicity of the effective potential $\{v_k\}_{k\ge 0}$. This implies
\begin{equation}
\tilde g_{n,k}=g_k
\end{equation}
and so
\begin{equation}
\label{E:7.33a}
\begin{aligned}
E(\texte^{-Y'_{n,k}})
&=E \bigl(E(\texte^{-Y'_{n,k}}|\FF_k)\bigr)
\\
&= E \biggl(\exp\Bigl\{-\sum_{z\in\Lambda_{n-k}} \tilde g_{n,k}\bigl([z]_{n-k},\varphi_z^{\GFF} - \lfloor m_{n}\rfloor\bigr)1_{F(z)}\Bigr\}\biggr)
\\
&\ge E\biggl(\exp\Bigl\{-\sum_{z\in\Lambda_{n-k}} g_k\bigl([z]_{n-k},\varphi_z^{\GFF} - \lfloor m_{n}\rfloor\bigr)\Bigr\}\biggr),
\end{aligned}
\end{equation}
where we used~$g_k\ge0$ to drop the indicators of~$F(z)$. 

For the complementary upper bound we note that the random variable~$1_{F(z)}$ can be thought of as product of tree-indexed Bernoulli's on the path from~$z$ to the tree root. For each~$\lambda\in(0,1)$, Lemma~\ref{lemma-6.5} applied to $E(\texte^{-Y'_{n,k}}|\FF_k)$ along with Jensen's inequality gives
\begin{equation}
\label{E:7.34a}
E(\texte^{-Y'_{n,k}})
\le \frac1{\lambda^2}\frac{1- q}{ q}+\Biggl[\,E\biggl(\exp\Bigl\{-\sum_{z\in\Lambda_{n-k}} g_k\bigl([z]_{n-k},\varphi_z^{\GFF} - \lfloor m_{n}\rfloor\bigr)\Bigr\}\biggr)\Biggr]^{(1-\lambda) q},
\end{equation}
 where $q=q_{n,k}:=\exp\{- \sum_{j=k}^{n} R_j\}$. 
Setting $\lambda:=(1- q)^{1/3}$, the first term on the right tends to zero  as $n\to\infty$ and  $k\to\infty$ while the exponent $(1-\lambda) q$ tends to one. Taking  these limits  in~\eqref{E:7.34a}, combined with the lower bound in \eqref{E:7.33a}, we get~\eqref{E:6.11}.  
\end{proofsect}

\subsection{Proofs of the main results}
We are now finally in a position to give formal proofs of our main results. We start with the convergence of the extremal process.

\begin{proofsect}{Proof of Theorem~\ref{thm-1.2}}
Let~$\lambda>0$ be such that $\supp(f)\subseteq[0,1] \times [-\lambda, \lambda]$.  The function~$g_k$ defined in \eqref{E:1.1e} is bounded and, thanks to the continuity of~$f$ and~$v_k$, continuous. 
While~$g_k$ does not have compact support, we claim that it has Gaussian tails in the second variable. To see this note that $f(v,\varphi_x^{\DG}-\varphi_z^{\DG}+h) > 0$ implies $\varphi_x^{\DG}-\varphi_z^{\DG}+h \in [-\lambda ,\lambda]$. Then $\varphi_x^{\GFF}-\varphi_z^{\GFF}+h \in [-\lambda - Ck,\lambda + Ck]$ a.s.\ for~$C$ as in Theorem~\ref{thm-4.1}(4). Since $\varphi_x^{\GFF}-\varphi_z^{\GFF}+h$ is independent of~$\varphi^{\GFF}_z$ with law~$\NN(h, k/\beta)$, for any $x\in \Lambda_k(z)$ the standard Gaussian estimate yields 

\begin{equation}
\begin{aligned}
P\Bigl( f(v,\varphi_x^{\DG}-\varphi_z^{\DG}+&h) > 0 \,\bigg |\, \varphi_z^{\GFF} = \varphi_z^{\DG} = h\Bigr) 
   \\
    &\le P\Bigl(\bigl|\NN(0, k/\beta)\bigr|\ge |h| - \lambda -Ck\Bigr) \le 2 \texte^{-\frac{\beta ( |h| - \lambda - Ck)^2 }{2k}} 
\end{aligned}
\end{equation} 
whenever $|h| \ge \lambda + Ck$. Since~$\Lambda_k(z)$ has cardinality~$b^k$, the union bound then shows 
\begin{equation}
1-\texte^{-g_k(v,h)}\le 2 b^k \texte^{-\frac{\beta (|h| - \lambda - Ck)^2 }{2k}}.
\end{equation}
It follows that $g_k(v,h)$ has Gaussian tails. 

With~$g_k$ conforming to the condition of Lemma~\ref{lemma-6.2} (which feeds into Theorem~\ref{thm-6.1}), for any sequence $n_j\to\infty$ such that
\begin{equation}
\label{E:1.7}
s:=\lim_{j\to\infty}\bigl(m_{n_j}-\lfloor m_{n_j}\rfloor\bigr)\,\,\text{\rm\ exists},
\end{equation}
we claim to get
\begin{equation}
\begin{aligned}
\label{E:1.3}
&E_{n_j-k,\beta}'\biggl(\exp\Bigl\{-\sum_{z\in\Lambda_{n_j-k}} g_k\bigl([z]_{n_j-k},\varphi_z^{\GFF} - \lfloor m_{n_j}\rfloor
\bigr)\Bigr\}\biggr)
\\&
\quad\underset{j\to\infty}\longrightarrow\,\,
E\biggl(\exp\Bigl\{-\int Z(\textd x)\otimes\texte^{-\alpha h}\textd h\otimes\nu(\textd\chi)\bigl(1-\texte^{-\int g_k(x,\,s+\beta^{-1/2}(h-\alpha k+\cdot))\,\textd\chi})\Bigr\}\biggr).
\end{aligned}
\end{equation}
This would follow from Theorem~\ref{thm-6.1} (or, more precisely, Lemma~\ref{lemma-6.2}) if $\varphi_z^{\GFF} - \lfloor m_{n_j}\rfloor$ were replaced by $\varphi_z^{\GFF} - m_{n_j-k}+s-\beta^{-1/2} \alpha k$, so a key point to address here is the implicit limit $m_{n_j-k}-s+\beta^{-1/2} \alpha k-\lfloor m_{n_j}\rfloor\to0$ inside the second argument of~$g_k$. This is handled by setting $g_{k,\delta}^-(v,h)$, resp, $g_{k,\delta}^+(v,h)$ to be the infimum, resp., supremum of $g_k(v,h')$ over $h'\in[h-\delta,h+\delta]$ and noting that then
\begin{equation}
\label{E:7.40i}
\begin{aligned}
g_{k,\delta}^-\bigl([z]_{n_j-k},\varphi_z^{\GFF} -  & m_{n_j-k}+s-\beta^{-1/2} \alpha k\bigr)
\\
&\le g_k\bigl([z]_{n_j-k},\varphi_z^{\GFF} - \lfloor m_{n_j}\rfloor\bigr)
\\
&\qquad\qquad\le g_{k,\delta}^+\bigl([z]_{n_j-k},\varphi_z^{\GFF} - m_{n_j-k}+s-\beta^{-1/2} \alpha k\bigr)
\end{aligned}
\end{equation}
as soon as $j$ is so large that $|m_{n_j-k}-s+\beta^{-1/2} \alpha k-\lfloor m_{n_j}\rfloor|\le\delta$. As~$g_{k,\delta}^\pm$ are both continuous with a Gaussian upper bound, Lemma~\ref{lemma-6.2} implies sub-subsequential convergence of the type \eqref{E:1.3} with~$g_k$ replaced by the corresponding~$g_{k,\delta}^\pm$ in~\eqref{E:7.40i} and thus bound the \textit{limes superior}, resp., \textit{inferior} of the expectations on the left of \eqref{E:1.3} as~$j\to\infty$ by the right-hand side with~$g_k$ replaced by~$g_{k,\delta}^+$, resp., $g_{k,\delta}^-$. Using that $g_{k,\delta}^\pm\to g_k$ monotonically as~$\delta\downarrow0$, \eqref{E:1.3} is concluded with the help of the Dominated Convergence Theorem. 

With \eqref{E:1.3} in place, a simple change of variables allows us to get rid of the shift by~$s$ at the cost of multiplying the integral by~$\texte^{\alpha\sqrt\beta\, s}$.
Next note that each sample~$\chi$ from~$\nu$ can be written as
\begin{equation}
\label{E:6.35}
\chi = \sum_{i\ge1}\delta_{t_i}
\end{equation}
for a sequence of real-valued points~$t_1\ge t_2\ge\dots$. 
The definition of~$g_k$ then shows
\begin{equation}
\label{E:6.37}
\begin{aligned}
&\texte^{-\int g_k(v,\beta^{-1/2}(h+ \cdot))\textd\chi}
=\prod_{i\ge1}\texte^{- g_k(v,\beta^{-1/2}(h+ t_i))} 
\\&\qquad=\prod_{i\ge1}E\Bigl(\,\texte^{-\sum_{x\in\Lambda_k(z)} f(v,\varphi^{\DG}_x)}\,\Big|\,\varphi_z^{\GFF} = \varphi_z^{\DG} = \beta^{-1/2}(h+ t_i)\Bigr),
\end{aligned}
\end{equation}
where the conditional event was used to bring in the second argument of~$f$ to the stated form. Observe that the second argument of~$f$ receives an integer value, as it should  because the DG model will ultimately only see the restriction of~$f$ to integers. 

We now want to cast the right-hand side of \eqref{E:6.37} as exponential of an integral but for that we will need to borrow part of the integral with respect to~$h$ from \eqref{E:1.3}. Writing the conditional expectation on the right of \eqref{E:6.37} as $F_k(\beta^{-1/2}(h+ t_i))$, elementary manipulations with integrals show
\begin{equation}
\label{E:7.43}
\begin{aligned}
\int\textd h\,\texte^{-\alpha h}\Bigl(&1-\prod_{i\ge1}F_k\bigl(\beta^{-1/2}(h+ t_i)\bigr)\Bigr)
\\
&=\sqrt\beta\int\textd h\,\texte^{-\alpha\sqrt\beta\, h}\Bigl(1-\prod_{i\ge1}F_k\bigl(h+\beta^{-1/2} t_i\bigr)\Bigr)
\\
&=\sqrt\beta\sum_{n\in\Z}\texte^{-\alpha\sqrt\beta\, n}\int_0^1\textd u\,\texte^{-\alpha\sqrt\beta\,u}\Bigl(1-\prod_{i\ge1}F_k\bigl(n+u+\beta^{-1/2} t_i\bigr)\Bigr)
\\
&=\alpha^{-1}(1-\texte^{-\alpha\sqrt\beta})\sum_{n\in\Z}\texte^{-\alpha\sqrt\beta\, n}\,E\Bigl(1-\prod_{i\ge1}F_k\bigl(n+U+\beta^{-1/2} t_i\bigr)\Bigr),
\end{aligned}
\end{equation}
where the expectation is with respect to the random variable~$U$ distributed as
\begin{equation}
 P(U\in\textd u) =  1_{[0,1)}(u)\frac{\alpha\sqrt\beta\,\texte^{-\alpha\sqrt\beta\, u}}{1-\texte^{-\alpha\sqrt\beta}}\textd u
\end{equation}
 and independent of~$\{t_i\}_{i\ge1}$. 

In order to rewrite  the expectation in \eqref{E:7.43} further, for each~$h\in\R$ let~$\mu_h$ be the probability measure on~$\Z^{\Lambda_k}$ defined by
\begin{equation}
\mu_h\bigl(\varphi\in\cdot):=P\bigl(\,\varphi^{\DG}\in\cdot\,\big|\,\varphi_z^{\GFF} = \varphi_z^{\DG} = h\bigr),
\end{equation}
where~$z$ denotes the root of~$\Lambda_k$. (Again, we are relying on the Markovian nature of the coupling.)
Now sample~$U$ and~$\{t_i\}_{i\ge1}$ independently from the above laws and then let $\{\varphi^{(i)}\}_{i\ge1}$ be independent with~$\varphi^{(i)}$ sampled from~$\mu_h$ for $h:=U+\beta^{-1/2} (t_i - \alpha k)$. With these samples in hand, let  
$\nu_{\beta,k}$ be the law of
\begin{equation}
\label{E:7.46}
\zeta:=\sum_{i\ge1}\sum_{x\in\Lambda_k}\delta_{\varphi^{(i)}_x}
\end{equation} 
 on the space of  $\N\cup\{+\infty\}$-valued  Borel measures endowed with the $\sigma$-algebra that makes $\zeta(B)$ measurable for each Borel~$B\subseteq\R$. (We will later show that~$\nu_{\beta,k}$ is asymptotically concentrated on~$\MM_\N(\R)$, i.e., the subset of Radon measures.) 

Using the conditional independence of $\{\varphi^{(i)}\}_{i\ge1}$ we now readily check
\begin{equation}
\int\nu(\textd\chi)E\Bigl(1-\prod_{i\ge1}F_k\bigl(n+U+\beta^{-1/2} t_i\bigr)\Bigr)
=\int\nu_{\beta,k}(\textd \zeta)\,\bigl(1-\texte^{-\int  f(x,n+\cdot)\textd\zeta}\bigr)
\end{equation}
thus producing the desired integral form.
Combining the above observations,
the right-hand side of \eqref{E:1.3} becomes
\begin{equation}
\label{E:6.43}
E\,\biggl(\exp\Bigl\{-\tilde c\,\texte^{\alpha\sqrt\beta\, s}\int Z(\textd x)\otimes\nu_{\beta,k}(\textd\zeta)\sum_{n\in\Z}\texte^{-\alpha\sqrt\beta\, n}\bigl(1-\texte^{-\int f(x,n+\cdot)\textd\zeta})\Bigr\}\biggr),
\end{equation}
where
\begin{equation}
\tilde c:=\alpha^{-1}(1-\texte^{-\alpha\sqrt\beta})
\end{equation}
has no dependency on~$s$.

Our next task is to extract a weak (subsequential) limit as~$k\to\infty$.  For this we first need to address the difficulty posed by the fact that the integral in \eqref{E:6.43} is effectively with respect to an infinite measure due to the sum over~$n\in\Z$. We solve this by incorporating the sum into the measure~$\bar\nu_{\beta,k}$ defined by
\begin{equation}
\bar\nu_{\beta,k}(\textd\zeta) := \sum_{n\in\Z}\texte^{\alpha\sqrt\beta\,n}\nu_{\beta,k}\circ\tau_n^{-1}(\textd\zeta),
\end{equation}
where~$\tau_n$ is the shift by~$n$ acting on elements of~$\MM_\N(\Z)$ as~$\tau_n(\zeta)(A):=\zeta(n+A)$ and so on integrals as~$\int f\textd \tau_n(\zeta) = \int f(-n+\cdot)\textd\zeta$. Note that this gives
\begin{equation}
\label{E:7.51b}
\bar\nu_{\beta,k}\circ\tau_n^{-1} = \texte^{-\alpha\sqrt\beta\,n}\,\bar\nu_{\beta,k}
\end{equation}
for each~$n\in\Z$.
With this notation \eqref{E:6.43} becomes
\begin{equation}
\label{E:7.51a}
E\,\biggl(\exp\Bigl\{-\tilde c\,\texte^{\alpha\sqrt\beta\, s}\int Z(\textd x)\otimes\bar\nu_{\beta,k}(\textd\zeta)\bigl(1-\texte^{-\int f(x,\cdot)\textd\zeta})\Bigr\}\biggr),
\end{equation}
where \MB the implicit use of Tonelli's theorem will be justified by showing that, given integer-valued $r\le u$ there is $k_0\ge1$ such that
\begin{equation}
\label{E:7.53w}
\sup_{k\ge k_0}\bar\nu_{\beta,k}\bigl(\zeta([r,u])\ge1\bigr)<\infty.
\end{equation}
We will show this by the same argument that proves that~$\{\bar\nu_{\beta,k}\}_{k\ge }$ \eMB are tight relative to the vague topology. This will follow once we show that, for any~$r\in\Z$,
\begin{equation}
\label{E:6.45}
\lim_{u\to\infty}\limsup_{k\to\infty}\,\bar\nu_{\beta,k} \Bigl(\zeta\bigl([r,u]\bigr)>1/\epsilon\Bigr)\,\,\underset{\epsilon\downarrow0}\longrightarrow0.
\end{equation}
\MB (The limit $u\to\infty$ exists by monotonicity.) \eMB
 To prove this, take a continuous function~$f\colon[0,1]\times\R\to[0,1]$ with~$f(x,h)\ge 1_{ [r,u]}(h)$ and $\supp(f)\subseteq [0,1]\times (r-1,\infty)$. The integral in \eqref{E:6.43} with~$f$ replaced by~$\epsilon f$ is then at least $Z([0,1])(1-\texte^{-1})\bar\nu_{\beta,k}(\zeta([r,u])>1/\epsilon)$. 
For~$Y_n$ associated with~$f$, Proposition~\ref{prop-6.2} along with \eqref{E:1.3} give 
\begin{equation}
\label{E:7.51u}
 \liminf_{n\to\infty} E_{n,\beta}\bigl(\texte^{-\epsilon Y_{n}}\bigr)\le  
 \liminf_{k\to\infty} E\bigl(\texte^{-\tilde c Z([0,1])(1-\texte^{-1})\bar\nu_{\beta,k}(\zeta([r,u])>1/\epsilon)}\bigr).
\end{equation}
\MB Recalling the notation~$G_n(\lambda)$ for the size of the level set from \eqref{E:6.2w}, the conditions on~$f$ imply $Y_n\le G_n(r-2)$. By Theorem~\ref{thm-5.2}, the expectation on the left is thus uniformly positive as~$n\to\infty$. Taking~$\epsilon:=2$, the fact that $Z([0,1])>0$ a.s.\ shows \eqref{E:7.53w}. 
Similarly, taking instead~$u\to\infty$ followed by \eMB $\epsilon\downarrow0$, the right-hand side of \eqref{E:7.51u} tends to one. Since $Z([0,1])>0$ a.s., we get \eqref{E:6.45} as desired.

The  condition \eqref{E:6.45}  permits us to extract a  subsequence~$\{\bar\nu_{\beta,k_j}\}_{j\ge0}$ that converges relative to the vague topology to a measure~$\bar\nu_\beta'$. However, since the integrand in \eqref{E:7.51a} vanishes when~$\zeta(\supp(f))=0$, the mass that~$\bar\nu_\beta'$ puts on~$\{\zeta=0\}$ never enters and so we will discard it by putting $\bar\nu_\beta(\cdot):=\bar\nu_\beta'(\cdot\cap\{\zeta\ne0\})$. Along the  subsequence~$n_j\to\infty$ such that \eqref{E:1.7} holds, $E\bigl(\texte^{-Y_{n_j}}\bigr)$ thus converges to
\begin{equation}
\label{E:7.55}
E\,\biggl(\exp\Bigl\{-\tilde c\, \texte^{\alpha\sqrt\beta\, s}\int
Z(\textd x)\otimes \bar\nu_\beta(\textd\zeta)\bigl(1-\texte^{-\int f(x,\cdot)\textd\zeta})\Bigr\}\biggr)
\end{equation}
 with the integral inside the exponential finite for almost every~$Z$. 
  The measure~$\bar\nu_\beta$ is non-zero due to the fact that, for~$f(x,h)\ge1_{[-u,u]}(h)$, the associated~$Y_n$ will have a non-trivial limit for~$u$ sufficiently large, thanks to the tightness of the maximum in Theorem~\ref{thm-5.1}.
 
We will now rewrite \eqref{E:7.55} back to the form with an infinite sum. First observe that,
for any~$\epsilon>0$ and~$f(x,h)\ge1_{[r,u]}(h)$ with~$\supp(f)\subseteq[0,1]\times (r-1,\infty)$, the left-hand side of \eqref{E:7.51u} is then at least $P_{n,\beta}(\max_{x\in\Lambda_n}\varphi^{\DG}_x< \lfloor m_n\rfloor +r)$ which by Theorem~\ref{thm-5.1} tends to~$1$ as $n\to\infty$ and~$r\to\infty$. For~$\epsilon:=2$, the reasoning underlying \eqref{E:7.51u} then gives
\begin{equation}
\bar\nu_{\beta} \Bigl(\zeta\bigl([r,\infty)\bigr)\ge1\Bigr)\,\,\underset{r\to\infty}\longrightarrow0,
\end{equation}
thus showing that~$\bar\nu_\beta$ is concentrated on point measures that charge at least one point but whose support is upper bounded. In particular, $\max\zeta:=\sup\{n\in\Z\colon\zeta(\{n\})>0\}$ is well defined~$\bar\nu_\beta$-a.s. Noting that \eqref{E:7.51b} gives
 \begin{equation}
\label{E:7.57}
\bar\nu_\beta\circ\tau_n^{-1} = \texte^{-\alpha\sqrt\beta\,n}\,\bar\nu_\beta
\end{equation}
for all~$n\in\Z$, we get $\bar\nu_\beta(\max\zeta =0)>0$ by the fact that~$
\max\zeta$ is well defined~$\bar\nu_\beta$-a.s. This permits us to define the probability measure
\begin{equation}
\nu_\beta(\cdot):=\frac{\bar\nu_\beta(\cdot\cap\{\max\zeta =0\})}{\bar\nu_\beta(\max\zeta =0)}
\end{equation}
which by \eqref{E:7.57} brings~$\bar\nu_\beta$ to the form
\begin{equation}
\bar\nu_\beta(\cdot) = \bar\nu_\beta(\max\zeta =0)\sum_{n\in\Z}\texte^{\alpha\sqrt\beta\, n}\,\nu_\beta\circ\tau_n^{-1}(\cdot).
\end{equation}
Using this we rewrite \eqref{E:7.55} as
\begin{equation}
E\,\biggl(\exp\Bigl\{-\tilde c_\beta\, \texte^{\alpha\sqrt\beta\, s}\int
Z(\textd x)\otimes\nu_\beta(\textd\zeta)\sum_{n\in\Z}\texte^{-\alpha\sqrt\beta\, n}\bigl(1-\texte^{-\int f(x,n+\cdot)\textd\zeta})\Bigr\}\biggr),
\end{equation}
where $\tilde c_\beta:= \bar\nu_\beta(\max\zeta =0)\tilde c$.
 Expressing the result  in process language, this now gives the desired claim \eqref{E:1.9}.
\end{proofsect}

It remains to give:

\begin{proofsect}{Proof of Theorem~\ref{thm-1.1}}
This is a relatively straightforward consequence of Theorem~\ref{thm-1.2} and tightness of the maximum. Indeed, for any integers~$u<t$ 
\begin{equation}
\label{E:6.49}
\begin{aligned}
0\le E_{n,\beta}\Bigl(\texte^{-\lambda\sum_{x\in\Lambda_n}1_{(u,t]}(\varphi_x^{\DG}- \lfloor m_n\rfloor)}&\Bigr)-P_{n,\beta}\biggl(\,\max_{x\in\Lambda_n}\varphi_x^{\DG}\le \lfloor m_n\rfloor+u\Bigr)
\\
&\le\texte^{-\lambda}+P_{n,\beta}\biggl(\,\max_{x\in\Lambda_n}\varphi_x^{\DG}> \lfloor m_n\rfloor+t\Bigr)
\end{aligned}
\end{equation}
which uses the observation that the sum is at least one when the maximum lies in the interval $\lfloor m_n\rfloor+(u,t]$. 
Since~$\lambda 1_{(u,t]}$ coincides with a bounded continuous function on the integers, along the sequence $\{n_j\}_{j\ge1}$ such that \eqref{E:1.7} holds, Theorem~\ref{thm-1.2} shows that the expectation in \eqref{E:6.49} tends to
\begin{equation}
E\biggl(\exp\Bigl\{-\tilde c_\beta\,\texte^{\alpha\sqrt\beta\,s}Z\bigl([0,1]\bigr)\int \nu_\beta(\textd \zeta)\sum_{n\in\Z}\texte^{-\alpha\sqrt\beta\, n}\bigl(1-\texte^{-\lambda\int 1_{(u,t]}(n+\cdot)\textd\zeta})\Bigr\}\biggr).
\end{equation}
Taking~$t\to\infty$ followed by $\lambda\to\infty$, this becomes
\begin{equation}
\label{E:7.56a}
E\biggl(\exp\Bigl\{-\tilde c_\beta\,\texte^{\alpha\sqrt\beta\,s}Z\bigl([0,1]\bigr)\int \nu_\beta(\textd \zeta)\sum_{n\in \Z}\texte^{-\alpha\sqrt\beta\, n}1_{\{\zeta((-n+u,\infty))\ge1\}}\Bigr\}\biggr).
\end{equation}
Thanks to Theorem~\ref{thm-5.1}, all the errors in \eqref{E:6.49} are wiped out in these limits, so this expectation is also the limit of the CDF of the centered maximum. 
 Relabeing $-n+u$ for~$n$ and invoking  the notations \twoeqref{E:1.13u}{E:1.14u}, we  get the claim.
\end{proofsect}

For completeness' sake, we also give:

\begin{proofsect}{Proof of Corollary~\ref{cor-1.2}}
We will adhere to the notation in the statement in which~$\varphi'$ denotes the GFF and~$\varphi$ the DG model. As will be discussed in the proof of Theorem~\ref{thm-6.1}, a routine calculation along with \eqref{E:1.13u} shows
\begin{equation}
\label{E:7.53a}
P_{n,\beta}'\Bigl(\,\max_{x\in\Lambda_{n}}\varphi_x'\le  m_n + u\Bigr)\,\,\underset{n\to\infty}\longrightarrow\,\, E\bigl(\texte^{-(\alpha\sqrt\beta)^{-1}\ZZ\texte^{-\alpha\sqrt\beta\, u}}\bigr),\quad u\in\R,
\end{equation}
where~$\ZZ:=Z([0,1])$.
Given~$r\in\R$ and~$u\in\R$,  \eqref{E:1.6} then gives
\begin{equation}
\begin{aligned}
\liminf_{n\to\infty}\,
\biggl[P_{n,\beta}\Bigl(\,\max_{x\in\Lambda_{n}}&\,\varphi_x\le \lfloor m_n \rfloor + u\Bigr)
- P_{n,\beta}'\Bigl(\,\max_{x\in\Lambda_{n}}\varphi_x'\le m_n + u+r\Bigr)\biggr]
\\
&\ge E\Bigl(\texte^{-\hat c_\beta(s)\ZZ\texte^{-\alpha\sqrt\beta\,\lfloor u\rfloor}}\Bigr)-E\Bigl(\texte^{-(\alpha\sqrt\beta)^{-1}\ZZ\texte^{-\alpha\sqrt\beta\,(u+r)}}\Bigr),
\end{aligned}
\end{equation}
where~$s$ is a value obtained by taking a sequence~$n_k\to\infty$ achieving the \textit{limes inferior} such that $s:=\lim_{k\to\infty}(m_{n_k}-\lfloor m_{n_k}\rfloor)$ exists. Thanks to the explicit form of~$\hat c_\beta(s)$, the right-hand side is non-negative regardless of~$s$ once~$r$ obeys $(\alpha\sqrt\beta)^{-1}\texte^{-\alpha\sqrt\beta (r+1)}\ge \hat c_\beta(0)\texte^{\alpha\sqrt\beta}$. Solving for $r$, we get $r \le -2 -(\alpha\sqrt\beta)^{-1}\log[\hat c_\beta(0)\alpha\sqrt\beta]$.

A completely analogous argument shows that
\begin{equation}
\begin{aligned}
\limsup_{n\to\infty}\,
\biggl[P_{n,\beta}\Bigl(\,\max_{x\in\Lambda_{n}}&\,\varphi_x\le \lfloor m_n \rfloor + u\Bigr)
- P_{n,\beta}'\Bigl(\,\max_{x\in\Lambda_{n}}\varphi_x'\le m_n + u+r\Bigr)\biggr]
\\
&\le E\Bigl(\texte^{-\hat c_\beta(0)\ZZ\texte^{-\alpha\sqrt\beta\,\lceil u\rceil}}\Bigr)-E\Bigl(\texte^{-(\alpha\sqrt\beta)^{-1}\ZZ\texte^{-\alpha\sqrt\beta\,(u+r)}}\Bigr)\le0
\end{aligned}
\end{equation}
once~$\hat c_\beta(0)\ge(\alpha\sqrt\beta)^{-1}\texte^{-\alpha\sqrt\beta r}$. Solving for $r$, we get $r \ge -(\alpha\sqrt\beta)^{-1}\log[\hat c_\beta(0)\alpha\sqrt\beta]$.  Replacing~$\lfloor m_n\rfloor+u$ by~$u$ and using that $|m_n-\lfloor m_n\rfloor|\le1$  along with the monotonicity of the involved CDFs, the claim follows with~$a:=1+(\alpha\sqrt\beta)^{-1}|\log[\hat c_\beta(0)\alpha\sqrt\beta]|$.
\end{proofsect}

\subsection{Proof of Theorem~\ref{thm-6.1}}
\label{sec-7.4}\noindent
As our final task, we  will  provide the details how to add the spatial component to convergence of the extremal process of the Branching Random Walk with normal step distribution. The argument is based on manipulations that are likely standard in the theory of Mandelbrot's multiplicative cascades.
To reduce clutter, we will assume that~$\beta=1$ throughout. (The scaling by~$\beta$ can be added to the statement at the very end.) Departing from our earlier notation, we will also write~$\varphi^{(n)}$ for the GFF on~$\Lambda_n$.

Let us start by recalling Madaule's result~\cite{Madaule2}  that we cast in the form of convergence of Laplace transforms: There exists an a.s.-positive and finite random variable~$\ZZZ$ and a law~$\nu$ on $\MM_\N(\R)$ with
\begin{equation}
\label{E:7.58}
\nu\Bigl(\chi\in\MM_\N(\R)\colon \chi\bigl((-\infty,0]\bigr)=\infty,\, \chi\bigl(\{0\}\bigr)\ge1,\,\chi\bigl((0,\infty)\bigr)=0 \Bigr)=1
\end{equation}
such that for all continuous $f\colon \R\to[0,\infty)$ with compact support,
\begin{equation}
\label{E:7.59}
\begin{aligned}
E_{n,1}'\biggl(\exp\Bigl\{-&\sum_{x\in\Lambda_n}f\bigl(\varphi_x^{(n)}-m_n\bigr)\Bigr\}\biggr)
\\
&\underset{n\to\infty}\longrightarrow\,\,
E\biggl(\exp\Bigl\{- \ZZZ \int \texte^{-\alpha h}\textd h\otimes\nu(\textd\chi)(1-\texte^{-\int f(h+\cdot)\textd\chi}\bigr)\Bigr\}\biggr),
\end{aligned}
\end{equation}
where $E_{n,1}'$ is expectation with respect to the law $P_{n,1}'$ of GFF on~$\Lambda_n$ with~$\beta:=1$ and $\alpha:=\sqrt{2\log b}$. The expectation on the right is with respect to the law of~$\ZZZ$. (We deliberately use a different letter than~$\ZZ$ as a connection with a random measure yet needs to be shown.)  We start with the observation that~$\ZZZ$ obeys a so called cascade relation:

\begin{lemma}[Cascade relation]
\label{lemma-7.6}
Given a natural~$n\ge1$, let~$\varphi^{(n)}$ be a GFF on~$\Lambda_n$ and let $\{\ZZZ_{n,x}\colon x\in\Lambda_n\}$ be i.i.d.\ copies of~$\ZZZ$ that we assume are independent of~$\varphi^{(n)}$. Then 
\begin{equation}
\label{E:7.60a}
\ZZZ \laweq\sum_{x\in\Lambda_n}\texte^{\alpha\varphi^{(n)}_x - \alpha^2 n}\,\ZZZ_{n,x}.
\end{equation}
\end{lemma}

\begin{proofsect}{Proof}
Using suitable approximations along with the fact that the integral on the right of \eqref{E:7.59} is insensitive to changes of~$f$ at a single point the limit \eqref{E:7.59} applies to the function~$f(h):=\lambda 1_{(u,\infty)}(h)$, for any~$u\in\R$. Since the difference of the left-hand side of \eqref{E:7.59} with~$f$ as above and its limit as $\lambda\to \infty$ is bounded by $\texte^{-\lambda}$ independent of $n$, we may take~$\lambda\to\infty$ and invoke the properties \eqref{E:7.58} to get
\begin{equation}
\label{E:7.60}
P_{n,1}'\Bigl(\,\max_{x\in\Lambda_n}\varphi_x^{(n)}\le m_n+u\Bigr)
\,\,\underset{n\to\infty}\longrightarrow\,\, E\bigl(\texte^{-\alpha^{-1}\texte^{-\alpha u}\ZZZ}\bigr)
\end{equation}
thus recovering the main result of A\"idekon~\cite{Aidekon} on which \cite{Madaule2} draws heavily.

Let~$M_n$ be a random variable with $M_n\laweq\max_{x\in\Lambda_n}\varphi_x^{(n)}$ and let us write $\zeta_1,\dots,\zeta_b$ for i.i.d. copies of~$\NN(0,1)$ and~$M_{n-1}^{(1)},\dots,M^{(b)}_{n-1}$ for i.i.d.\ copies of~$M_{n-1}$ that are assumed to be independent of~$\zeta_1,\dots,\zeta_b$. Using that the GFF is a Branching Random Walk, decomposing according to the first ``step'' shows
\begin{equation}
M_n \laweq \max_{i=1,\dots,b}\bigl(\zeta_i+M_{n-1}^{(i)}\bigr).
\end{equation}
This equates the probability on the left of \eqref{E:7.60} with~$b$-th power of the probability that~$\zeta+M_{n-1}\le m_n+u$, for~$\zeta=\NN(0,1)$ independent of~$M_{n-1}$. Passing this through the $n\to\infty$ limit while using that~$m_n = m_{n-1}+\alpha+o(1)$ gives 
\begin{equation}
E\bigl(\texte^{-\alpha^{-1}\texte^{-\alpha u}\ZZZ}\bigr) = \Bigl[E\bigl(\texte^{-\alpha^{-1}\texte^{-\alpha (u-\zeta+\alpha)}\ZZZ}\bigr)\Bigr]^b,
\end{equation}
where~$\zeta=\NN(0,1)$ is independent of~$\ZZZ$ on  the right-hand side.  Interpreting the right-hand side expectation with respect to a product measure, the fact that the Laplace transform determines the law yields
\begin{equation}
\label{E:7.63}
\ZZZ\laweq \sum_{i=1}^b\texte^{\alpha \zeta_i-\alpha^2}\ZZZ_i,
\end{equation}
where~$\zeta_1,\dots,\zeta_b$ are i.i.d.\ $\NN(0,1)$ and~$\ZZZ_1,\dots,\ZZZ_b$ are i.i.d. copies of~$\ZZZ$ that are independent of~$\zeta_1,\dots,\zeta_b$. This is the desired cascade relation \eqref{E:7.60a} for~$n=1$. The general case is now readily proved by induction.
\end{proofsect}

The argument underlying the proof of the previous lemma can be bolstered to construct a full multiplicative cascade associated with random variable~$\ZZZ$:

\begin{lemma}[Multiplicative cascade]
Let~$\Lambda_0$ be a singleton designating an empty sequence. There exists random variables $\{\ZZZ_{n,x},\zeta_{n,x}\colon n\ge0,\,x\in\Lambda_n\}$ such that
\begin{enumerate}
\item[(1)] $\{\zeta_{n,x}\colon n\ge0,\,x\in\Lambda_n\}$ are i.i.d.~$\NN(0,1)$,
\item[(2)] for each~$n\ge0$, $\{\ZZZ_{n,x}\colon x\in\Lambda_n\}$ are i.i.d.\ copies of~$\ZZZ$ that are independent of  the  random variables $\{\zeta_{k,x}\colon k=0,\dots,n,\,x\in\Lambda_k\}$,
\item[(3)] for all~$n\ge0$ and all~$x\in\Lambda_n$,
\begin{equation}
\label{E:7.66}
\ZZZ_{n,x} = \sum_{\begin{subarray}{c}
z\in\Lambda_{n+1}\\ m(z)=x
\end{subarray}}
\texte^{\alpha\zeta_{n+1,z}-\alpha^2 }\ZZZ_{n+1,z}
\end{equation}
holds pointwise almost surely.
\end{enumerate}
\end{lemma}

\begin{proofsect}{Proof}
The one-step cascade relation \eqref{E:7.63} can be realized as an almost sure identity by simply declaring~$\ZZZ$ to be the right-hand side of \eqref{E:7.63}. Given~$n\ge0$ and i.i.d.\ standard normals $\{\zeta_{k,x}\colon  k = 0,\dots,n,  \,\,x\in\Lambda_k\}$, we now construct the joint law of random variables $\{\ZZZ_{k,x},\zeta_{n,x}\colon k=0,\dots,n,\,x\in\Lambda_k\}$ by taking independent i.i.d.\ copies $\{\ZZZ_{n,x}\colon x\in\Lambda_n\}$ of~$\ZZZ$ and using \eqref{E:7.66} iteratively to define $\{\ZZZ_{k,x}\colon x\in\Lambda_k\}$ for all~$k<n$. Noting that the cascade relation from Lemma~\ref{lemma-7.6} makes these laws consistent, the Kolmogorov Extension Theorem extends this to a law on the full infinite collection.
\end{proofsect}

Using the family $\{\ZZZ_{n,x},\zeta_{n,x}\colon n\ge0,\,x\in\Lambda_n\}$, we can couple all of the GFFs on the same probability space by setting
\begin{equation}
\label{E:7.67}
\varphi^{(n)}_x:= \sum_{k=0}^n  \zeta_{k,m^{n-k}(x)},\quad x\in\Lambda_n.
\end{equation}
We now represent the random variables using a random Borel measure:

\begin{lemma}
\label{lemma-7.8}
There exists a unique random Borel measure~$Z$ on~$\R$ concentrated on~$[0,1)$ such that for each~$g\colon[0,1]\to[0,\infty)$ and each~$n\ge0$,
\begin{equation}
\int Z(\textd x) g\bigl( \lfloor x\rfloor_n\bigr) = \sum_{x\in\Lambda_n}g\bigl([x]_n\bigr)\texte^{\alpha\varphi_x^{(n)}-\alpha^2 n}\ZZZ_{n,x}
\end{equation}
holds almost surely.  Here $\lfloor x\rfloor_n:=b^{-n}\lfloor b^n x\rfloor$.
\end{lemma}

\begin{proofsect}{Proof}
The identity \eqref{E:7.66} upgrades the equality in law \eqref{E:7.60a} to the pointwise identity
\begin{equation}
\label{E:7.68}
\texte^{\alpha\varphi_x^{(n)}-\alpha^2 n}\ZZZ_{n,x} = \sum_{\begin{subarray}{c}
z\in\Lambda_{n+k}\\ m^k(z)=x
\end{subarray}}\texte^{\alpha\varphi_z^{(n+k)}-\alpha^2 (n+k)}\ZZZ_{n+k,z},\quad x\in\Lambda_n,
\end{equation}
for each~$n\ge0$. Let~$\AA$ be the set of finite disjoint unions of sets from the collection
\begin{equation}
\AA_0:=\bigl\{[kb^{-n},(k+1)b^{-n})\colon k=0,\dots,b^{n}-1,\, n\ge0\bigr\}.
\end{equation}
 Define~$Z\colon\AA_0\to[0,\infty)$ by putting
\begin{equation}
\label{E:7.70}
Z\bigl([kb^{-n},(k+1)b^{-n})\bigr):=\texte^{\alpha\varphi_x^{(n)}-\alpha^2 n}\ZZZ_{n,x}
\end{equation}
for the unique~$x\in\Lambda_n$ such that~$[x]_n\in [kb^{-n},(k+1)b^{-n})$. Since \eqref{E:7.68} makes~$Z$ finitely additive on~$\AA_0$, a standard argument extends~$Z$ uniquely to an additive set function on~$\AA$. But~$\AA$ is an algebra such that any decreasing sequence $\{A_n\}_{n\ge1}$ of elements from~$\AA$ with~$A_n\downarrow\emptyset$ is eventually empty. This shows that~$Z$ is also trivially countably subadditive on~$\AA$ and so, by the Carath\'eodory Extension Theorem, extends uniquely to a Borel measure on~$[0,1)$. Setting the measure to zero in the complement of~$[0,1)$, the claim follows directly from \eqref{E:7.70}.
\end{proofsect}

With these in hand, we are ready to give:

\begin{proofsect}{Proof of Theorem~\ref{thm-6.1}}
Let $\lambda>0$ and let~$f\colon[0,1]\times\R\to[0,\infty)$ be continuous with support in~$[0,1]\times[-\lambda,\lambda]$. For~$k\ge0$ define $f_k\colon [0,1] \times \mathbb R\to[0,\infty)$ by
\begin{equation}
f_k(x,h):=f\biggl(\,\sum_{j=1}^k x_j b^{-j},h\biggr). 
\end{equation}
where~$x_j:=\lfloor b^j  x\rfloor\text{ mod }b$ for each~$j=1,\dots, k$. Note that since~$f$ is uniformly continuous and $m_n- (m_{n-k}+\alpha k)\to0$ as~$n\to\infty$, the quantity 
\begin{equation}
\epsilon_k:=\sup_{n\ge 2k}\,\sup_{x\in[0,1]}\,\sup_{h\in\R}\,\Bigl|\,f(x,h) - f_k (x,h+m_n-m_{n-k}-\alpha k)\Bigr|
\end{equation}
obeys~$\epsilon_k\to0$ as~$k\to\infty$. Using that also $\supp(f_k)\subseteq [0,1]\times[-\lambda,\lambda]$, we have
\begin{equation}
\label{E:7.73}
\begin{aligned}
\biggl|\,\sum_{x\in\Lambda_n}f\bigl([x]_n,\varphi^{(n)}_x-m_n\bigr)
-\sum_{x\in\Lambda_n}f_k\bigl([x]_n, &\varphi^{(n)}_x-m_{n-k}-\alpha k\bigr)\biggr|
\\
&\le \epsilon_k\sum_{x\in\Lambda_n}1_{[-\lambda-1,\infty)}\bigl(\varphi^{(n)}_x-m_n\bigr)
\end{aligned}
\end{equation}
once~$n\ge 2k$ and~$k$ is so large that $|m_n- m_{n-k}-\alpha k|\le1$. Since the convergence \eqref{E:7.59} entails that the level sets of the GFF are tight, we can swap the first sum on the left for the second sum provided we take the limit~$k\to\infty$ after taking~$n\to\infty$.

Invoking the coupling \eqref{E:7.67} and using that $f_k([x]_n,\cdot)=f([x]_k,\cdot)$ when $n\ge k$, for~$n$ in excess of~$k$ we get
\begin{equation}
\begin{aligned}
\sum_{x\in\Lambda_n}f_k\bigl([x]_n,&\varphi^{(n)}_x-m_{n-k}-\alpha k\bigr)
\\
&=\sum_{x\in\Lambda_k}\sum_{\begin{subarray}{c}
z\in\Lambda_n\\ m^{n-k}(z)=x
\end{subarray}}
f\bigl([x]_k,\varphi^{(n)}_z-\varphi^{(k)}_x-m_{n-k}+\varphi^{(k)}_x-\alpha k\bigr).
\end{aligned}
\end{equation}
Noting that the sets of values $\{\varphi^{(n)}_z-\varphi^{(k)}_x\colon z\in\Lambda_n,\,m^{n-k}(z)=x\}$ are independent for distinct~$x\in\Lambda_k$, and are all independent of~$\varphi^{(k)}$ with law of~$\{\varphi^{(n-k)}_z\colon z\in\Lambda_{n-k}\}$, Madule's limit result \eqref{E:7.59} along with the Bounded Convergence Theorem show
\begin{equation}
\begin{aligned}
&\lim_{n\to\infty}\,E_{n,1}'\biggl(\exp\Bigl\{-\sum_{x\in\Lambda_n}f_k\bigl([x]_n,\varphi_x^{(n)}-m_{n-k} - \alpha k\bigr)\Bigr\}\biggr) 
\\
&\,\,=
E_{k,1}'\otimes E\biggl(\exp\Bigl\{-  \sum_{x\in\Lambda_k}\ZZZ_{k,x}\int \texte^{-\alpha h}\textd h\otimes\nu(\textd\chi)(1-\texte^{-\int f([x]_k,h+\varphi^{(k)}_x-\alpha k+\cdot)\textd\chi}\bigr)\Bigr\}\biggr),
\end{aligned}
\end{equation}
where the product expectation indicates that the law of~$\varphi^{(k)}$ is independent of the independent copies $\{\ZZZ_{k,x}\colon x\in\Lambda_k\}$ of~$\ZZZ$. 
Performing a routine change of variables, the latter expectation equals
\begin{equation}
\begin{aligned}
E_{k,1}'\otimes E\biggl(&\exp\Bigl\{-  \int \texte^{-\alpha h}\textd h\otimes\nu(\textd\chi)\sum_{x\in\Lambda_k}\texte^{\alpha \varphi^{(k)}_x-\alpha^2 k}\ZZZ_{k,x}(1-\texte^{-\int f([x]_k,h+\cdot)\textd\chi}\bigr)\Bigr\}\biggr)
\\
&=
E\biggl(\exp\Bigl\{-  \int Z(\textd x)\otimes\texte^{-\alpha h}\textd h\otimes\nu(\textd\chi)(1-\texte^{-\int f_k(x,h+\cdot)\textd\chi}\bigr)\Bigr\}\biggr),
\end{aligned}
\end{equation}
where Lemma~\ref{lemma-7.8} was used to rewrite the expression inside exponential as an integral with respect to~$Z$-measure. Since~$f_k\to f$ pointwise and the indicator $1_{[-\lambda,\infty)}(h)$ can be added to the integral over~$h$ thanks to the restriction on~$\supp(f)$, the Bounded Convergence Theorem shows that the expectation tends to that for~$f$ as~$k\to\infty$. In combination with \eqref{E:7.73}, this proves the desired claim.
\end{proofsect}


\vbox{
\section*{Acknowledgments}
\nopagebreak\nopagebreak\noindent
This project has been supported in part by the NSF award DMS-1954343 and the BSF award 2018330.  The authors are  grateful to anonymous referees for detailed comments on a first version of this manuscript. 

}

\section*{Declarations}
\nopagebreak\nopagebreak\noindent

\noindent
\textbf{Data availability statement: } No data is associated with this research.

\medskip
\noindent
\textbf{Conflict of interest: }The authors declare no conflicts of interest related to this work.

\bibliographystyle{abbrv}

\end{document}